%% file: stats_sdp.tex
\documentclass[11pt, a4paper]{article}

\input{packages_notes}

\input{commandes_guillaume}

\newcommand{\Aobs}{A}
\newcommand{\proj}[1]{\cP(#1)}
\newcommand{\porth}[1]{\cP^\perp(#1)}
\newcommand{\floorsup}[1]{\lceil #1 \rceil}
\newcommand{\tr}[1]{\mathrm{Tr}(#1)}
\newcommand{\normop}[1]{\norm{#1}_{\mbox{op}}}

\usepackage{marginnote}
\newcommand{\bluenote}[1]{ {\color{blue} GL: #1}  \marginnote{({\color{blue} G-note})}}

\newcommand{\rednote}[1]{  {\color{red}MS: #1} \marginnote{({\color{red} MS-note})}}

\usepackage{mathrsfs}
\newtheorem{lemma}{Lemma}
\newtheorem{cor}{Corollary}

\usepackage{comment}
\usepackage{subcaption}
\usepackage{caption}
\usepackage[skip=-5pt, font={small}]{subcaption}
\usepackage{algorithmic}
\usepackage{mathabx}

\usepackage{multirow,bigdelim}
\usepackage{array}     % needed for the nice table of figures with rownames and column names 
\newcommand{\mb}[1]{\mbox{\boldmath$#1$}}
\newcommand{\bs}[1]{\mathbf{#1}}

\newcommand\mydef{\mathrel{\overset{\makebox[0pt]{\mbox{\normalfont\tiny\sffamily def}}}{=}}}

\usepackage{authblk}

\numberwithin{equation}{section}
\begin{document}
\title{Learning with Semi-Definite Programming: new statistical bounds based on fixed point analysis and excess risk curvature}
\author{St{\'e}phane Chr{\'e}tien, Mihai Cucuringu,  Guillaume Lecu{\'e} and Lucie Neirac}
\author[1]{St{\'e}phane Chr{\'e}tien, Mihai Cucuringu, Guillaume Lecu{\'e}, Lucie Neirac \\ email:   \href{mailto:stephane.chretien@npl.co.uk}{stephane.chretien@npl.co.uk},
\href{mailto:mihai.cucuringu@stats.ox.ac.uk}{mihai.cucuringu@stats.ox.ac.uk}, \href{mailto:lecueguillaume@gmail.com}{guillaume.lecue@ensae.fr}, \href{mailto:Lucie.NEIRAC@ensae.fr}{Lucie.neirac@ensae.fr},\\  Universit{\'e} Lyon 2, Oxford University, The Alan Turing Institute, CREST, ENSAE, IPParis.}

\date{}                     %% if you don't need date to appear
\setcounter{Maxaffil}{0}
\renewcommand\Affilfont{\itshape\small}

\maketitle

\begin{abstract}
Many statistical learning problems have recently been shown to be amenable to Semi-Definite Programming (SDP), with \textit{community detection} and \textit{clustering in Gaussian mixture models} as the most striking instances \cite{javanmard2016phase}. Given the growing range of applications of SDP-based techniques to machine learning problems, and the rapid progress in the design of efficient algorithms for solving SDPs, an intriguing  question is to understand how the recent advances from empirical process theory can be put to work in order to provide a precise statistical analysis of SDP estimators. 
    
    In the present paper, we borrow cutting edge techniques and concepts from the learning theory literature, such as fixed point equations and excess risk curvature arguments, which yield general estimation and prediction results for a wide class of SDP estimators. From this perspective, we revisit some classical results in community detection from \cite{guedon2016community} and \cite{chen2016statistical}, and we obtain statistical guarantees for SDP estimators used in signed clustering, group synchronization and MAXCUT.
\end{abstract}

{
  \hypersetup{linkcolor=black}
  \tableofcontents
}
% \tableofcontents

	\section{Introduction} % (fold)
	\label{sec:introduction}
	Many statistical learning problems have recently been shown to be amenable to Semi-Definite Programming (SDP), with \textit{community detection} and \textit{clustering in Gaussian mixture models} as the most striking instances where SDP performs significantly better than other current approaches \cite{javanmard2016phase}. SDP is a class of convex optimisation problems generalising linear programming to linear problems over semi-definite matrices \cite{todd2001semidefinite}, \cite{wolkowicz2012handbook}, \cite{boyd2004convex}, and which proved an important tool in the computational approach to difficult challenges in automatic control, combinatorial optimisation, polynomial optimisation, data mining, high dimensional statistics and the numerical solution to partial differential equations. The goal of the present paper is to introduce a new fixed point approach to the statistical analysis of SDP-based estimators and illustrate our method on four current problems of interest, namely  \textit{MAX-CUT}, \textit{community detection},  \textit{signed clustering},  \textit{angular group synchronization}. The rest of this section gives some historical background and presents the mathematical definition of SDP based estimators. 
	
	\subsection{Historical background}
	SDP is a class of optimisation problems which includes linear programming as a particular case and can be written as the set of problems over symmetric (resp. Hermitian) positive semi-definite matrix variables, with linear cost function and affine constraints, i.e. optimization problems of the form
	\begin{align}
	\max_{Z \succeq 0} \left( \inr{A,Z}:  \inr{B_j,Z} =b_j \mbox{ for } j=1,\ldots,m\right)
	\end{align}where $A, B_1, \ldots, B_m$ are given matrices. SDPs are convex programming problems which can be solved in polynomial time when the constraint set is compact and it plays a paramount role in a large number of convex and nonconvex problems, for which it often appear as a convex relaxation \cite{anjos2011handbook}. 
	
	We will sometimes use the notation $\mathbb S_{n,+}$ (resp. $\mathbb S_{n,-}$) for the cone of positive semi-definite matrices (resp. negative semi-definite matrices).
	
	\subsubsection{Early history}
	Early use of Semi-Definite programming to statistics can be traced back to  \cite{scobey1978vector} and \cite{fletcher1981nonlinear}. In the same year, Shapiro used SDP in factor analysis \cite{shapiro1982weighted}. The study of the mathematical properties of SDP then gained momentum with the introduction of Linear Matrix Inequalities (LMI) and their numerous applications in control theory, system identification and signal processing. The book \cite{boyd1994linear} is the standard reference of these type of results, mostly obtained in the 90's. 
	
	\subsubsection{The Goemans-Williamson SDP relaxation of \textit{Max-Cut} and its legacy}
	A notable turning point is the publication of \cite{goemans1995improved} where SDP was shown to provide a 0.87 approximation to the NP-Hard problem known as \textit{Max-Cut}. The \textit{Max-Cut} problem is a clustering problem on graphs which consists in finding two complementary subsets $S$ and $S^c$ of nodes such that the sum of the weights of the edges between $S$ and $S^c$ is maximal. In \cite{goemans1995improved}, the authors approach this difficult combinatorial problem by using what is now known as the Goemans-Williamson \textit{SDP relaxation} and use the Choleski factorization of the optimal solution to this SDP in order to produce a randomized scheme achieving the .87 bound in expectation. Moreover, this problem can be seen as a first instance where the Laplacian of a graph is employed in order to provide an optimal bi-clustering in a graph and certainly represents the first chapter of a long and fruitful relationship between clustering, embedding and Laplacians. Other SDP schemes for approximating hard combinatorial problems are, to name a few, for the graph coloring problem  \cite{karger1998approximate}, for satisfiability problem \cite{goemans1995improved,goemans1994new}. These results were later surveyed in \cite{lemarechal1995new,goemans1997semidefinite} and \cite{wolkowicz1999semidefinite}. The randomized scheme introduced by Goemans and Williamson was then further improved in order to study more general Quadratically Constrained Quadratic Programmes (QCQP) in various references, most notably \cite{nesterov1997semidefinite,zhang2000quadratic} and further extended in \cite{he2008semidefinite}. Many applications to signal processing are discussed in \cite{olsson2007solving}, \cite{ma2010semidefinite}; one specific reduced complexity implementation in the form of an eigenvalue minimisation problem and its application to binary least-squares recovery and denoising is presented in \cite{chretien2009using}. 
	
	\subsubsection{Relaxation of machine learning and high dimensional statistical estimation problems}
	Applications of SDP to problems related with machine learning is more recent and probably started with the SDP relaxation of $K$-means in \cite{peng2005new,peng2007approximating} and later in \cite{ames2014guaranteed}. This approach was then further improved using a refined statistical analysis by \cite{royer2017adaptive} and \cite{giraud2018partial}. Similar methods have also been applied to community detection \cite{hajek2016achieving,abbe2015exact} and for the weak recovery viewpoint,
	\cite{guedon2016community}. This last approach was also  re-used via the kernel trick for the point cloud clustering \cite{chretien2016semi}. Another incarnation of SDP in machine learning is the extensive use of nuclear norm-penalized least-square costs as a surrogate for rank-penalization in low-rank recovery problems such as matrix completion in recommender systems, matrix compressed sensing, natural language processing and quantum state tomography; these topics are surveyed in  \cite{davenport2016overview}. 
	
	The problem of manifold learning was also addressed using SDP and is often mentioned as one of the most accurate approaches to the problem, let aside its computational complexity; see \cite{weinberger2005nonlinear,weinberger2006unsupervised,weinberger2006introduction,hegde2012near}. Connections with the design of fast converging  Markov-Chains were also exhibited in \cite{sun2006fastest}.
	
In a different direction, A. Singer and collaborators have recently promoted the use of SDP relaxation for estimation under group invariance, an active area with many applications \cite{singer2011angular,bandeira2014multireference}. SDP-based relaxations have also been considered in \cite{syncZ2} in the context of synchronization over $\mathbb{Z}_2$ in signed multiplex networks with constraints, and \cite{syncRank} in the setting of ranking from inconsistent and incomplete pairwise comparisons where an SDP-based relaxation of angular  synchronization over SO(2) outperformed a suite of state-of-the-art algorithms from the literature. Phase recovery using SDP was studied in e.g. \cite{waldspurger2015phase} and \cite{demanet2014stable}.  An extension to multi-partite clustering based on SDP was then proposed in \cite{karger1998approximate}. Other important applications of SDP are, information theory \cite{lovasz1979shannon}, estimation in power networks \cite{lavaei2011zero}, quantum tomography \cite{mazziotti2011large}, \cite{gross2010quantum} and  polynomial optimization via Sums-of-squares relaxations \cite{lasserre2015introduction, blekherman2012semidefinite}. Sums of squares relaxations were recently applied to statistical problems in \cite{de2019approximate,Hopkins2018SOS,yohann2017SOS}. Extension to the field of complex numbers, with $\inr{\cdot,\cdot}$ denoting the Hermitian inner product, has been less extensively studied but has many interesting applications and comes with efficient algorithms \cite{goemans1995improved,gilbert2017plea}.

    \subsection{Mathematical formulation of the problem}
	The general problem we want to study can be stated as follows. Let $A$ be a random matrix in $\R^{n\times n}$  and $\cC\subset \R^{n\times n}$ be a constraint. The object that we want to recover, for instance, the community membership vector in community detection, is related to an \textit{oracle} defined as
	\begin{equation}\label{eq:oracle_0}
	Z^*\in\argmax_{Z\in\cC} \;    \inr{\E A, Z}
	\end{equation}where $\inr{A,B}={\rm Tr}(A\bar B^\top) = \sum A_{ij}\bar B_{ij}$ when $A,B\in\bC^{n\times n}$ where $\bar z$ is the conjugate of $z\in\bC$.
	We would like to estimate $Z^*$, from which we can ultimately retrieve the object that really matters to us (for instance, by considering a singular vector associated to the largest singular value of $Z^*$). To that end, we consider the following estimator
	\begin{equation}\label{eq:estimator_0}
	\hat Z \in\argmax_{Z\in\cC} \;    \inr{A, Z}, 
	\end{equation} 
	which is simply obtained by replacing the unobserved quantity $\E A$ by the observation $A$.
	
	As pointed out, in many situations, $Z^*$ is not the object we want to estimate, but there is a straightforward relation between $Z^*$ and this object. For instance, consider the community detection problem, where the goal is to recover the class community vector $x^*\in\{-1, 1\}^n$ of $n$ nodes. Here, when $\cC$ is well chosen,  there is a close relation between $Z^*$ and $x^*$, given by $Z^* = x^* (x^*)^\top$. We therefore need a final step to estimate $x^*$ using $\hat Z$, for instance, by letting $\hat x$ denote a top eigenvector of $\hat Z$, and then using the Davis-Kahan "sin-theta" Theorem \cite{daviskahan,dkuseful} to control the estimation of $x^*$ by $\hat x$ from the one of $Z^*$ by $ \hat Z $.
	
	When the constraint $\cC$ is of the form $\cC=\{Z\in\bR^{n\times n}: Z\succeq0, \inr{Z, B_j}= b_j, j=1,\ldots,m\}$ where $B_1, \ldots, B_m\in\bR^{n\times n}$ and $Z\succeq 0$ is notation for `` $Z$ is positive semidefinite'' then \eqref{eq:estimator_0} is a semidefinite programming (SDP) \cite{boyd2004convex}.
	%for which polynomial-time algorithm exist since \cite{MR1315703,MR1258086}. % attention was focused on SDP after the work of \cite{Goemans} who introduced a SDP to construct approximation algorithms for several NP-hard problems such as MAX-CUT, MAX 2-SAT and MAX 3-SAT. Following this approach several SDP were designed to solve or approximatively solve problems in graph (coloring, clique, cut, etc.). More importantly, a general concept of `` SDP relaxation'' became a classical tool in optimization, computer science, signal processing, discrete geometry, Banach spaces, etc. \cite{MR2121604,MR2181661,MR2298255,MR2350049,DBLP:journals/siamcomp/BansalDG19}. 

	\subsection{Goal of the paper}
	The aim of the present work is to present a general approach to the study of the statistical properties of SDP-based estimators defined in \eqref{eq:estimator_0}. In particular, using our framework, one is able  to obtain new (non-asymptotic) rates of convergence or exact reconstruction properties for a wide class of estimators obtained as a solution of a semidefinite program like \eqref{eq:estimator_0}. Specifically, our goal is to show that the solution to \eqref{eq:estimator_0} can be analyzed in a statistical way when $\E A$ is only partially and noisily observed in $A$. Even though the constraint $\cC$ may not necessarily be the intersection of the set of PSD (or Hermitian) matrices with linear spaces -- such as in the definition of SDP -- in the following, a solution $\hat Z$ of \eqref{eq:estimator_0} will be called a SDP estimator because in all our examples  $\hat Z$ will be solution of a SDP. But for the sake of generality we will only assume only minimal requirement on the shape of $\cC$. 	
	We also illustrate our results on a number of specific machine learning problems such as various forms of \textit{clustering problems} and \textit{angular group synchronization}.
	% in physics, people refer something else as "synchronization"
	Three out of the four examples worked out here are concerned with real-valued matrices. Only the angular synchronization problem is approached using complex matrices.

	\section{Main general results for the statistical analysis of SDP estimators} % (fold)
	\label{sec:main_general_results_for_the_statistics_of_sdp_estimators}
	From a statistical point of view, the task remains to estimate in the most efficient way the oracle $Z^*$, and to that end $\hat Z$ is our candidate estimator. The point of view we will use to evaluate how far $\hat Z$ is from $Z^*$ is coming from the learning theory literature. We therefore see $\hat Z$ as an empirical risk minimization (ERM) procedure built on a single observation $A$, where the loss function is the linear one $Z \in \cC \to \ell_Z(A) = -\inr{A,Z}$, and the oracle $Z^*$ is indeed the one minimizing the risk function $Z \in \cC \to \E \ell_Z(A)$ over $\cC$. Having this setup in mind, we can use all the machinery developed in learning theory (see for instance \cite{MR1641250,MR2329442,MR2319879,MR1739079}) to obtain rates of convergence for the ERM  (here $\hat Z$) toward the oracle (here $Z^*$).
	
	There is one key quantity driving the rate of convergence of the ERM: a fixed point complexity parameter. This type of parameter carries all the statistical complexity of the problem, and even though it is usually easy to set up, its computation can be tedious since it requires to control, with large probability, the supremum of empirical processes indexed by ``localized classes''. We now define this complexity fixed point related to the problem we are considering here.   
	
	\begin{Definition} 
		\label{def:fixed_point}
		Let $0<\Delta<1$. The fixed point complexity parameter at deviation $1-\Delta$ is 
		\begin{equation}\label{eq:fixed_point}
		r^*(\Delta) = \inf\left(r>0 : \bP\left[\sup_{Z\in\cC: \inr{\E A,Z^*-Z}\leq r}\inr{A-\E A, Z-Z^*}\leq (1/2) r\right]\geq 1-\Delta\right).
		\end{equation}
	\end{Definition}
	Fixed point complexity parameters have been extensively used in learning theory since the introduction of the localization argument \cite{MR2319879,MR2829871,MR1739079,MR1240719}. When they can be computed, they are preferred to the (global) analysis developed by Chervonenkis and Vapnik \cite{MR1641250} to study ERM, since the latter analysis always yields slower rates given that the Vapnik-Chervonenkis bound is a high probability bound on the non-localized empirical process $\sup_{Z\in\cC}\inr{A-\E A, Z-Z^*}$, which is an upper bound for $r^*(\Delta)$ since $\{Z\in\cC: \inr{\bE A, Z^*-Z}\leq r\}\subset \cC$. The gap between the two global and local analysis can be important since fast rates cannot be obtained using the VC approach, whereas the localization argument resulting in fixed points such as the one in Definition~\ref{def:fixed_point} may yield fast rates of convergence or even exact recovery results. 
	
	An example of a Vapnik-Chervonenkis's type of analysis of SDP estimators can be found in \cite{guedon2016community} for the community detection problem. An improvement of the latter approach has been obtained in \cite{MR3901009} thanks to a localization argument -- even though it is not stated in these words (we elaborate more on the two approaches from \cite{guedon2016community,MR3901009} in Section~\ref{sec:application_to_the_community_detection_problem}). Somehow, a fixed point such as  \eqref{eq:fixed_point} is a sharp way to measure the statistical performances of ERM estimators and in particular for the SDP estimators that we are considering here. They can even be proved to be optimal (in a minimax sense) when the noise $A-\E A$ is Gaussian \cite{LM13} and under mild conditions on the structure of $\cC$.

	Before stating our general result, we first recall a definition of a minimal structural assumption on the constraint $\cC$.
	
	\begin{Definition}\label{def:star_shapped}
		We say that the set $\cC$ is \textit{star-shaped in $Z^*$} when for all $Z\in\cC$, the segment $[Z, Z^*]$ is in $\cC$.
	\end{Definition}
	This is a pretty mild assumption satisfied, for instance, when $\cC$ is convex, which is the setup we will always encounter in practical applications, given that SDP estimators are usually introduced after a ``convex relaxation'' argument. Our main general statistical bound on SDP estimators is as follows.
	
	\begin{Theorem}\label{theo:main}
		We assume that the constraint $\cC$ is star-shaped in $Z^*$. Then, for all $0<\Delta<1$, with probability at least $1-\Delta$, it holds true that $\inr{\E A, Z^*-\hat Z}\leq r^*(\Delta)$.
	\end{Theorem}
	
	Theorem~\ref{theo:main} applies to any type of setup where an oracle $Z^*$ is estimated by an estimator $\hat Z$ such as in \eqref{eq:estimator_0}. Its result shows that $\hat Z$ is almost a maximizer of the true objective function $Z\to\inr{\bE A, Z}$ over $\cC$ up to $r^*(\Delta)$. In particular, when $r^*(\Delta)=0$, $\hat Z$ is exactly a maximizer such as $Z^*$ and, in that case, we can work with $\hat Z$ as if we were working with $Z^*$ without any loss.  
	
	Theorem~\ref{theo:main} may be applied in many different settings; in the following, we study four such instances. We will apply Theorem~\ref{theo:main} (or one of its corollary stated below) to some  popular problems in the networks and graph signal processing literatures, namely,   
	community detection \cite{fortunato2010community} (we will mostly revisit the results in \cite{guedon2016community} and \cite{MR3901009} from our perspective), 
	signed clustering \cite{SPONGE}, group synchronization \cite{sync} and MAX-CUT.

	The proof of Theorem~\ref{theo:main} is straightforward (mostly because the loss function is linear). Its importance stems from the fact that it puts forward two important concepts originally introduced in Learning Theory, namely that the complexity of the problem comes from the one of the local subset $\cC\cap\{Z:\inr{\E A, Z^*-Z}\leq r^*(\Delta)\}$ and that the ``radius'' $r^*(\Delta)$ of the localization is solution of a fixed point equation. For a setup given by a random matrix $A$ and a constraint $\cC$, we should try to understand how these two ideas apply  to obtain estimation properties of SDP estimators such as $\hat Z$. That is to understand the shape of the local subsets $\cC\cap\{Z:\inr{\E A, Z^*-Z}\leq r\}, r>0$ and the maximal oscillations of the empirical process $Z\to\inr{A-\E A, Z-Z^*}$ indexed by these local subsets.
	% We will do it in three cases. 
	We will consider this task in three distinct problem instances. 
	We now provide a proof for Theorem~\ref{theo:main}.
	
	\textbf{Proof of Theorem~\ref{theo:main}.} Denote $r^*=r^*(\Delta)$. Assume first that $r^*>0$ (the case $r^*=0$ is analyzed later). Let $\Omega^*$  be the event onto which for all $Z\in\cC$ if $\inr{\E A, Z^*-Z}\leq r^*$ then $\inr{A-\E A, Z-Z^*}\leq (1/2)r^*$. By Definition of $r^*$, we have $\bP[\Omega^*]\geq 1-\Delta$.

	Let $Z \in \cC$ be such that $\inr{\E A, Z^*-Z} > r^* $ and define $Z^\prime$ such that $Z^\prime -Z^* = \left(r^*/\inr{\E A, Z^*-Z}\right)(Z-Z^*)$. We have $\inr{\E A, Z^*-Z^\prime} = r^*$ and $Z^\prime\in\cC$ because $\cC$ is star-shaped in $Z^*$. Therefore, on the event $\Omega^*$,  $\inr{A-\E A, Z^\prime-Z^*}\leq (1/2)r^*$ and so $\inr{A-\E A, Z-Z^*}\leq (1/2)\inr{\E A, Z^*-Z}$. It therefore follows that on the event $\Omega^*$, if $Z \in \cC$ is such that $\inr{\E A, Z^*-Z} > r^*$ then
	\begin{equation*}
	\inr{A, Z - Z^*}\leq (-1/2)\inr{\E A, Z^* - Z}<-r^*/2 
	\end{equation*}
	which implies that $\inr{A, Z - Z^*} < 0$ and therefore $Z$ does not maximize $Z\to\inr{A, Z}$ over $\cC$. As a consequence, we necessarily have $\inr{\E A, Z^* - \hat{Z}} \leq r^*$ on the event $\Omega^*$ (which holds with probability at least $1-\Delta$). 
	
	Let us now assume that $r^*=0$. There exists a decreasing sequence $(r_n)_n$ of positive real numbers tending to $r^*=0$ such that for all $n\geq0$, $\bP[\Omega_n]\geq 1-\Delta$ where $\Omega_n$ is the event $\Omega_n = \left\{\psi(r_n)\leq \theta/2\right\}$ where for all $r>0$,
	\begin{equation*}
	\psi(r) = \frac{1}{r}\sup_{Z\in\cC: \inr{\E A, Z^*-Z}\leq r}\inr{A-\E A, Z-Z^*}.
	\end{equation*}Since $\cC$ is star-shapped in $Z^*$, $\psi$ is a non-increasing function and so  $(\Omega_n)_n$ is a decreasing sequence (i.e. $\Omega_{n+1}\subset \Omega_n$ for all $n\geq0$). It follows that $\bP[\cap_{n}\Omega_n] = \lim_n \bP[\Omega_n]\geq 1-\Delta$. Let us now place ourselves on the event $\cap_{n}\Omega_n$. For all $n$, since $\Omega_n$ holds and $r_n>0$, we can use the same argument as in  first case to conclude that $\inr{\E A , Z^*-\hat Z }\leq r_n$. Since the latter inequality is true for all $n$ (on the event $\cap_{n}\Omega_n$) and $(r_n)_n$ tends to zero, we conclude that $\inr{\E A , Z^*-\hat{Z}} \leq 0= r^*$. 
	\endproof
	
	The main conclusion of Theorem~\ref{theo:main} is that all the information for the problem of estimating $Z^*$ via $\hat Z$ is contained in the fixed point $r^*(\Delta)$. We therefore have to compute or upper bound such a fixed point. This might be difficult in great generality but there are some tools that can help to find upper bounds on $r^*(\Delta)$. 
	
	A first approach is to understand the shape of the local sets  $\cC\cap\{Z:\inr{\E A, Z^*-Z}\leq r\}, r>0$ and to that end it is helpful to characterize the \textit{curvature} of the excess risk $Z\to \inr{\E A, Z^*-Z}$ around its maximizer $Z^*$. This type of local characterization of the excess risk is also a tool used in Learning theory that goes back to classical conditions such as the  Margin assumption \cite{MR2051002,MR1765618} or the Bernstein condition \cite{MR2240689}. The latter condition was initially introduced as an upper bound of the variance term by its expectation: for all $Z \in \cC, \; \E(\ell_A(Z)-\ell_A(Z^*))^2\leq c_0 \E(\ell_A(Z)-\ell_A(Z^*))$ for some absolute constant $c_0$, but it has now been better understood as a way to discriminate the oracle from the other points in the model $\cC$. These assumptions were \textit{global} assumption in the sense that they concern all $Z$ in $\cC$. It has been recently shown \cite{chinot2018statistical} that only the \textit{local curvature}  of the excess risk needs to be understood. We now introduce this tool in our setup.
	
	We characterize the local curvature of the excess risk by some function $G:\bR^{n\times n}\to \bR$. Most of the time the $G$ function is a norm like the $\ell_1$-norm or a power of a norm such as the $\ell_2$ norm to the square. The radius defining the local subset onto which we need to understand the curvature of the excess risk is also solution of a fixed point equation:
	\begin{equation}\label{eq:fixed_point_G}
	r^*_G(\Delta) = \inf\left(r>0 : \bP\left[\sup_{Z\in\cC: G(Z^*-Z)\leq r}\inr{A-\E A, Z-Z^*}\leq (1/2) r\right]\geq 1-\Delta\right).
	\end{equation}The difference between the two fixed points $r^*(\Delta)$ and $r^*_G(\Delta)$ is that the local subsets are not defined using the same proximity function to the oracle $Z^*$; the first one uses the excess risk as a proximity function while the second one uses the $G$ function as a proximity function. The $G$ function should play the role of a \textit{simple} description of the curvature of the excess risk function locally around $Z^*$; that is formalized in the next assumption.

	\begin{Assumption}\label{ass:curvature}
		For all $Z\in\cC$, if $\inr{\E A ,Z^*-Z}\leq r^*_G(\Delta)$ then $\inr{\E A, Z^* - Z}\geq G(Z^*-Z)$.
	\end{Assumption}Typical examples of curvature function $G$ will have the form $G(Z^*-Z)=\theta\norm{Z^*-Z}^\kappa$ for some $\kappa\geq1$, $\theta>0$ and some norm $\norm{\cdot}$. In that case, the parameter $\kappa$ was initially called the \textit{margin parameter}  \cite{TsyCOLT07,MR1765618}. %We will derive estimation upper bounds for $\norm{\hat Z-Z^*}$, where the $[0,1]$-homogeneous function $\norm{\cdot}$ is the one characterizing the ``curvature'' of the objective function $Z\in\cC\to\inr{\E A, Z}$ around its maximum $Z^*$ in $\cC$,  as defined in Assumption~\ref{ass:curvature}. 
	% Even though margin condition or Bernstein condition have been the classical way to call a relation such as \eqref{eq:curvature} in the learning theory literature 
	Even though the relation given in Assumption~\ref{ass:curvature} has been typically referred to as a margin condition or Bernstein condition in the learning theory literature, 
	we will rather call it a \textit{local curvature assumption}, following \cite{guedon2016community} and \cite{chinot2018statistical}, since this type of relation describes the behavior of the risk function locally around its oracle. The main advantage for finding a local curvature function $G$ is that $r^*_G(\Delta)$ should be easier to compute than $r^*(\Delta)$ and $r^*(\Delta)\leq r^*_G(\Delta)$ because of the definition of $r_G^*(\Delta)$ and $\{Z\in\cC: \inr{\E A, Z^*-Z}\leq r_G^*(\Delta)\}\subset \{Z\in\cC: G(Z^*-Z)\leq r_G^*(\Delta)\}$ (thanks to Assumption~\ref{ass:curvature}). We can therefore state the following corollary.
	
	\begin{Corollary}\label{coro:main_coro}
		We assume that the constraint $\cC$ is star-shaped in $Z^*$ and that the ``local curvature'' Assumption~\ref{ass:curvature} holds  for some $0<\Delta<1$. With probability at least $1-\Delta$, it holds true that 
		\begin{equation*}
		r^*_G(\Delta)\geq \inr{\E A, Z^*-\hat Z}\geq G(Z^*-\hat Z).
		\end{equation*}
	\end{Corollary}
	When it is possible to describe the local curvature of the excess risk around its oracle by some $G$ function and when some estimate of $r^*_G(\Delta)$ can be obtained, Corollary~\ref{coro:main_coro} applies and estimation results of $Z^*$ by $\hat Z$ (w.r.t. to both the ''excess risk'' metric $\inr{\E A, Z^*-\hat Z}$ and the $G$ metric)  follow. If not, either because understanding the local curvature of the excess risk or the computation of $r^*_G(\Delta)$ is difficult, it is still possible to apply Theorem~\ref{theo:main} with the \textit{global} VC approach which boils down to simply upper bound the fixed point $r^*(\Delta)$ used in Theorem~\ref{theo:main} by a global parameter that is a complexity measure of the entire set $\cC$:
	\begin{equation}\label{eq:global_approach}
	r^*(\Delta)\leq \inf\left(r>0 : \bP\left[\sup_{Z\in\cC}\inr{A-\E A, Z-Z^*}\leq (1/2) r\right]\geq 1-\Delta\right).
	\end{equation}Interestingly, if the latter last resort approach is used then, following the approach form \cite{guedon2016community}, Grothendieck's inequality \cite{MR94682,pisier2012grothendieck} appears to be a powerful tool to upper bound the right-hand side of  \eqref{eq:global_approach} in the case of the community detection problem such as in \cite{MR3520025} as well as in the MAX-CUT problem. Of course, when it is possible to avoid this ultimate global approach one has to do it because the local approach will always provide better results.
	
	Finally, proving a ``local curvature'' property such as in Assumption~\ref{ass:curvature} may be difficult because it requires to understand the shape of the local subsets $\cC\cap\{Z:\inr{\E A, Z^*-Z}\leq r\}, r>0$. It is however possible to simplify this assumption if getting estimation results of $Z^*$ only w.r.t. the $G$ function (and not necessarily an upper bound on the excess risk $\inr{\E A, Z^*-\hat Z}$) is enough. In that case, Assumption~\ref{ass:curvature} may be replaced by the following one.
	\begin{Assumption}\label{ass:curvature_2}
		For all $Z\in\cC$, if $G(Z^*-Z)\leq r^*_G(\Delta)$ then $\inr{\E A, Z^* - Z}\geq G(Z^*-Z)$.
	\end{Assumption}Assumption~\ref{ass:curvature_2} assumes a curvature of the excess risk function in a $G$ neighborhood of $Z^*$ unlike Assumption~\ref{ass:curvature} which grants this curvature in an `` excess risk neighborhood''. The shape of a neighborhood defined by the $G$ function may be easier to understand (for instance when $G$ is a norm, a neighborhood defined by $G$ is the ball of a norm centered at $Z^*$ with radius $r_G^*(\Delta)$). In general, the latter assumption and Assumption~\ref{ass:curvature} do not compare. If Assumption~\ref{ass:curvature_2} holds then $\hat Z$ can estimate $Z^*$ w.r.t. the $G$ function.
	
	\begin{Theorem}\label{theo:main_2}
		We assume that the constraint $\cC$ is star-shaped in $Z^*$ and that the ``local curvature'' Assumption~\ref{ass:curvature_2} holds  for some $0<\Delta<1$. We assume that the $G$ function is continuous, $G(0)=0$ and $G(\lambda (Z^*-Z))\leq \lambda G(Z^*-Z)$ for all $\lambda\in[0,1], Z\in Z^*-\cC$. With probability at least $1-\Delta$, it holds true that $G(Z^*-\hat Z)\leq r^*_G(\Delta)$.
	\end{Theorem}
	\textbf{Proof of Theorem~\ref{theo:main_2}.} Let $r^* = r^*_G(\Delta)$. First assume that $r^*>0$. Let $Z\in\cC$ be such that $G(Z^*-Z)>r^*$. Let $f:\lambda\in[0,1]\to G(\lambda(Z^*-Z))$. We have $f(0)=G(0)=0$, $f(1)=G(Z^*-Z)>r^*$ and $f$ is continuous. Therefore, there exists $\lambda_0\in(0,1)$ such that $f(\lambda_0)=r^*$. We let $Z^\prime$ be such that $Z^\prime - Z^* = \lambda_0(Z-Z^*)$. Since $\cC$ is star-shapped in $Z^*$ and $\lambda_0\in[0,1]$ we have $Z^\prime\in \cC$. Moreover, $G(Z^*-Z^\prime) = r^*$. As a consequence, on the event $\Omega^*$ such that for all $Z\in\cC$ if $G(Z^*-Z)\leq r^*$ then $\inr{A-\bE A, Z-Z^*}\leq (1/2)r^*$, we have $\inr{A-\bE A, Z^\prime-Z^*}\leq (1/2)r^*$. The latter and Assumption~\ref{ass:curvature_2} imply that, on $\Omega^*$, 
	\begin{equation*}
	(1/2)r^*\geq \inr{A, Z^\prime - Z^*} + \inr{\bE A, Z^*-Z^\prime}\geq \inr{A, Z^\prime - Z^*} + G(Z^*-Z^\prime)\geq \inr{A, Z^\prime - Z^*} + r^*
	\end{equation*}and so $\inr{A, Z^\prime - Z^*}\leq -r^*/2$. Finally, using the definition of $Z^\prime$, we obtain 
	\begin{equation*}
	\inr{A, Z - Z^*} = (1/\lambda_0)\inr{A, Z^\prime - Z^*}\leq -r^*/(2\lambda_0)<0.
	\end{equation*} In particular, $Z$ cannot be a maximizer of $Z\to\inr{A, Z}$ over $\cC$ and so necessarily, on the event $\Omega^*$, $G(Z^*-\hat Z)\leq r^*$.
	
	Let us now consider the case where $r^*=0$. Using the same approach as in the proof of Theorem~\ref{theo:main}, we only have to check that the function
	\begin{equation*}
	\psi:r>0 \to \frac{1}{r}\sup_{Z\in\cC: G(Z^*-Z)\leq r}\inr{A-\E A, Z-Z^*}
	\end{equation*}is non-increasing. Let $0<r_1<r_2$. W.l.o.g. we may assume that there exists some $Z_2\in\cC$ such that $G(Z^*-Z_2)\leq r_2$ and $\psi(r_2)=\inr{A-\E A, Z_2-Z^*} / r_2$. If $G(Z^*-Z_2)\leq r_1$ then $\psi(r_2)\leq (r_1/r_2)\psi(r_1)\leq \psi(r_1)$. If $G(Z^*-Z_2)> r_1$ then there exists $\lambda_0\in(0,1)$ such that for $Z_1 = Z^* + \lambda_0(Z_2-Z^*)$ we have $G(Z^*-Z_1)=r_1$ and $Z_1\in\cC$. Moreover, $r_1=G(\lambda_0(Z^*-Z_2))\leq \lambda_0 G(Z^*-Z_2)\leq \lambda_0 r_2$ and so $\lambda_0\geq r_1/r_2$. It follows that 
	\begin{equation*}
	\psi(r_2) = \frac{1}{r_2}\inr{A-\E A, Z_2-Z^*} = \frac{1}{\lambda_0 r_2}\inr{A-\E A, Z_1-Z^*}\leq\frac{r_1}{\lambda_0 r_2}\psi(r_1)\leq \psi(r_1) 
	\end{equation*}where we used that $\psi(r)>0$ for all $r>0$ because $Z^*\in\{Z\in\cC:G(Z^*-Z)\leq r\}$ for all $r>0$.
	\endproof
	
	As a consequence, Theorem~\ref{theo:main}, Corollary~\ref{coro:main_coro} and Theorem~\ref{theo:main_2} are the three  tools at our disposal to study the performance of SDP estimators depending on the deepness of understanding we have on the problem. The best approach is given by Theorem~\ref{theo:main} when it is possible to compute efficiently the complexity fixed point $r^*(\Delta)$. If the latter approach is too complicated (likely because understanding  the geometry of the local subset $\cC\cap\{Z:\inr{\bE A, Z^*-Z}\leq r\},r>0$ may be difficult) then one may resort to find a curvature function $G$ of the excess risk locally around $Z^*$. In that case, both Corollary~\ref{coro:main_coro} and Theorem~\ref{theo:main_2} may apply depending on the hardness to find a local curvature function $G$ on an ``excess risk neighborhood'' (see Assumption~\ref{ass:curvature}) or a ``$G$-neighborhood'' (see Assumption~\ref{ass:curvature_2}). Finally, if no local approach can be handled (likely because describing the curvature of the excess risk in any neighborhood of $Z^*$ or controlling the maximal oscillations of the empirical process $Z\to \inr{\bE A-A,Z^*-Z}$ locally are too difficult) then one may resort ultimately to a global approach which follows from Theorem~\ref{theo:main} as explained in \eqref{eq:global_approach}. In the following, we will use these tools for various problems.

	Results like Theorem~\ref{theo:main}, Corollary~\ref{coro:main_coro} and Theorem~\ref{theo:main_2} appeared in many papers on ERM in learning theory such as in \cite{MR2829871,MR2240689,MR2319879,LM13}. In all these results, typical loss functions such as the quadratic or logistic loss functions were not linear one such as the one we are using here. From that point of view our problem is easier and this can be seen by the simplicity to prove our three general results from this section. What is much more complicated here than in other more classical problems in Learning Theory is the computation of the fixed point because 1) the stochastic processes $Z\to \inr{A-\bE A, Z-Z^*}$ may be far from being a Gaussian process if the noise matrix $A-\bE A$ is complicated and 2) the local sets $\{Z\in\cC: \inr{\E A, Z^*-Z}\leq r\}$ or $\{Z\in\cC: G(Z^*-Z)\leq r\}$ for $r>0$ maybe very hard to describe in a simple way. Fortunately, we will rely on several previous papers such as \cite{MR3901009} to solve such problems.

\section{Revisiting two results from the community detection literature \cite{MR3901009,guedon2016community}}
	\label{sec:application_to_the_community_detection_problem}
	% \cite{guedon2016community} and \cite{MR3901009} for the community detection problem

	The rapid growth of social networks on the Internet has lead many statisticians and computer scientists to focus their research on data coming from graphs. One important topic that has attracted particular interest during the last decades is that of community detection \cite{fortunato2010community,porter2009communities}, where the goal is to recover mesoscopic structures in a network, the so-called called communities. A community consists of group of nodes that are relatively densely connected to each other, but sparsely connected to other dense groups present within the network. The motivation for this line of work stems not only from the fact that finding communities in a network is an interesting and challenging problem of its own as it leads to understanding structural properties of networks, but community detection is also used as a data pre-processing step for other statistical inference tasks on large graphs, as it facilitates parallelization and allows one to distribute time consuming processes on several smaller subgraphs (i.e., the extracted communities).
	
	One challenging aspect of the community detection problem arises in the setting of sparse graphs. Many of the existing algorithms, which enjoy theoretical guarantees, do so in the relatively dense regime for the edge sampling probability, where the expected average degree is of the order $\Theta(\log n)$. The problem becomes challenging in very sparse graphs with bounded average degree. To this end, Gu\'{e}don and Vershynin proposed a semidefinite relaxation for a  discrete optimization problem  \cite{guedon2016community}, an instance of which encompasses the community detection problem, and showed that it can recover a solution with any given relative accuracy even in the setting of very sparse graphs with average degree of order $O(1)$.
	
	A subset of the existing literature for community detection and clustering relies on spectral methods, which consider the adjacency matrix associated to a graph, and employ its eigenvalues, and especially eigenvectors, in the analysis process or to propose efficient algorithms to solve the task at hand. Along these lines, Can et al. \cite{MR3449772} proposed a general framework for optimizing a general function of the graph adjacency matrix over discrete label assignments by projecting onto a low-dimensional subspace spanned by vectors that approximate the top eigenvectors of the expected adjacency matrix. The authors consider the problem of community detection with $k=2$ communities, which they frame as an instance of their proposed framework, combined with a regularization step that shifts each entry in the adjacency matrix by a small constant $\tau$, which renders their methodology applicable in the sparse regime as well.
	
	In the remainder of this section, we  focus on the community detection problem on random graphs under the general stochastic block model. We will mostly revisit the work in \cite{guedon2016community} and \cite{MR3901009} from the perspective given by Theorem~\ref{theo:main}, which simplifies the proof in \cite{MR3901009} since the peeling argument is no longer required, and neither is the upper bound from \cite{guedon2016community} unlike \cite{MR3901009}, thanks to the homogeneity argument hidden in Theorem~\ref{theo:main} (which underlies the localization argument).
	
	We first recall the definition of the generalized stochastic block model (SBM). We consider a set of vertices $V=\{1, \cdots, n\}$, and assume it is partitioned into $K$ communities $\cC_1, \cdots, \cC_K$ of arbitrary sizes $|\cC_1|=l_1, \cdots, |\cC_K|= l_K$. 
	\begin{Definition}
		For any pair of nodes $i,j\in V$, we denote by $i \sim j$ when $i$ and $j$ belong to the same community (i.e., there exists $k\in\{1, \ldots, K\}$) such that $i,j\in\cC_k$), and we denote by $i\not\sim j$ if $i$ and $j$ do not belong to the same community.
	\end{Definition}
	
	For each pair $(i, j)$ of nodes from $V$, we draw an edge between $i$ and $j$ with a fixed probability $p_{ij}$ independently from the other edges. We assume that there exist numbers $p$ and $q$ satisfying $0<q<p<1$, such that
	\begin{equation}\label{eq:proba_com_detect}
	\left\{ \begin{array}{ll}
	p_{ij} > p \mbox{, if $i\sim j$ and $i\neq j$}, \\
	p_{ij} = 1 \mbox{, if } i=j, \\
	p_{ij} < q \mbox{, otherwise.}
	\end{array}
	\right.
	\end{equation}
	We denote by $A=(A_{i,j})_{1\leq i,j,\leq n}$ the observed symmetric adjacency matrix, such that, for all $1\leq i\leq j\leq n$,  $A_{ij}$ is distributed according to a Bernoulli of parameter $p_{ij}$.
	The community structure of such a graph is captured by the membership matrix $\bar{Z} \in \R^{n\times n}$, defined by $\bar{Z}_{ij}=1$ if $i\sim j$, and $\bar{Z}_{ij}=0$ otherwise. The main goal in community detection is to reconstruct $\bar{Z}$ from the observation $A$.
	
	Spectral methods for community detection are very popular in the literature  \cite{guedon2016community,MR3901009,vershynin2018high,blondel2008fast,clauset2004finding}. There are many ways to introduce such methods, one of which being via convex relaxations of certain graph cut problems aiming to minimize a modularity function such as the RatioCut \cite{newman2006finding}. Such relaxations often lead to SDP estimators, such as those introduced in Section~\ref{sec:introduction}.
	
	Considering a random graph distributed according to the generalized stochastic block model, and its associated adjacency matrix $A$ (i.e. $A=A^\top$ and $A_{ij}\sim {\rm Bern}(p_{ij})$ for $1\leq i\leq j\leq n$ and $p_{ij}$ as defined in \eqref{eq:proba_com_detect}), we will estimate its membership matrix $\bar{Z}$ via the following SDP estimator
	\begin{equation*}
	\hat{Z} \in \argmax_{Z \in \mathcal{C}}  \inr{A, Z}, 
	\end{equation*}
	where $\mathcal{C} = \{ Z \in \R^{n\times n}, Z \succeq 0, Z \geq 0, {\rm diag}(Z)\preceq I_n, \sum_{i, j=1}^n Z_{ij} \leq \lambda \}$ and  $\lambda = \sum_{i, j=1}^n \bar{Z}_{ij} = \sum_{k=1}^K |\mathcal{C}_k|^2$ denotes the number of nonzero elements in the membership matrix $\bar{Z}$. The motivation for this approach stems from the fact that the membership matrix  $\bar{Z}$ is actually the oracle, i.e., $Z^*=\bar Z$ (see  Lemma~7.1 in \cite{guedon2016community} or Lemma~\ref{lem:curvature_com_detect} below), where
	\begin{equation*}
	Z^* \in \argmax_{Z \in \mathcal{C}} \inr{\E A, Z}.
	\end{equation*}

	Following the strategy from Theorem~\ref{theo:main} and from our point of view, the upper bound on $r^*(\Delta)$ from \cite{guedon2016community} is the one that is based on the global approach -- that is, without localization. Indeed, \cite{guedon2016community} uses the observation that, for all $r>0$, it holds true that 
	\begin{equation}\label{eq:non-local}
	\sup_{Z\in\cC: \inr{\E A , Z^*-Z}\leq r}  \inr{A-\bE A, Z-Z^*}  \overset{(a)}{\leq}\sup_{Z\in\cC}  \inr{A-\bE A, Z-Z^*} \overset{(b)}{\leq}  2K_G \norm{A-\bE A}_{\mathrm{cut}},
	\end{equation}
	where  $\norm{\cdot}_{{\rm cut}}$ is the cut-norm\footnote{The cut-norm  $\norm{\cdot}_{{\rm cut}}$ of a real matrix $A = (a_{ij})_{i \in R, j \in C}$ with a set of rows indexed by $R$ and a set of columns indexed by $C$, is the maximum, over all $I \subset R$ and $J \subset C$, of the quantity $ | \sum_{i \in I, j \in J} a_{ij}| $. It is also the operator norm of $A$ from $\ell_\infty$ to $\ell_1$ and the ``injective norm'' in the orginal Grothendieck ``r{\'e}sum{\'e}'' \cite{grothendieck1956resume,pisier2012grothendieck}}
	and $K_G$ is the Grothendieck constant (Grothendieck's inequality is used in (b), see \cite{pisier2012grothendieck,vershynin2018high}). Therefore, the localization around the oracle $Z^*$ by the excess risk ``band'' $B^*_r:=\{Z:\inr{\E A, Z^*-Z}\leq r\}$  is simply removed in inequality~{(a)}. As a consequence, the resulting statistical bound is based on the complexity of the entire class $\cC$ whereas, in a localized approach, only the complexity of $\cC\cap B^*_r$ matters. Next step in the proof of \cite{guedon2016community} is a high probability upper bound on $\norm{A-\bE A}_{\mathrm{cut}}$  which follows from Bernstein's inequality and a union bound since one has $\norm{A-\bE A}_{\mathrm{cut}} = \max_{x, y \in \{-1, 1\}^n} \inr{A-\bE A, xy^\top}$, then for all $t>0$, $\norm{A - \E A}_{\mathrm{cut}} \leq tn(n-1)/2$ with probability at least $1 - \exp\left(2n\log 2 - (n(n-1)t^2)/(16\bar{p} + 8t/3)\right)$ where $\bar{p} \overset{def}{=} 2/[n(n-1)] \sum_{i<j} p_{ij}(1 - p_{ij})$. The resulting upper bound on the fixed point obtained in \cite{guedon2016community} is, 
	\begin{equation}\label{eq:guedon_vershynin_upper_bound_r_star}
	r^*(\Delta) \leq (8/3) K_G (2n\log(2)+\log(1/\Delta)).
	\end{equation}
	Finally, under the assumption of Theorem~1 in \cite{guedon2016community} (i.e., for some some $\epsilon \in (0, 1)$, $n\geq 5.10^4/\eps^2$, $\max(p(1-p), q(1-q))\geq20/n$,  $p = a/n>b/n = q$ and  $(a-b)^2 \geq 2.10^4 \epsilon^{-2}(a+b)$), for $\Delta = e^35^{-n}$ we obtain (using the general result in Theorem~\ref{theo:main}) with probability at least  $1-\Delta$, the bound $\inr{\E A , Z^* - \hat Z}\leq r^*(\Delta)\leq \eps n^2 = \eps \norm{Z^*}_2^2$, which is the result from Theorem~1 in \cite{guedon2016community}.  Finally, \cite{guedon2016community} uses a (global) curvature property of the excess risk in its Lemma~7.2: 
	\begin{lemma}[Lemma~7.2 in \cite{guedon2016community}]\label{lem:curvature_com_detect}
		For all $Z \in \mathcal{C}$, $\inr{\E A, Z^* - Z} \geq [(p - q)/2] \norm{Z^*-Z}_1$. 
	\end{lemma}
	Therefore, a (global-- that is for all $Z\in\cC$) curvature assumption holds for a $G$ function which is here the $\ell_1^{n\times n}$ norm, a margin parameter $\kappa=1$ and $\theta=(p-q)/2$ for the community detection problem. However this curvature property is not use to compute a ``better'' fixed point parameter but only to have a $\ell_1^{n\times n}$ estimation bound since
	\begin{equation*}
	\norm{\hat Z - Z^*}_1\leq \left(\frac{2}{p-q}\right)\inr{\E A, Z^*-\hat Z}\leq \frac{16 K_G (2n\log(2)+\log(1/\Delta))}{3(p-q)}.
	\end{equation*}The latter bound together with the sin-Theta theorem allow the authors from \cite{guedon2016community} to obtain estimation bound for the community membership vector $x^*$.

	The approach from \cite{MR3901009} improves upon the one in \cite{guedon2016community} because it uses a localization argument: the curvature property of the excess risk function from Lemma~\ref{lem:curvature_com_detect} is used to improve the upper bound in \eqref{eq:guedon_vershynin_upper_bound_r_star} obtained following a global approach. Indeed, the authors from  \cite{MR3901009} obtain high probability upper bound on the quantity
	\begin{equation*}
	\sup_{Z\in\cC: \norm{ Z^*-Z}_1\leq r}  \inr{A-\bE A, Z-Z^*}
	\end{equation*}depending on $r$. This yields to exact reconstruction result in the ``dense'' case and exponentially decaying rates of convergence in the ``sparse'' case. This is a typical example where the localization argument shows its advantage upon the global approach. The price to pay is usually a more technical proof for the local approach compare with the global one. However, the argument from \cite{MR3901009} also uses an unnecessary peeling argument together with an unnecessary a priori upper bound on $\norm{\hat Z-Z^*}_1$ (which is actually the one from  \cite{guedon2016community}). It appears that this peeling argument and this a priori upper bound on $\norm{\hat Z-Z^*}_1$ can be avoided thanks to our approach from Theorem~\ref{theo:main}. This improves the probability estimate and simplifies the proofs (since the result from \cite{guedon2016community} is not required anymore neither is the peeling argument). For the sign clustering problem we consider below as an application of our main results, we will mostly adapt the probabilistic tools from \cite{MR3901009} (in the ``dense'' case) to the methodology associated with Theorem~\ref{theo:main} (without this two unnecessary arguments).

	%%%%%%%%%%%%%%%%%%%%%%%%%%%%%%%%%%%%%%%%%%%%%%%%%%%%%%%%%%%%%%%%%%%%%%%%%%-----------------------
	%%%%%%%%%%%%%%%%%%%%%%%%%%%%%%%%%%%%%%%%%%%%%%%%%%%%%%%%%%%%%%%%%%%%%%%%%%-----------------------
	
	\section{Contributions of the paper}
	\subsection{Application to signed clustering} % (fold)
	\label{sec:application_to_signed_clustering}
	
	Much of the clustering literature, including both spectral and non-spectral methods, has focused on unsigned graphs, where each edge carries a non-negative scalar weight that encodes a measure of affinity (similarity, trust) between pairs of nodes. However, in numerous instances, the above-mentioned affinity takes negative values, and encodes a measure of dissimilarity or distrust. Such applications arise in social networks where users relationships denote trust-distrust or friendship-enmity, shopping bipartite networks which capture like-dislike  relationships between users and products \cite{banerjee2012partitioning}, 
	online news and review websites, such as Epinions \cite{epinions} and Slashdot \cite{slashdot}, that allow users to approve or denounce others  \cite{Leskovec_2010_PPN}, and clustering financial or economic time series data \cite{aghabozorgi2015time}.
	Such applications have spurred interest in the analysis of signed networks, which has recently become an increasingly important research topic \cite{Leskovec_2010_SNS}, with relevant lines of work in the context of clustering signed networks including, in chronological order,  \cite{kunegis2010spectral,Chiang_2012_Scalable,SPONGE}.
	
	The second application of our proposed methodology is an extension of the community detection and clustering problems to the setting of signed graphs, where, for simplicity, we assume that an edge connecting two nodes can take either $-1$ or $+1$ values.

	\subsubsection{A Signed Stochastic Block Model (SSBM)}
	\label{sec:SSBM}
	We focus on the problem of clustering a K-weakly balanced graphs\footnote{A signed graph is K-weakly balanced if and only if all the edges are positive, or the nodes can be partitioned into $K \in \mathbb{N}$ disjoint sets such that positive edges exist only within clusters, and negative edges are only present across clusters \cite{DavisWeakBalance}.}.
	We consider a signed stochastic block model (SSBM) similar to the one introduced in \cite{SPONGE}, where we are given a graph $G$ with $n$ nodes $\{1, \ldots, n\}$ which are divided into $K$ communities, $\{\cC_1, \cdots, \cC_K\}$, such that, in the noiseless setting, edges within each community are positive and edges between communities are negative.
	
	The only information available to the user is given by a $n\times n$ sparse  adjacency matrix $A$ constructed as follows: $A$ is symmetric, with $A_{ii}=1$ for all $i=1, \ldots, n$,  and for all $1\leq i<j\leq n$, $A_{ij} = s_{ij}(2B_{ij}-1)$ where  
	\begin{equation*}
	\mathrm{B}_{ij} \sim \left\{ \begin{array}{l}
	\mathrm{Bern}(p) \mbox{ if } i\sim j \\
	\mathrm{Bern}(q) \mbox{ if } i \not\sim j 
	\end{array}
	\right. \;\; \mbox{ and } \;\; s_{ij}\sim {\rm Bern}(\delta),
	\end{equation*}
	for some $0\leq q<1/2<p\leq 1$ and $\delta\in(0,1)$. Moreover, all the variables $B_{ij}, s_{ij}$ for $1\leq i<j\leq n$ are independent.
	
	We remark that this SSBM model is similar to the one considered in \cite{SPONGE}, which was governed by two parameters, the sampling probability $\delta$ as above, and the noise level $\eta$, which may flip entries of the adjacency matrix.

	Our aim is to recover the community membership matrix or cluster matrix $\bar Z = (\bar Z_{ij})_{n\times n}$where $\bar Z_{ij}=1$ when $i\sim j$ and $\bar Z_{ij}=0$ when $i\not\sim j$ using only the observed censored adjacency matrix $A$.
	
	Our approach is similar in nature to the one used by spectral methods in community detection. We first observe that for $\alpha :=\delta(p+q-1)$ and $J=(1)_{n\times n}$ we have $\bar Z = Z^*$ where
	\begin{equation}\label{eq:oracle_egal_cluster_mat}
	Z^*\in\argmax_{Z\in \cC}\inr{\E A - \alpha J, Z}
	\end{equation}and $\cC = \{Z\in\bR^{n\times n}: Z\succeq 0, Z_{ij}\in[0,1], Z_{ii}=1, i=1, \ldots, n\}$. The proof of \eqref{eq:oracle_egal_cluster_mat} is recalled in Section~\ref{sub:annex_2_in_signed_clustering}.
	
	Since we do not know $\bE A$ and $\alpha$, we should estimate both of them. We will estimate $\bE A$ with $A$ but, for simplicity, we will assume that $\alpha$ is known. The resulting estimator of the cluster matrix $\bar Z$ is
	\begin{equation}\label{eq:esti_signed_cluster}
	\hat{Z} \in\argmax_{Z\in\cC}  \inr{A - \alpha J, Z}
	\end{equation}which is indeed a SDP estimator and therefore Theorem~\ref{theo:main} (or Corollary~\ref{coro:main_coro} and Theorem~\ref{theo:main_2}) may be used to obtain statistical bounds for the estimation of $Z^*$ from \eqref{eq:oracle_egal_cluster_mat} by $\hat Z$.

	We will use the following notations: $s := \delta(p-q)^2$, $\theta :=\delta(p-q)$, $\rho := \delta\max\{1-\delta(2p-1)^2, 1-\delta(2q-1)^2\}$, $\nu := \max\{2p-1, 1-2q\} $, $[m] := \{1, \cdots, m \}$ for all $m \in \bN$, $l_k := |\cC_k|$  for all $k \in [K]$, $\lambda^2 := \sum_{k=1}^K l_k^2$, $\cC^+ := \underset{k=1}{\overset{K}{\cup}} (\cC_k \times \cC_k)$ and  $\cC^- := \underset{k \neq k^\prime}{\overset{}{\cup}} (\cC_k \times \cC_{k^\prime})$. We also use the notation $c_0,c_1, \ldots,$ to denote absolute constants whose values may change from one line to another.

	\subsubsection{Main result for the estimation of the cluster matrix in signed clustering}
	Our main result concerns the reconstitution of the $K$ communities from the observation of the matrix $A$. In order to avoid solutions with some communities of degenerated size (too small or too large) we consider the following assumption.
	
	\begin{Assumption}\label{ass1}
		Up to constant, the elements of the partition $\cC_1\sqcup\cdots\sqcup\cC_K$ of $\{1, \ldots, n\}$ are of same size (up to constants): there are absolute constant $c_0,c_1>0$ such that for any $k \in [K]$,  $n/(c_1K)\leq |\cC_k| = l_k \leq c_0n/K$. 
	\end{Assumption}

	\noindent We are now ready to state the main result on the estimation of the cluster matrix $Z^*$  from \eqref{eq:oracle_egal_cluster_mat} by the SDP estimator $\hat Z$ from \eqref{eq:esti_signed_cluster}.
	\begin{Theorem}\label{theo:main_signed_cluster}
		There is an absolute positive constant $c_0$ such that the following holds.
		Grant Assumption \ref{ass1}. Assume that 
		\begin{equation}
		n \nu \delta \geq \log n,    
		\label{eq:thm_cond_1_density}
		\end{equation}
		\begin{equation}
		(p-q)^2n \delta \geq c_0 K^2 \nu
		\label{eq:thm_cond_2_sep}
		\end{equation}
		\begin{equation}
		\frac{K \log(2 e K n)}{n}\leq  \max \left(\frac{\theta^2}{\rho}, \frac{9\rho}{32}\right).
		\label{eq:thm_cond_3_largeClust}
		\end{equation}
		Then, with probability at least $1 - \exp(-\delta \nu n)-3/(2e Kn)$, exact recovery holds true, i.e.,  $\hat Z= Z^*.$
	\end{Theorem}
	Therefore, we have exact  reconstruction in the dense case (that is under assumption \eqref{eq:thm_cond_1_density}),  which stems from condition  \eqref{eq:thm_cond_3_largeClust}.  
	% \rednote{?do you agree on that? It is more restrictive than \eqref{eq:thm_cond_2_sep}. It seems that $n$ is of the order $\nu$ and $1 / \nu$ at the same time. Does that imply that $\nu$ is somewhat independent of $n$?}   
	% \blue{Yes, I agree I look at $\nu$ more like a constant. Maybe I'm wrong do you think it is an issue.} 
	The latter condition is in the same spirit as the one in Theorem~1 of \cite{MR3901009}, it measures the SNR (signal-to-noise ratio) of the model which captures the hardness of the SSBM. As mentioned in \cite{MR3901009}, it is related to the Kesten-Stigum threshold \cite{MR3383334}. The last condition \eqref{eq:thm_cond_3_largeClust} basically requires that the number of clusters $K$ is at most $n/\log n$. If this condition is dropped out then we do not have anymore exact reconstruction but only a certified rate of convergence:   with probability at least $1 - \exp(-\delta \nu n)-3/(2e Kn)$,
	\begin{equation*}
	\norm{Z^*-\hat Z}_1\leq \frac{n}{c_1\delta(p-q)K}.
	\end{equation*} This shows that there is a phase transition in the dense case: exact reconstruction is possible when $K\lesssim n/\log n$ and, otherwise, when $K\gtrsim n/\log n$ we only have a control of the estimation error.
	%\bluenote{Do not hesitate do add more comment on our result comparing to others in the literature.}

	\subsection{Application to angular group synchronization} % (fold)
	\label{sec:application_to_synchronization}
	In this section, we introduce the group  synchronization problem as well as a  stochastic model for this problem. We consider a SDP relaxation of the original problem (which is exact) and construct the associated SDP estimator such as \eqref{eq:estimator_0}.
	
	The angular synchronization problem consists of estimating $n$ unknown angles $\theta_1, \cdots, \theta_n$ (up to a global shift angle) given a noisy subset of their pairwise  offsets $(\theta_i - \theta_j) [2 \pi]$, where $[2\pi]$ is the modulo $2\pi$ operation. The pairwise measurements can be realized as the edge set of a graph $G$, typically modeled as an Erd\"os-Renyi random graph \cite{sync}. 
	
	The aim of this section is to show that the angular synchronization problem can be analyzed using our methodology. In order to keep the presentation as simple as possible, we assume that all pairwise offsets are observed up to some Gaussian noise: we are given $(\theta_i - \theta_j + \sigma g_{ij}) [2 \pi]$ for all $1\leq i<j\leq n$ where $(g_{ij}:1\leq i<j\leq n)$ are $n(n-1)/2$ i.i.d. standard Gaussian variables and $\sigma>0$ is the noise variance. We may rewrite the problem as follows: we observe a $n\times n$ complex matrix $A$ defined by
	\begin{equation}  \label{eq:sync_noise_model}
	A = S \circ [x^*(\overline{x^*})^\top] \mbox{ where } S=(S_{ij})_{n\times n}, S_{ij}=\left\{
	\begin{array}{cc}
	e^{\iota \sigma g_{ij}} & \mbox{ if } i<j\\
	1 & \mbox{ if } i=j\\
	e^{-\iota \sigma g{ij}} & \mbox{ if } i>j
	\end{array}
	\right.,
	\end{equation}$\iota$ denotes the imaginary number such that $\iota^2=-1$, 
	$x^*=(x^*_i)_{i=1}^n\in\bC^n$, $x^*_i = e^{\iota \theta_i}, i=1, \ldots, n$, $\bar{x}$ denotes the conjugate vector of $x$ and $S \circ [x^*(\overline{x^*})^\top]$ is the element-wise product $(S_{ij}x_i \bar{x}_j)_{n\times n}$. In particular, $S$ is a Hermitian matrix  (i.e. $\bar S^\top = S$) and $\E S_{ij} = \exp(-\sigma^2/2)$ for $i\neq j$ and $\E S_{ii}=1$ if $i=j$.
	% and  $W \in \mathbb{C}^{n\times n}$ is from the GUE Ensemble (i.e. $W=(W_{ij})_{n\times n}$ is Hermitian and $(W_{ii})_{1}^n$, $(\sqrt{2}\frak{R}(W_{ij}))_{i<j}$ and $(\sqrt{2}\frak{I}(W_{ij}))_{i<j}$ are $n^2$ i.i.d. $\cN(0,1)$, where $\frak{R}$ and $\frak{I}$ stand for the real and imaginary parts), see Chapter~6 in \cite{MR2129906}.  In other words, we observe for all $i,j\in[n]$, $A_{ij} = e^{\iota \delta_{ij}} + \sigma g_{ij}$, where $\delta_{ij} \overset{def}{=} \theta_i - \theta_j$ are the offsets of the angles, and $g_{ij}\sim\cN(0,1)$ are the Gaussian entries of $W$. \rednote{Just to clarify here, the added noise does not necessarily have to be real - essentially, the angle measurement is perturbed, which makes $A$ as a Hermitian random perturbation of the rank-1 outer product clean matrix. Does this affect any the current analysis?}
	We want to estimate $(\theta_1, \ldots, \theta_n)$ (up to a global shift) from the matrix of data $A$. Unlike classical statistical models the noise here is multiplicative; we show that our approach covers this type of problem as well. 
	
	The first step is to find an (vectorial) optimization problem which solutions are given by $(\theta_i)_{i=1}^n$ (up to global angle shift) or some bijective function of it.  Estimating $(\theta_i)_{i=1}^n$ up to global angle shift is  equivalent to estimating the vector $x^*=(e^{\iota \theta_i})_{i=1}^n$. The latter is, up to a global rotation of its coordinates, the unique solution of the following maximization problem 
	\begin{equation}
	\label{eq:sync_x_star}
	\argmax_{x \in \mathbb{C}^n:|x_i| = 1} \left\{ \bar{x}^\top \; \bE A \; x \right\} = \{(e^{\iota (\theta_i+\theta_0)})_{i=1}^n:\theta_0\in[0,2\pi)\}.
	\end{equation} 
	A proof of \eqref{eq:sync_x_star} is given in Section~\ref{sec:proof_of_Theorem_in_synchronization}. Let us now rewrite \eqref{eq:sync_x_star} as a SDP problem. For all $x\in\bC^n$, we have $\bar{x}^\top\E Ax = \mathrm{tr}(\E AX) = \inr{\E A, X}$ where $X=x\bar{x}^\top$ and $ \{Z \in \mathbb{C}^{n \times n} : Z = x\bar{x}^T, |x_i|=1 \} = \{ Z \in \mathbb{H}_n: Z \succeq 0, \mathrm{diag}(Z) = \mathbf{1}_n, \mathrm{rank}(Z) = 1 \}$ where $\mathbb{H}_n$ is the set of all $n\times n$ Hermitian matrices and $\mathbf{1}_n\in\bC^n$ is the vector with all coordinates equal to $1$. It is therefore straightforward to construct a SDP relaxation of \eqref{eq:sync_x_star} by dropping the rank constraint. It appears that this relaxation is exact since, for $\cC = \{ Z \in \mathbb{H}_n: Z \succeq 0, \mathrm{diag}(Z) = \mathbf{1}_n \}$,
	\begin{equation}
	\label{eq:sync_X_star}
	\argmax_{Z\in \cC} \inr{\bE A, Z} = \{Z^*\},
	\end{equation}
	where $Z^* = x^*(\overline{x^*})^\top$. A proof of \eqref{eq:sync_X_star} can be found in  Section~\ref{sec:proof_of_Theorem_in_synchronization}. Finally, as we only observe $A$, we consider the following SDP estimator of $Z^*$
	\begin{equation}\label{pb:sync_SDP_est} 
	\hat{Z} \in\argmax_{Z\in \cC}\inr{A,Z}.
	\end{equation}
	In the next section, we use the strategy from Corollary~\ref{coro:main_coro} to obtain statistical guarantees for the estimation of $Z^*$ by $\hat Z$. 
	
	Intuitively, the above maximization problem \eqref{eq:sync_X_star} attempts to preserve the given angle offsets as best as possible
	\begin{equation}
	\argmax_{\theta_1,\ldots,\theta_n \in [0,2\pi)  } \sum_{i,j=1}^{n} e^{-\iota \theta_i} A_{ij} e^{\iota \theta_j},
	\label{eq:maxSync}
	\end{equation}
	where the objective function is incremented by $+1$ whenever an assignment of angles $\theta_i$ and $\theta_j$ perfectly satisfies the given edge constraint $  \delta_{ij} = (\theta_i - \theta_j) [2\pi]$ (i.e., for a \textit{clean} edge for which $\sigma=0$), while the contribution of an incorrect assignment (i.e., of a \textit{very noisy} edge) will be almost uniformly distributed on the unit circle in the complex plane. Due to non-convexity of optimization in \eqref{eq:maxSync}, it is difficult to solve computationally \cite{zhang2006complex}; one way to overcome this problem is to consider the SDP relaxation from \eqref{eq:sync_X_star} but it is also possible to consider a  spectral relaxation such as the one proposed by Singer \cite{sync} which replaces the individual constraints that all $z_i$'s should have unit magnitude by the much weaker single constraint $\sum_{i=1}^n |z_i|^2 = n $, leading to 
	\begin{equation}
	\argmax_{  z_1,\ldots,z_n \in \mathbb{C};  \;\;\; \sum_{i=1}^n |z_i|^2 = n }  \sum_{i,j=1}^{n}  \bar{z_i} A_{ij} z_j.
	\label{eq:maxSyncRelax}
	\end{equation}
	The solution to the resulting maximization problem is simply given by a top eigenvector of the Hermitian matrix $A$, followed by a normalization step. We remark that the main advantage of the SDP relaxation \eqref{eq:sync_X_star} is that it explicitly imposes the unit magnitude constraint for $e^{\iota \theta_i}$, which we cannot otherwise enforce in the spectral relaxation solved via the eigenvector method in \eqref{eq:maxSyncRelax}.
	The above SDP program  \eqref{eq:sync_X_star} is very similar to the well-known Goemans-Williamson SDP relaxation for the seminal \textsc{MAX-CUT} problem of 
	finding the maximal cut of a graph (the  \textsc{MAX-CUT} problem is one of the four applications considered in this work, see Section~\ref{sec:application_to_max_cut}), with the main difference being that here we optimize over the cone of complex-valued Hermitian positive semidefinite matrices, not just real symmetric matrices.

	\subsubsection{Main results for phase recovery in the synchronization problem}
	
	Our main result concerns the estimation of the matrix of offsets $Z^*= x^*(\overline{x^*})^\top$ from the observation of the matrix $A$. This result is then used to estimate (up to global phase shift) the angular vector $x^*$.

	\begin{Theorem}\label{theo:main_synchro}Let $0<\eps<1$. If $\sigma\leq \sqrt{\log(\eps n^4)}$ then, with probability at least $1 - \exp(-\eps \sigma^4n(n-1)/2)$, it holds true that 
		\begin{equation}
		(e^{-\sigma^2/2}/2)\norm{Z^*-Z}_2^2\leq \inr{\E A, Z^*-Z} \leq (128/3)\sqrt{\eps}\sigma^4 N.
		\end{equation}
	\end{Theorem}
	Once we have an estimator $\hat{Z}$ for the oracle $Z^*$, we can extract an estimator $\hat{x}$ for the vector of phases $x^*$ by considering a  top eigenvector (i.e. an eigenvector associated with the largest eigenvalue) of $\hat{Z}$. It is then possible to quantify the estimation properties of $x^*$ by $\hat x$ using a sin-theta theorem and Theorem~\ref{theo:main_synchro}.
	
	\begin{Corollary}\label{cor:sync}
		Let $\hat{x}$ be a top eigenvector of $\hat{Z}$ with Euclidean norm $\norm{\hat{x}}_2 = \sqrt{n}$. Let $0<\eps<1$ and assume that $\sigma\leq \sqrt{\log(\eps n^4)}$. We have the existence of a universal constant $c_0>0$ (which is the constant in the Davis-Kahan Theorem for Hermitian matrices) such that, with probability at least $1 - \exp(-\eps \sigma^4n(n-1)/2)$, it holds true that
		\begin{equation}
		\min_{z\in\bC:|z|=1}\norm{\hat{x}-zx^*}_2 \leq 8c_0\sqrt{2/3}\eps^{1/4}e^{\sigma^2/4} \sigma^2 \sqrt{n}.
		\end{equation}
	\end{Corollary}
	It follows from Corollary~\ref{cor:sync} that we can estimate $x^*$ (up to a global rotation $z\in\bC:|z|=1$) with a $\ell_2^n$-estimation error of the order of $\sigma^2 \sqrt{n}$ with exponential deviations. Given that $\norm{x^*}_2=\sqrt{n}$, this means that a constant proportion of the entries are well estimated. For a values of $\eps\sim1/n^2$, the rate of estimation is like $\sigma^2$, we therefore get a much better estimation of $x^*$ but only with constant probability. It is important to recall that $\hat Z$ and $\hat x$ can be both efficiently computed by solving a SDP problem and then by considering a top eigenvector of its solution (for instance, using the power method).
	
	\subsection{Application to the MAX-CUT problem} % (fold)
	\label{sec:application_to_max_cut}
	Let $A^0 \in \{0, 1 \}^{n\times n}$ be the adjacency (symmetric) matrix of an undirected graph $G = (V, E^0)$, where $V := \{1, \ldots, n\}$ is the set of the vertices and the set of edges is $E^0:=E\cup E^\top\cup\{(i,i):A_{ii}^0=1\}$ where $E:=\{(i, j)\in V^2~:~i<j \mbox{ and } A^0_{ij}=1\}$ and $E^\top = \{(j,i):(i,j)\in E\}$. We assume that $G$ has no self loop so that $A_{ii}^0=0$ for all $i\in V$. A \textit{cut} of $G$ is any subset $S$ of vertices in $V$. For a cut $S\subset V$, we define its weight by $\mathrm{cut}(G,S):= (1/2)\sum_{(i, j) \in S\times \bar{S}} A^0_{ij}$, that is the number of edges in $E$ between $S$ and its complement $\bar{S}=V\backslash S$. The \textit{MAX-CUT} problem is to find the cut with maximal weight:
	\begin{align}\label{pb:MAX-CUT0}
	S^* \in \underset{S\subset V}{\mathrm{argmax}}  \, \mathrm{cut}(G, S).
	\end{align}
	
	The MAX-CUT problem is a NP-complete problem but \cite{goemans1995improved} constructed a $0.878$ approximating solution via a SDP relaxation. Indeed, one can write the MAX-CUT problem in the following way.  For a cut $S\subset V$, we define the membership vector $x\in \{-1, 1 \}^n$ associated with $S$ by setting $x_i := 1$ if $i\in S$ and $x_i=-1$ if $i\notin S$ for all $i\in V$. We have $\mathrm{cut}(G,S) = (1/4)\sum_{i, j=1}^n A^0_{ij}(1-x_ix_j):=\mathrm{cut}(G, x)$ and so solving \eqref{pb:MAX-CUT0} is equivalent to solve
	\begin{align}\label{pb:MAX-CUT1}
	x^* \in \underset{x\in \{-1, 1\}^n}{\mathrm{argmax}}\,  \mathrm{cut}(G, x).
	\end{align}Since $(x_ix_j)_{i,j}=xx^\top$, the latter problem is also equivalent to solve 
	\begin{equation}\label{eq:max_cut_pb_3}
	\max\left(\frac{1}{4}\sum_{i,j=1}^n A^0_{ij}(1-Z_{ij}):{\rm rank}(Z)=1, Z\succeq 0, Z_{ii} = 1\right)
	\end{equation}
	which admits a SDP relaxation by removing the rank $1$ constraint. This yields the following SDP relaxation problem of MAX-CUT from \cite{goemans1995improved}:
	\begin{equation}\label{eq:SDP_relax_max_cut}
	Z^* \in \underset{Z\in \mathcal{C}}{\mathrm{argmin}} \inr{A^0, Z} 
	\end{equation}
	where $\mathcal{C} := \{Z \in \bR^{n\times n}:Z \succeq 0, Z_{ii} =  1, \forall i=1,\ldots, n\}$.
	
	Unlike the other examples from the previous sections, the SDP relaxation in \eqref{eq:SDP_relax_max_cut} is not exact, except for bipartite graphs; see \cite{khot2009approximate,gartner2012semidefinite} for more details. Nevertheless, thanks to the approximation result from \cite{goemans1995improved} we can use our methodology to estimate $Z^*$ and then deduce some approximate optimal cut. But first, we introduce a stochastic model because in many situations the adjacency matrix $A^0$ is only partially  observed but still it might be interesting to find an approximating solution to the MAX-CUT problem.

	We observe  $A = S\circ A^0=(s_{ij}A_{ij}^0)_{1\leq i,j\leq n}$ a ``masked'' version of $A^0$, where $S\in\bR^{n\times n}$ is symmetric with upper triangular matrix filled with i.i.d. Bernoulli entries: for all $i,j\in V$ such that $i\leq j$, $S_{ij}=S_{ji} = s_{ij}$ where $(s_{ij})_{i\leq j}$ is a family of i.i.d. Bernoulli random variables with parameter $p \in (1/2, 1)$. Let $B:=-(1/p)A$ so that $\bE[B] = -A^0$. We can write $Z^*$ as an oracle since $Z^*\in \argmax_{Z\in\cC}\inr{\bE B, Z}$ and so we estimate $Z^*$ via the SDP estimator $\hat{Z}\in \argmax_{Z\in\cC}\inr{B, Z}$. Our first aim is to quantify the cost we pay by using $\hat Z$ instead of $Z^*$ in our final choice of cut. It appears that the fixed point used in Theorem~\ref{theo:main} may be used to quantify this loss:
	\begin{equation}\label{eq:fixed_point_max_cut}
	r^*(\Delta) = \inf\left(r>0 : \bP\left[\sup_{Z\in\cC: \inr{\E B,Z^*-Z}\leq r}\inr{B-\E B, Z-Z^*}\leq (1/2) r\right]\geq 1-\Delta\right).
	\end{equation}Our second result is an explicit high probability upper bound on the latter fixed point.

	\subsubsection{Main results for the MAX-CUT problem}
	In this section, we gather the two results on the estimation of $Z^*$ from $\hat Z$ and on the approximate optimality of the final cut constructed from $\hat Z$. Let us now explicitly provide the construction of this cut.  We consider the same strategy as in \cite{goemans1995improved}. Assume that $\hat Z$ has been constructed. Let $\hat G$ be a centered Gaussian vector with covariance matrix $\hat Z$. Let $\hat x$ be the sign vector of $\hat G$. Using the statistical properties of $\hat Z$, it is possible to prove near optimality of $\hat x$. 
	
	We denote the optimal values of the maxcut problem and its SDP relaxation by
	\begin{equation*}
	\mathrm{SDP}(G) := (1/4)\inr{A^0,J-Z^*} = \max_{Z\in\cC}\frac{1}{4}\sum_{i,j}A_{i,j}^0(1-Z_{ij}) \mbox{ and } \mathrm{MAXCUT}(G) := \mathrm{cut}(G, S^*) 
	\end{equation*}
	where $S^*$ is a solution of \eqref{pb:MAX-CUT0} and $J=(1)_{n\times n}$. Our first result is to show how the $0.878$ approximating result from \cite{goemans1995improved} is downgraded by the incomplete information we have on the graph (we only partially observed the adjacency matrix $A^0$ via $A$).
	
	\begin{Theorem}\label{theo:goemans_williamson}
		For all $0<\Delta<1$. With probability at least $1-\Delta$ (with respect to the masked $S$), 
		\begin{equation*}
		\mathrm{SDP}(G) \geq \E \left[\mathrm{cut}(G, \hat x)|\hat Z\right]\geq 0.878\mathrm{SDP}(G) - \frac{0.878 r^*(\Delta)}{4}.
		\end{equation*}
	\end{Theorem}
	To precise the notation, $\hat x$ is the sign vector of $\hat G$ which is a centered Gaussian variable with covariance $\hat Z$. In that context, $\E \left[\mathrm{cut}(G, \hat x)|\hat Z\right]$ is the conditional expectation according to $\hat G$ for a fixed $\hat Z$. Moreover, the probability at least $1-\Delta$ that we obtain is w.r.t. the mask that is to the randomness in $A$.
	
	Let us put Theorem~\ref{theo:goemans_williamson} into some perspective. If we had known the entire adjacency matrix (which is the case when $p=1$), then we could have use $Z^*$ instead of $\hat Z$. In that case, for $x^\star$ the sign vector of $G^\star\sim \cN(0, Z^*)$, we know from \cite{goemans1995improved} that 
	\begin{equation}\label{eq:goemans_williamson}
	\mathrm{SDP}(G) \geq \E \left[\mathrm{cut}(G, x^\star)\right]\geq 0.878\mathrm{SDP}(G).
	\end{equation}Therefore, Theorem~\ref{theo:goemans_williamson} characterizes the price we pay for not observing the entire adjacency matrix $A^0$ but only a masked version $A$ of it. It is an interesting output of Theorem~\ref{theo:goemans_williamson} to observe that the fixed point $r^*(\Delta)$ measures, in a quantitative way, this loss.  If we were able to identify scenarii of $p$ and $E$ for which $r^*(\Delta)=0$ that would prove that there is no loss for partially observing $A^0$ in the MAX-CUT problem. Unfortunately, the approach we use to control $r^*(\Delta)$ is the global one which does not allow for exact reconstruction (that is to show that  $r^*(\Delta)=0$).  
	
	Let us now turn to an estimation result of $Z^*$ by $\hat Z$ via an upper bound on $r^*(\Delta)$.
	
	\begin{Theorem}\label{theo:max_cut_1}
		With probability at least $1-4^{-n}$:
		\begin{align*}
		\inr{\E B, Z^*-\hat Z}\leq r^*(4^{-n})\leq 2n\sqrt{\frac{(2\log4)(1-p)(n-1)}{p}} + \frac{8n\log 4}{3}.
		\end{align*}
	\end{Theorem}
	
	In particular, it follows from the approximation result from Theorem~\ref{theo:goemans_williamson} and the high probability upper bound on $r^*(\Delta)$ from Theorem~\ref{theo:max_cut_1} that, with probability at least $1-4^{-n}$,
	\begin{equation}\label{eq:final_eq_max_cut}
	\E \left[\mathrm{cut}(G, \hat x)|\hat Z\right]\geq 0.878\mathrm{SDP}(G) - \frac{0.878}{4}\left(2n\sqrt{\frac{(2\log4)(1-p)(n-1)}{p}} + \frac{8n\log 4}{3}\right).
	\end{equation} This result is none trivial only when the right-hand side term is strictly larger than $(0.5)\mathrm{SDP}(G)$ which is the performance of a random cut. As a consequence, \eqref{eq:final_eq_max_cut} shows that one can still do better than randomness even in an incomplete information setup for the MAX-CUT problem when $p$, $n$ and $\mathrm{SDP}(G)$ are such that
	\begin{equation*}
	0.378 \mathrm{SDP}(G)>\frac{0.878}{4}\left(2n\sqrt{\frac{(2\log4)(1-p)(n-1)}{p}} + \frac{8n\log 4}{3}\right).
	\end{equation*}For instance, when $p$ is like a constant, it requires $\mathrm{SDP}(G)$ to be larger than $c_0n^{3/2}$ (for some absolute constant $c_0$) and when $p=1-1/n$, it requires $\mathrm{SDP}(G)$ to be at least $c_0 n$ (for some absolute constant $c_0$).

	% subsection a_model_for_max_cut_and_goemans_and_williamson_sdp_relaxation (end)
	% section application_to_max_cut (end)

	\section{Proof of Theorem~\ref{theo:main_signed_cluster} (signed clustering)} % (fold)
	\label{sec:proof_of_Theorem_theo:main_signed_cluster_in_signed_clustering}
	The aim of this paper is to put forward a methodology developed in learning theory for the study of SDP estimators. In each example, we follow this methodology. 
	For a problem, such as the signed clustering, where it is possible to characterize the curvature of the excess risk, we start to identify this curvature because the curvature function $G$, coming out of it, defines the local subsets of $\cC$ driving the complexity of the problem. Then, we turn to the stochastic part of the proof which is entirely summarized into the complexity fixed point $r^*_G(\Delta)$ from \eqref{eq:fixed_point_G}. Finally, we put the two pieces together and apply the main general result from Corollary~\ref{coro:main_coro} to obtain estimation results for the SDP estimator \eqref{eq:esti_signed_cluster} in the signed clustering problem; which is Theorem~\ref{theo:main_signed_cluster}.

	\subsection{Curvature equation}
	In this section, we show that the objective function $Z\in\cC \rightarrow \inr{Z, \bE A - \alpha J}$ satisfies a curvature assumption around its maximizer $Z^*$ with respect to the $\ell_1^{n\times n}$-norm given by $G(Z^*-Z) = \theta\norm{Z^*-Z}_1$ with parameter $\theta=\delta(p-q)$ (and margin exponent $\kappa= 1$).

	\begin{Proposition}\label{prop:curvature_signed_cluster} For $\theta = \delta(p-q)$, we have for all $Z\in\cC$, $\inr{\bE A - \alpha J,Z^* - Z} = \theta \norm{Z^* - Z}_1$.
	\end{Proposition}

	\begin{proof}
		Let $Z$ be in  $\cC$. We have
		\begin{align*}
		\inr{Z^* - Z, \bE A - \alpha J} & = \sum_{i, j = 1}^n (Z^*-Z)_{ij}(\bE A_{ij} - \alpha) \\
		& = \sum_{(i, j) \in \cC^+} (Z^*_{ij} - Z_{ij})(\delta(2p-1) - \alpha) + \sum_{(i, j) \in \cC^- }(Z^*_{ij}-Z_{ij})(\delta(2q-1)-\alpha) \\
		& = \delta(p-q) \left[ \sum_{(i, j) \in \cC^+} (Z^* - Z)_{ij} - \sum_{(i, j) \in \cC^-} (Z^* - Z)_{ij}\right].
		\end{align*} 
		Moreover, for all $(i, j) \in \cC^+, Z^*_{ij} = 1$  and  $0\leq Z_{ij} \leq 1$, so $(Z^*-Z)_{ij} = |(Z^*-Z)_{ij}|$. We also have for all $(i, j) \in \cC^-$, $(Z^*-Z)_{ij} = -Z_{ij} = - |(Z^*-Z)_{ij}|$ because in that case $Z^*_{ij} = 1$  and $0\leq Z_{ij} \leq 1$.  Hence, 
		\begin{equation*}
		\inr{Z^* - Z, \bE A - \alpha J} = \delta(p-q) \left[\sum_{(i, j) \in \cC^+} |(Z^*-Z)_{ij}| + \sum_{(i, j) \in \cC^-} |(Z^*-Z)_{ij}| \right] = \theta \norm{Z^*-Z}_1.
		\end{equation*}
	\end{proof}
	
	\subsection{Computation of the complexity fixed point $r^*_G(\Delta)$}
	Define $W := \Aobs - \bE \Aobs$ the noise matrix of the problem. Since $W$ is symmetric, its entries are not independent. In order to work only with independent random variables, we define the following matrix $\Psi \in \bR^{n\times n}$:
	\begin{equation}\label{def:Psi}
	\Psi_{ij} = \left\{\begin{array}{l}
	W_{ij} \mbox{ if } i \leq j \\
	0 \mbox{ otherwise,}
	\end{array} \right.
	\end{equation} where $0$ entries are considered as independent Bernoulli variables with parameter $0$ and therefore, $\Psi$ has independent entries, and satisfies the relation $W = \Psi + \Psi^\top$. 
	
	In order to obtain upper bounds on the fixed point complexity parameter $r^*_G(\Delta)$ associated with the signed clustering problem, we need to prove a high probability upper bound on the quantity
	\begin{equation} 
	\label{eq:qty_W_ZZstar}
	\sup_{Z\in\cC:\norm{Z-Z^*}_1\leq r}\inr{W, Z-Z^*}, 
	\end{equation}
	and then find a radius $r$ as small as possible such that the quantity in \eqref{eq:qty_W_ZZstar} is less than $(\theta/2)r$. We denote $\cC_r := \cC\cap(Z^*+r B_1^{n\times n})= \left\{ Z \in \cC: \norm{Z-Z^*}_1\leq r \right\}$ where $B_1^{n\times n}$ is the unit $\ell_1^{n\times n}$-ball of $\bR^{n\times n}$.
	
	We follow the strategy from \cite{MR3901009} by decomposing the inner product $\inr{W, Z-Z^*}$ into two parts according to the SVD of $Z^*$. This observation is a key point in the work of \cite{MR3901009} compared to the analysis from  \cite{guedon2016community}. This allows to perform the localization argument efficiently. Up to a change of index of the nodes, $Z^*$ is a block matrix with $K$ diagonal blocks of $1$'s. It therefore admits $K$ singular vectors $U_{\bullet k}:=I(i\in\cC_k)/\sqrt{|\cC_k|}$ with multiplicity $l_k$ associated with the singular value $l_k$ for all $k\in[K]$. We can therefore write 
	\begin{equation*}
	Z^* = \sum_{k=1}^K l_k U_{\bullet k} \otimes U_{\bullet k} = UDU^\top,
	\end{equation*} where $U\in\bR^{n\times K}$ has $K$ column vectors given by $U_{\bullet k}, k=1, \ldots, K$ and $D={\rm diag}(l_1, \ldots, l_K)$. We define the following projection operator
	\begin{equation*}
	\cP:M \in \bR^{n \times n} \rightarrow UU^TM + MUU^T - UU^TMUU^T
	\end{equation*}and its orthogonal projection $\cP^\perp$ by
	\begin{equation*}
	\cP^\perp:M \in \bR^{n \times n} \rightarrow M-\proj{M} = (\mathrm{I}_n - UU^T)M(\mathrm{I}_n - UU^T) = \sum_{k=K+1}^n \inr{M, U_{\bullet k}\otimes U_{\bullet k}}U_{\bullet k}\otimes U_{\bullet k}
	\end{equation*}where $U_{\bullet k}\in\bR^n, k=K+1, \ldots, n$ are such that $(U_{\bullet k}:k=1, \ldots, n)$ is an orthonormal basis of $\bR^n$.
	
	We use the same decomposition as in \cite{MR3901009}: for all $Z\in\cC$,
	\begin{align*}
	\inr{W, Z-Z^*} = \inr{W, \proj{Z-Z^*} + \porth{Z-Z^*}} = \underset{S_1(Z)}{\underbrace{\inr{\proj{Z-Z^*}, W}}} + \underset{S_2(Z)}{\underbrace{\inr{\porth{Z-Z^*}, W}}}.
	\end{align*}The next step is to control with large probability the two terms $S_1(Z)$ and $S_2(Z)$ uniformly for all $Z\in\cC\cap(Z^*+r B_1^{n\times n})$. To that end we use the two following propositions where we recall that $\rho = \delta\max(1-\delta(2p-1)^2, 1-\delta(2q-1)^2)$ and $\nu =\max(2p-1,1- 2q)$. The proof of Proposition~\ref{prop:control_S1} and \ref{prop:S2_dense_case} can be found in Section~\ref{sub:annex_2_in_signed_clustering}, it is based on \cite{MR3901009}.
	
	\begin{Proposition}\label{prop:control_S1}
		There are absolute positive constants $c_0,c_1,c_2$ and $c_3$ such that the following holds. If $\floorsup{c_1 r K/n}\geq 2eKn \exp(-(9/32)n\rho/K)$ then we have
		\begin{equation*} 
		\bP\left[\sup_{Z\in\cC\cap (Z^*+r B_1^{n\times n})} S_1(Z)\leq c_2  r \sqrt{\frac{K\rho}{n}\log\left(\frac{2e K n}{\floorsup{\frac{c_1 r K}{n}}}\right) }\right]\geq 1-3\left(\frac{\floorsup{\frac{c_1 r K}{n}}}{2eKn}\right)^{\floorsup{\frac{c_1 r K}{n}}}.
		\end{equation*}
	\end{Proposition}
	
	\begin{Proposition}\label{prop:S2_dense_case}There exists an absolute constant $c_0>0$ such that the following holds. When $n\nu \delta\geq \log n$, with probability at least $1-\exp(-\delta\nu n)$,
		\begin{equation*}
		\sup_{Z\in\cC\cap(Z^*+rB_1^{n\times n})} S_2(Z)\leq c_0 Kr \sqrt{\frac{\delta \nu}{n}}.
		\end{equation*}
	\end{Proposition}
	It follows from Proposition~\ref{prop:control_S1} and Proposition~\ref{prop:S2_dense_case} that when $n\nu \delta\geq \log n$, for all $r$ such that $\floorsup{c_1 r K/n}\geq 2eKn \exp(-(9/32)n\rho/K)$ we have, for $\Delta=\Delta(r) := \exp(-\delta\nu n) - 3\left(\floorsup{\frac{c_1 r K}{n}}/(2eKn)\right)^{\floorsup{\frac{c_1 r K}{n}}}$, with probability at least $1-\Delta$,
	\begin{equation*}
	\sup_{Z\in\cC\cap(Z^*+rB_1^{n\times n})} \inr{W, Z-Z^*}\leq c_0 Kr \sqrt{\frac{\delta \nu}{n}} + c_2  r \sqrt{\frac{K\rho}{n}\log\left(\frac{2e K n}{\floorsup{\frac{c_1 r K}{n}}}\right) }.
	\end{equation*} Moreover, we have
	\begin{equation}\label{eq:fixed_point_eq_signed_clustering}
	c_0 Kr \sqrt{\frac{\delta \nu}{n}} + c_2  r \sqrt{\frac{K\rho}{n}\log\left(\frac{2e K n}{\floorsup{\frac{c_1 r K}{n}}}\right) }\leq \frac{\theta}{2}r
	\end{equation}for $\theta = \delta(p-q)$ when $K\sqrt{\nu}\lesssim\sqrt{n\delta}(p-q)$ and $\floorsup{c_1 r K/n}\geq 2eKn \exp(-\theta^2n /(K\rho))$. In particular, when $(p-q)^2 n\delta \geq K^2\nu$ and $1\geq 2eKn \max\big(\exp(-\theta^2n /(K\rho)), \exp(-(9/32)n\rho/K)\big)$, we conclude that for all $0<r\leq n/(c_1 K)$ \eqref{eq:fixed_point_eq_signed_clustering} is true. Therefore, one can take $r_G^*(\Delta(0))=0$ meaning that we have exact reconstruction of $Z^*$:  if $(p-q)^2 n\delta \geq K^2\nu$ and $n\gtrsim K \max(\rho/\theta^2 \log(2eK^2\rho/\theta^2), (1/\rho)\log(2eK^2/\rho))$ then with probability at least $1-\exp(-\delta \nu n) - 3/(2eKn)$, $\hat Z = Z^*$.

	If $(p-q)^2 n\delta \geq K^2\nu$ and  $1< 2eKn \max\big(\exp(-\theta^2n /(K\rho)), \exp(-(9/32)n\rho/K)\big)$ then we do not have  exact reconstruction anymore and we have to take $r$ such that $c_1 r K/n \geq1$. In that case, \eqref{eq:fixed_point_eq_signed_clustering} holds when $r=n/(c_1 K)$ and so $r_G^*(\Delta(n/(c_1 K)))\leq n/(c_1 K)$. Therefore, it follows from Corollary~\ref{coro:main_coro} that 
	\begin{equation*}
	\norm{Z^*-\hat Z}_1\leq \frac{n}{c_1\delta(p-q)K}.
	\end{equation*}

	\section{Proofs of Theorem \ref{theo:main_synchro} and Corollary~\ref{cor:sync} (angular synchronization)} % (fold)
	\label{sec:proof_of_Theorem_in_synchronization}
	\textbf{Proof of \eqref{eq:sync_x_star}:}
	For all $\gamma_1, \ldots, \gamma_n \in [0, 2\pi)$ we have $\gamma_i - \gamma_j = \delta_{ij}$ for all $i\neq j\in[n]$ if and only if $e^{\sigma^2/2}\E A_{ij}e^{\iota \gamma_j} - e^{\iota \gamma_i} = 0$ for all $i \neq j\in[n]$. We therefore have
	\begin{equation}\label{eq:sync_optimal_vector_x_star}
	\argmin_{x \in \bC^n:|x_i| = 1} \left\{ \sum_{i\neq j} |e^{\sigma^2/2}\E A_{ij}x_j - x_i |^2 \right\} = \{(e^{\iota (\theta_i+\theta_0)})_{i=1}^n:\theta_0\in[0,2\pi)\}.
	\end{equation}Moreover, for all $x=(x_i)_{i=1}^n\in \mathbb{C}^n$ such that $|x_i|=1$ for $i=1, \ldots, n$, we have 
	\begin{align*}
	&\sum_{i\neq j } |e^{\sigma^2/2}\E A_{ij}x_j - x_i|^2  = \sum_{i\neq j} |x^*_i\bar x_j^*x_j - x_i|^2 =\sum_{i,j=1}^n |x_i^* \bar x_j^* - x_i \bar x_j|^2 \\
	& = 2n^2 - 2 \Re\left(\sum_{i,j=1}^n x_i^*\bar x_j^*x_j \bar x_i\right)= 2n^2-2|\inr{x^*, x}|^2
	\end{align*}where $\Re(z)$ denotes the real part of $z\in\bC$.  On the other side, we have
	\begin{align*}
	&\bar x^\top (e^{\sigma^2/2}\E A)x = \sum_{i\neq j} \bar x_i x_i^* \bar x_j^* x_j + \sum_{i=1}^n \bar x_i e^{\sigma^2/2} x_i = n(e^{\sigma^2/2} - 1) + |\inr{x^*, x}|^2.
	\end{align*}Hence,  minimizing $x\to\sum_{i\neq j}^n |e^{\sigma^2/2}\E A_{ij}x_j - x_i|^2$ over all $x=(x_i)_i\in\bC^n$ such that $|x_i|=1$ is equivalent to maximize $x\to\bar{x}^\top\E A x$ over all $x=(x_i)_i\in\bC^n$ such that $|x_i|=1$. This concludes the proof with \eqref{eq:sync_optimal_vector_x_star}.\endproof
	
	\noindent\textbf{Proof of \eqref{eq:sync_X_star}:} Let $\cC^\prime = \{ Z \in \mathbb{C}^{n\times n}: |Z_{ij}| \leq 1, \forall i, j\in[n]\}$. We first prove that $\cC\subset \cC^\prime$.  Let $Z \in \cC$. Since $Z \succeq 0$, there exists $X \in \bC^{n\times n}$ such that $Z = X \bar{X}^\top$. For all $i \in \{1, \ldots, n\}$, denote by $X_{i\bullet}$ the $i$-th row vector of $X$. We have $\norm{X_{i\bullet}}_2^2 = \inr{X_{i\bullet}, X_{i\bullet}} = Z_{ii} =  1$ since $\mathrm{diag}(Z) = \mathbf{1}_n$. Moreover, for all $i,j\in[n]$, we have $|Z_{ij}|=|\inr{X_{i\bullet}, X_{j\bullet}}|\leq\norm{X_{i\bullet}}_2 \norm{X_{j\bullet}}_2\leq1$. This proves that $Z\in\cC^\prime$ and so $\cC\subset \cC^\prime$.

	Let $ Z^\prime \in \argmax\left( \Re(\inr{\E A, Z}): Z \in\cC^\prime\right\}$. Since $\cC^\prime$ is convex and  the objective function $Z \rightarrow \Re(\inr{\E A, Z})$ is linear (for real coefficients), $Z^\prime$ is one of the extreme points of $\cC^\prime$. Extreme points of $\cC^\prime$ are matrices $Z \in \bC^{n\times n}$ such that $|Z_{ij}| = 1$ for all $i, j\in[n]$. We can then write each entry of $Z^\prime$ as $Z^\prime_{ij} = e^{\iota \beta_{ij}}$ for some $0 \leq \beta_{ij} < 2\pi$ and now we obtain
	\begin{align*}
	&\Re(\inr{\E A, Z^\prime})  = \Re{\left(\sum_{i, j = 1}^n \E A_{ij} \overline{Z^\prime_{ij}}\right)} = \Re{\left(\sum_{i\neq  j }^n e^{-\sigma^2/2}e^{\iota \delta_{ij}} e^{-\iota \beta_{ij}}\right)} + \Re{\left(\sum_{i=1}^n e^{\iota \delta_{ii}} e^{-\iota \beta_{ii}}\right)}\\
	& = \sum_{i\neq j} e^{\sigma^2/2}\cos(\delta_{ij} - \beta_{ij}) + \sum_{i} \cos(\delta_{ii} - \beta_{ii}) \leq e^{\sigma^2/2}(n^2-n)+n.
	\end{align*}
	The maximal value $e^{\sigma^2/2}(n^2-n) + n$ is attained only for  $\beta_{ij} = \delta_{ij}$ for all $i, j\in[n]$, that is for $Z^\prime = (e^{\iota \delta_{ij}})_{i, j=1, \ldots, n} = Z^*$. But we have $Z^* \in \cC$ and $\cC \subset \cC^\prime$, so $Z^*$ is the only maximizer of $Z\to\Re(\inr{\E A, Z} )$ on $\cC$. But for all $Z \in \cC$ we have $ \inr{\E A, Z} = \overline{x^*}^\top Z x^*  \in \bR$, then $Z^*$ is the only maximizer of   $Z\to  \inr{\E A, Z}$ over $\cC$.\endproof
	
	\subsection{Curvature of the objective function}\label{sec:sync:curvature}
	
	\begin{Proposition}\label{prop:curvature_synchro} For $\theta = e^{-\sigma^2/2}/2$, we have for all $Z\in\cC$, $\inr{\bE A, Z^* - Z} \geq \theta \norm{Z^* - Z}_2^2$.
	\end{Proposition}
	\begin{proof}
		Let  $Z = (z_{ij}e^{\iota \beta_{ij}})_{i, j=1}^n) \in \cC$ where $z_{ij}\in\bR$ and $0 \leq \beta_{ij} < 2\pi$ for all $i,j\in[n]$. Since $Z_{ii}^* = Z_{ii}=1$ for all $i\in[n]$, we have, on one side, $\inr{\E A, Z^* - Z} = e^{-\sigma^2/2}\overline{x^*}^\top (Z^*-Z)x^* \in \bR$, and so
		\begin{align}\label{eq:sync_curvature_1}
		\notag\inr{\E A, Z^* - Z}  & = \Re{\left( \inr{\E A, Z^* - Z} \right) }= \Re{\left( \sum_{i, j = 1}^n \E A_{ij}\overline{(Z^* - Z)_{ij}} \right)} = \Re{\left( \sum_{i, j = 1}^n e^{-\sigma^2/2}e^{\iota \delta_{ij}} (e^{-\iota \delta_{ij}} - z_{ij}e^{-\iota \beta_{ij}}) \right) } \\
		& = e^{-\sigma^2/2} \Re{\left( \sum_{i, j = 1}^n 1 - z_{ij}e^{\iota (\delta_{ij} - \beta_{ij})}\right) } = e^{-\sigma^2/2}\sum_{i, j = 1}^n (1 - z_{ij}\cos(\delta_{ij} - \beta_{ij})).
		\end{align}On the other side, we proved in the proof of  \eqref{eq:sync_X_star} that $\cC\subset \{ Z \in \mathbb{C}^{n\times n}: |Z_{ij}| \leq 1, \forall i, j\in[n]\}$. So we have $|z_{ij}|\leq1$ for all $i,j\in[n]$ and
		\begin{align}\label{eq:sync_curvature_2}
		\notag&\norm{Z^* - Z}_2^2  = \sum_{i, j = 1}^n |(Z^* - Z)_{ij}|^2 = \sum_{i, j = 1}^n |e^{\iota \delta_{ij}} - z_{ij}e^{\iota \beta_{ij}}|^2  =\sum_{i, j = 1}^n |1 - z_{ij}e^{\iota(-\delta_{ij} + \beta_{ij})}|^2 \\
		\notag    & =\sum_{i, j = 1}^n (1 - z_{ij} \cos(\beta_{ij} - \delta_{ij} ))^2 + z_{ij}^2 \sin^2(\beta_{ij} - \delta_{ij}) = \sum_{i, j = 1}^n 1 - 2z_{ij} \cos(\beta_{ij} - \delta_{ij}) + z_{ij}^2 \\
		& \leq 2 \sum_{i, j = 1}^n (1 - z_{ij}\cos(\beta_{ij} - \delta_{ij})).
		\end{align}We conclude with \eqref{eq:sync_curvature_1} and \eqref{eq:sync_curvature_2}.
	\end{proof}
	In fact, it follows from the proof of Proposition~\ref{prop:curvature_synchro} that we have the following equality: for all $Z\in\cC$,
	\begin{equation*}
	\inr{\bE A, Z^* - Z} = \theta \left(\norm{Z^* - Z}_2^2 + \norm{|Z^*|^2 - |Z|^2}_1 \right)
	\end{equation*}where $|Z|^2=(|Z_{ij}|^2)_{1\leq i,j\leq n}$ (in particular, $|Z^*|^2=(1)_{n\times n}$). We therefore know exactly how to characterize the curvature of the excess risk for the angular synchronization problem in terms of the $\ell_2$ (to the square) and the $\ell_1$ norms.  Nevertheless, we will not use the extra term $\norm{|Z^*|^2 - |Z|^2}_1$ in the following.

	\subsection{Computation of the complexity fixed point $r^*_G(\Delta)$}
	\label{sec:sync_fixed_point}
	It follows from the (global) curvature property of the excess risk for the angular synchronization problem obtained in Proposition~\ref{prop:curvature_synchro} that for the curvature $G$ function defined by $G(Z^*-Z)=\theta \norm{Z^*-Z}_2^2, \forall Z\in\cC$, we just have to compute the $r^*_G(\Delta)$ fixed point and then apply Corollary~\ref{coro:main_coro} in order to obtain statistical properties of $\hat Z$ (w.r.t. to both the excess risk and the $G$ function). In this section, we compute the complexity fixed point $r_G^*(\Delta)$ for $0<\Delta<1$.
	
	Following Proposition~\ref{prop:curvature_synchro},  the natural ``local'' subsets of $\cC$ around $Z^*$ which drive the statistical complexity of the synchronization problem are defined for all $r>0$ by $ \cC_r = \{Z \in \cC: \norm{Z-Z^*}_2\leq r\}=\cC\cap(Z^* + r B_2^{n\times n})$. 
	
	Let $Z\in\cC_r$. Denote by $b_{ij}^R$ (resp. $b_{ij}^I$) the real (resp. imaginary) part of  $b_{ij} = Z_{ij}^*\overline{Z_{ij} - Z_{ij}^*}$ for all $i,j\in[n]$. Since $\norm{Z-Z^*}_2\leq r$ we also have $\sum_{i,j}(b_{ij}^R)^2 + (b_{ij}^I)^2\leq r^2$ and so
	\begin{align*}
	&\inr{A-\bE A, Z-Z^*} = \inr{(S-\bE S)\circ Z^*, Z-Z^*} = 2 \Re\left(\sum_{i<j}(S_ij-\bE S_{ij})b_{ij}\right)\\
	& = 2 \sum_{i<j}(\cos(\sigma g_{ij}) - \E \cos(\sigma g_{ij}))b_{ij}^R - \sin(\sigma g_{ij}) b_{ij}^I\leq 2 r \sqrt{\sum_{i<j}(\cos(\sigma g_{ij}) - \E \cos(\sigma g_{ij}))^2 +  (\sin(\sigma g_{ij}))^2}\\
	&\leq 2r \sqrt{1-e^{-\sigma^2} + 2e^{-\sigma^2/2}\left(\sum_{i<j}\bE \cos(\sigma g_{ij})  - \cos(\sigma g_{ij})\right) }
	\end{align*}where we used that $\E \cos(\sigma g) = \Re(\bE e^{\iota g}) = e^{-\sigma^2/2}$  for $g\sim \cN(0,1)$. Now its remains to get a high probability upper bound on the sum of the centered cosinus of  $\sigma g_{ij}$. We use Bernstein's inequality (see Equation~\ref{eq:Bernstein} below) to get such a bound. For all $t>0$, with probability at least $1-\exp(-t)$,
	\begin{equation*}
	\frac{1}{\sqrt{N}}\sum_{i<j}\bE \cos(\sigma g_{ij})-\cos(\sigma g_{ij})\leq \sqrt{2Vt} + \frac{2t}{3\sqrt{N}}\leq (1-e^{-\sigma^2})\sqrt{t} + \frac{2t}{3\sqrt{N}}
	\end{equation*} for $N=n(n-1)/2$ and $V=\E \cos^2(\sigma g) - (\E \cos(\sigma g))^2 = (1/2)(1-e^{-\sigma^2})^2$ (because $\E \cos^2(\sigma g) =(1/2)\bE(1+\cos(2\sigma g)) = (1/2)(1+e^{-2\sigma^2})$ when $g\sim \cN(0,1)$).
	
	We now have all the ingredients to compute the fixed point $r_G^*(\Delta)$ for $0<\Delta<1$: for $\theta = e^{-\sigma^2/2}/2$ and $t=\log(1/\Delta)$,
	\begin{equation*}
	r^*_G(\Delta) \leq \frac{4}{\theta}\left(1-e^{-\sigma^2} + 2e^{-\sigma^2/2}\left((1-e^{-\sigma^2})\sqrt{tN} + \frac{2t}{3}\right)\right) = \frac{32t}{3} + 8(1 - e^{-\sigma^2})(e^{\sigma^2/2} + 2\sqrt{tN}).
	\end{equation*}In particular, using $1-e^{-\sigma^2}\leq \sigma^2$ and for $t=\eps\sigma^4 N$ (where $N=n(n-1)/2$) for some $0<\eps<1$, if $e^{\sigma^2/2}\leq 2\sigma^2\sqrt{\eps}N$ then $ r^*_G(\Delta)\leq (128/3)\sigma^4N\sqrt{\eps}$.

	%  Since $W$ is symmetric and $Z^*-Z$ is Hermitian, we have
	% \begin{align}\label{eq:ineq_621}
	% \inr{A-\bE A, Z-Z^*} &=  \sum_{i<j} W_{ij} \overline{(Z - Z^*)_{ij}} + \sum_{i>j} W_{ij} \overline{(Z - Z^*)_{ij}} = 2\sum_{i<j} g_{ij} \Re((Z - Z^*)_{ij}) = 2 \sum_{i<j} g_{ij} b_{ij}\in\bR
	% \end{align}
	%  where we set $b_{ij}=\Re((Z - Z^*)_{ij})$ for $1\leq i<j\leq n$.  It follows from Borell's concentration inequality for Gaussian processes \cite{MR1849347} that for all $u>0$, with probability at least $1-\exp(-u^2/2)$,
	% \begin{equation*}
	% \sup_{Z\in\cC_r}\inr{W, Z - Z^*}\leq \bE[\sup_{Z\in\cC_r}\inr{ W, Z - Z^*}] + v u
	% \end{equation*}where $v^2 = \sup_{Z\in\cC_r}\{\mathrm{Var}(\inr{W, Z - Z^*})\} = 4\sup_{Z\in\cC_r}\sum_{i<j} \Re((Z - Z^*)_{ij})^2 \leq 4r^2$. It also follows from \eqref{eq:ineq_621} that, for  $N=n(n-1)/2$,
	% \begin{align*}
	% \bE\sup_{Z\in\cC_r}\inr{W, Z - Z^*} \leq 2r \bE \sqrt{\sum_{i<j}g_{ij}^2}\leq 2r \sqrt{N}.
	% \end{align*}
	
	% We now have all the necessary tools to obtain an upper bound for $r^*_2(\Delta, \theta)$. Let $0<\Delta<1$. With probability at least $1-\Delta$, 
	% \begin{equation*}
	% \sup_{Z\in\cC_r}\inr{A-\bE A, Z-Z^*}\leq 2\sigma r\sqrt{N} + 2\sigma \sqrt{2}r \sqrt{\log(1/\Delta)}\leq (\theta/2)r^2
	% \end{equation*}for $\theta=1/2$ as long as $r\geq 8 \sigma(\sqrt{N} + \sqrt{2\log(1/\Delta)})$. Therefore, we have $r^*_2(\Delta, \theta)\leq 8 \sigma(\sqrt{N} + \sqrt{2\log(1/\Delta)})$.

	\subsection{End of the proof of Theorem \ref{theo:main_synchro} and Corollary~\ref{cor:sync}: application of Corollary~\ref{coro:main_coro}}
	Take $\Delta = \exp(-\eps\sigma^4N)$ (for $N=n(n-1)/2$), we have $r^*_G(\Delta)\leq (128/3)\sqrt{\eps} \sigma^4N$ when $e^{\sigma^2/2}\leq 2\sqrt{\eps}\sigma^2N$ (which holds for instance when $\sigma\leq \sqrt{\log(\eps N^2)}$)  and so it follows from Corollary~\ref{coro:main_coro} (together with the curvature property in Section~\ref{sec:sync:curvature} and the computation of the fixed point $r^*_G(\Delta)$ from Section~\ref{sec:sync_fixed_point}), that with probability at least $1-\exp(-\eps\sigma^4n(n-1)/2)$, $\theta\inr{Z^*-Z}_2^2\leq \inr{\E A, Z^*-Z} \leq (128/3)\sqrt{\eps}\sigma^4 N$, which is the statement of Theorem \ref{theo:main_synchro}.
	
	\textbf{Proof of Corollary~\ref{cor:sync}:}
	The oracle $Z^*$ is the rank one matrix $x^*\overline{x^*}^\top$ which has $n$ for largest eigenvalue and associated eigenspace $\{\lambda x^*:\lambda\in\bC\}$. In particular, $Z^*$ has a spectral gap $g=n$. Let $\hat{x} \in \bC^n$ be a top eigenvector of $\hat{Z}$ with norm $\norm{\hat{x}}_2=\sqrt{n}$. 
	%Now, we can see $\hat{X}$ as a perturbation of $x^*\bar{x^*}^T$:
	%\begin{equation*}
	%\hat{X} = x^*\bar{x^*}^T + (\hat{X} - x^*\bar{x^*}^T)
	%\end{equation*}
	It follows from Davis-Kahan Theorem (see, for example, Theorem~4.5.5 in \cite{vershynin2018high} or Theorem~4 in \cite{VanVu}) that there exists an universal constant $c_0>0$ such that
	\begin{equation*}
	\min_{z\in\bC:|z|=1}\norm{\frac{\hat{x}}{\sqrt{n}} - z\frac{x^*}{\sqrt{n}}}_2 \leq \frac{c_0}{g}\norm{\hat Z - Z^*}_2
	\end{equation*}where $g=n$  is the spectral gap of $Z^*$. We conclude the proof of Corollary~\ref{cor:sync} using the upper bound on $\norm{\hat Z - Z^*}_2$ from Theorem \ref{theo:main_synchro}. \endproof
	
	\section{Proofs of Theorem~\ref{theo:goemans_williamson} and \ref{theo:max_cut_1} (MAX-CUT)} % (fold)
	\label{sec:proofs_for_the_max_cut_problem_}
	In this section, we prove the two main results from Section~\ref{sec:application_to_max_cut} using our general methodology for Theorem~\ref{theo:max_cut_1} and the technique from \cite{goemans1995improved} for Theorem~\ref{theo:goemans_williamson}. 
	
	\subsection{Proof of Theorem~\ref{theo:goemans_williamson}} % (fold)
	\label{sub:proof_of_theorem_theo:goemans_williamson}
	The proof of Theorem~\ref{theo:goemans_williamson} follows the one from \cite{goemans1995improved} up to a minor modification due to the fact that we use the SDP estimator $\hat Z$ instead of the oracle $Z^*$. It is based on two tools. The first one is Grothendieck's identity:  let $g\sim \cN(0, \mathrm{I}_n)$ and $u, v \in \cS_2^{n-1}$, we have:
	\begin{equation}\label{eq:grothendieck_identity}
	\bE[\mathrm{sign}(\inr{g, u})\mathrm{sign}(\inr{g, u})] = \frac{2}{\pi}\mathrm{arcsin}(\inr{u, v})
	\end{equation}and the identity: for all $t \in [-1, 1]$
	\begin{equation}\label{eq:identity_arcsin}
	1-\frac{2}{\pi}\mathrm{arcsin}(t) = \frac{2}{\pi}\mathrm{arccos}(t)\geq 0.878(1-t).
	\end{equation}
	
	We now have enough tools to prove Theorem~\ref{theo:goemans_williamson}.  The right-hand side inequality is trivial since $\mathrm{MAXCUT}(G) \leq \mathrm{SDP}(G)$. For the left-hand side, we denote by $\hat X_1, \ldots, \hat X_n$ (resp. $X_1^*, \ldots, X_n^*$) the $n$ columns vectors in $\cS_2^{n-1}$ of $\hat Z$ (resp. $Z^*$). We also consider the event $\Omega^*$ onto which
	\begin{equation*}
	\inr{\E B, Z^*-\hat Z}\leq r^*(\Delta)
	\end{equation*}which hold with probability at least $1-\Delta$ according to Theorem~\ref{theo:main}. On the event $\Omega^*$, we have
	\begin{align*}
	&\bE\left[\mathrm{cut}(G, \hat x)|\hat Z\right]  = \bE\left[\frac{1}{4}\sum_{i, j=1}^nA^0_{ij}(1-\hat x_i \hat x_j)\right] = \frac{1}{4}\sum_{i, j=1}^nA^0_{ij}\left(1-\bE[\mathrm{sign}(\inr{\hat X_i, g})\mathrm{sign}(\inr{\hat X_j, g})]\right) \\
	& \overset{(i)}{=} \frac{1}{4} \sum_{i, j=1}^nA^0_{ij} \left(1-\frac{2}{\pi} \arcsin(\inr{\hat X_i, \hat X_j})\right) =  \frac{1}{2\pi} \sum_{i, j=1}^nA^0_{ij} \arccos(\inr{\hat X_i, \hat X_j}) \\
	& \overset{(ii)}{\geq}\frac{0.878}{4} \sum_{i, j=1}^nA^0_{ij}(1-\inr{\hat X_i, \hat X_j})  = \frac{0.878}{4} \sum_{i, j=1}^nA^0_{ij}(1-\inr{X_i^*, X_j^*}) + \frac{0.878}{4} \sum_{i, j=1}^nA^0_{ij}(\inr{X_i^*, X_j^*} - \inr{\hat X_i, \hat X_j})\\
	& = \frac{0.878}{4} \inr{A^0,J-Z^*} + \frac{0.878}{4} \inr{A^0, Z^*-\hat Z} = 0.878\mathrm{SDP}(G) - \frac{0.878}{4} \inr{\E B, Z^*-\hat Z}\\
	&\geq 0.878\ \mathrm{SDP}(G) - \frac{0.878}{4}r^*(\Delta)
	\end{align*}where we used \eqref{eq:grothendieck_identity} in \textit{(i)} and \eqref{eq:identity_arcsin} in \textit{(ii)}.

	\subsection{Proof of Theorem \ref{theo:max_cut_1}}
	% subsection proof_of_theorem_theo:goemans_williamson (end)
	
	For the MAX-CUT problem, we do not use any localization argument; we therefore use the (likely sub-optimal) global approach. The methodology is very close to the one used in \cite{guedon2016community} for the community detection problem. In particular, we use both Bernstein and Grothendieck inequalities to compute high probability upper bound for $r^*(\Delta)$. We recall theses two tools now. First Bernstein's inequality: if $Y_1,\ldots, Y_N$ are $N$ independent centered random variables such that $|Y_i|\leq M$ a.s. for all $i=1, \ldots, N$ then for all $t>0$, with probability at least $1-\exp(-t)$,
	\begin{equation}\label{eq:Bernstein}
	\frac{1}{\sqrt{N}}\sum_{i=1}^N Y_i\leq \sigma \sqrt{2t} + \frac{2Mt}{3\sqrt{N}}
	\end{equation}  where $\sigma^2 = (1/N)\sum_{i=1}^N {\rm var}(Y_i)$. The second tool is Grothendieck inequality \cite{grothendieck1956resume} (see also \cite{pisier2012grothendieck} or Theorem~3.4 in \cite{guedon2016community}): if $C\in\bR^{n\times n}$ then
	\begin{equation}\label{eq:Grothendieck_ineq}
	\sup_{Z\in\cC}\inr{C, Z}\leq K_G \norm{C}_{cut} = K_G\max_{s,t\in\{-1,1\}^n}\sum_{i,j=1}^n C_{ij} s_i t_j
	\end{equation}where $\cC=\{Z\succeq0: Z_{ii} =  1, i=1, \ldots, n\}$ and $K_G$ is an absolute constant, called the Grothendieck constant.

	In order to apply Theorem~\ref{theo:main}, we just have to compute the fixed point $r^*(\Delta)$. As announced, we use the global approach and Grothendieck inequality \eqref{eq:Grothendieck_ineq} to get
	\begin{equation}\label{eq:global_approach_max_cut}
	\sup_{Z\in\cC:\inr{\E B, Z^*-Z}\leq r} \inr{B-\E B, Z-Z^*} \leq \sup_{Z\in\cC}\inr{B-\E B, Z-Z^*}\leq 2K_G \norm{B-\E B}_{cut}
	\end{equation}because $Z^*\in\cC$. It follows from Bernstein's inequality \eqref{eq:Bernstein} and a union bound that for all $t>0$, with probability at least $1-4^n\exp(-t)$,
	\begin{equation*}
	\norm{B-\E B}_{cut} = \sup_{s,t\in\{\pm1\}^n}\sum_{1\leq i<j\leq n} (B_{ij} - \E B_{ij})(s_it_j+s_jt_i)\leq 2\sqrt{\frac{(1-p)n(n-1)t}{p}} + \frac{4t}{3}.
	\end{equation*}Therefore, for $t=2n\log 4$,  with probability at least $1-4^{-n}$, 
	\begin{equation*}
	r^*(\Delta)\leq \norm{B-\E B}_{cut}\leq 2n\sqrt{\frac{(2\log4)(1-p)(n-1)}{p}} + \frac{8n\log 4}{3}
	\end{equation*}for  $\Delta = 4^{-n}$. Then the result follows from Theorem~\ref{theo:main}.

	\section{Annex~1: Signed clustering} % (fold)
	\label{sub:annex_2_in_signed_clustering}
	
	\subsection{Proof of Equation~\eqref{eq:oracle_egal_cluster_mat}} % (fold)
	\label{sub:eq_3.4}
	We recall that the cluster matrix $\bar Z \in \{0, 1\}^{n \times n}$ is defined by $Z_{ij} = 1$ if $i\sim j$ and $Z_{ij}=0$ when $i\not\sim j$ and $\alpha=\delta(p+q-1)$. For all matrix $Z \in [0,1]^{n \times n}$, we have
	\begin{align*}
	&\inr{Z, \bE A - \alpha J}  = \sum_{i, j = 1}^n Z_{ij}(\bE A_{ij} - \alpha)  = \sum_{(i, j) \in \cC^+} Z_{ij}(\bE A_{ij} - \alpha) + \sum_{(i, j) \in \cC^-} Z_{ij}(\bE A_{ij} - \alpha)  \\
	& = [\delta(2p-1)-\alpha] \sum_{(i, j) \in \cC^+:i\neq j} Y_{ij} + [\delta(2q-1)-\alpha]\sum_{(i, j) \in \cC^-} Z_{ij} + (1-\alpha)\sum_{i=1}^n Z_{ii} \\
	& = \delta(p-q) \left[ \sum_{(i, j) \in \cC^+} Y_{ij} - \sum_{(i, j) \in \cC^-} Y_{ij} \right] + (1-\alpha)\sum_{i=1}^n Z_{ii}.
	\end{align*}
	The latter quantity is maximal for  $Z \in [0,1]^{n \times n}$ such that $Z_{ij} = 1$ for $(i, j) \in \cC^+$ and $Z_{ij} = 0$ for $(i, j) \in \cC^-$, that is when $Z = \bar{Z}$. As a consequence $\{\bar Z\}=\argmax_{Z\in[0,1]^{n\times n}} \inr{Z, \bE A - \alpha J}$. Moreover, $\bar Z\in \cC\subset[0,1]^{n\times n}$ so we also have that $\bar{Z}$ is the only solution to the problem $\max_{Z\in\cC} \inr{Z, \bE A - \alpha J}$ and so $\bar Z = Z^*$.
	\endproof

	\subsection{Proof of Proposition~\ref{prop:control_S1}: control of $S_1(Z)$ adapted from \cite{MR3901009}}
	The noise matrix $W$ is symmetric and has been decomposed as $W=\Psi + \Psi^\top$ where $\Psi$ has been  defined in \eqref{def:Psi}. For all $Z\in\cC\cap (Z^*+r B_1^{n\times n})$, we have
	\begin{align}\label{eq:S1_begin}
	\notag&S_1(Z)  = \inr{\cP(Z-Z^*), W} = \inr{\proj{W}, Z - Z^*}\\
	\notag& = \inr{UU^\top W, Z-Z^*} + \inr{WUU^\top, Z-Z^*} - \inr{UU^\top WUU^\top, Z-Z^*} \\
	\notag& = 2 \inr{UU^\top W, Z-Z^*} - \inr{UU^\top WUU^\top, Z-Z^*} \\
	\notag& = 2 \inr{UU^\top\Psi, Z-Z^*}  + 2 \inr{UU^\top \Psi^\top, Z-Z^*} - \inr{UU^\top(\Psi + \Psi^\top)UU^\top, Z-Z^*} \\
	\notag& = 2 \inr{UU^\top \Psi, Z-Z^*}  + 2 \inr{UU^\top\Psi^\top, Z-Z^*} - 2\inr{UU^\top\Psi UU^\top, Z-Z^*} \\
	& = 2 \inr{UU^\top\Psi, Z-Z^*}  + 2 \inr{UU^\top\Psi^\top, Z-Z^*} - 2\inr{UU^\top\Psi, (Z-Z^*)UU^\top}.
	\end{align}An upper bound on $S_1(Z)$ follows from an upper bound on the three terms in the right side of \eqref{eq:S1_begin}. Let us show how to bound the first term. Similar arguments can be used to control the other two terms.
	
	Let $V:=UU^T\Psi$. Let us find a high probability upper bound on the term  $\inr{UU^\top\Psi, Z-Z^*} = \inr{V, Z-Z^*}$ uniformly over $Z\in\cC\cap (Z^*+r B_1^{n\times n})$.  For all $k\in [K]$, $i\in\cC_k$ and $j\in[n]$, we have 
	\begin{equation*}
	V_{ij} = \sum_{t=1}^n(UU^T)_{it} \Psi_{tj} = \sum_{t\in \cC_k} \frac{1}{l_k} \Psi_{tj} = \frac{1}{l_k} \sum_{t\in \cC_k}  \Psi_{tj} = \frac{1}{l_k} \sum_{t\in \cC_k}  \Psi_{tj}.
	\end{equation*}Therefore, given $j\in[n]$ the $V_{ij}$'s are all equal for $i \in \cC_k$. We can therefore fix some arbitrary index $i_k \in \cC_k$ and have $V_{ij}=V_{i_k j}$ for all $i \in \cC_k$. Moreover, $(V_{i_k j}:k\in[K], j\in[n])$ is a family of independent random variables. We now have
	\begin{align*}
	\inr{V, Z-Z^*} = \sum_{k\in [K]} \sum_{i\in \cC_k} \sum_{j\in [n]} V_{ij}(Z-Z^*)_{ij}= \sum_{k\in [K]} \sum_{j\in[n]}l_k V_{i_kj} \sum_{i\in \cC_k}\frac{(Z-Z^*)_{ij}}{l_k}=\sum_{k\in [K]} \sum_{j\in[n]}l_k V_{i_kj}w_{kj}
	\end{align*} which is a weighted sum of $nK$ independent centered random variables $X_{k,j}:=l_k V_{i_kj}$ with weights $w_{k,j} = (1/l_k)\sum_{i\in \cC_k}(Z-Z^*)_{ij}$ for $ k \in [K],j \in [n]$. We now idenfity some properties on the weights $w_{kj}$.

	The weights are such that  
	\begin{align*}
	\sum_{k \in [K]}\sum_{j \in [n]} |w_{k,j}| \leq \frac{c_1 K}{n} \sum_{k \in [K]}\sum_{j \in [n]} \sum_{i \in \cC_k} |(Z-Z^*)_{ij}| = \norm{Z-Z^*}_{1} \frac{c_1K}{n} \leq \frac{c_1rK}{n}
	\end{align*}which is equivalent to say that the weights vector $w=(w_{kj}:k\in[K], j\in[n])$ is in the $\ell_1^{Kn}$-ball $(c_1rK/n)B_1^{Kn}$. It is also in the unit $\ell_\infty^{Kn}$-ball since for all $k \in [K]$ and $j\in[n]$, 
	\begin{equation*}
	|w_{k,j}|\leq \sum_{i\in \cC_k} \frac{|(Z-Z^*)_{ij}|}{l_k} \leq \norm{Z-Z^*}_{\infty} \leq 1.
	\end{equation*}We therefore have $w\in B_\infty^{Kn}\cap (c_1rK/n)B_1^{Kn}$ and so
	\begin{equation}\label{eq:inter_S1_rearrangement}
	\sup_{Z\in \cC\cap (Z^*+r B_1^{n\times n})}\inr{V,Z-Z^*}\leq \sup_{w\in B_\infty^{Kn}\cap (c_1rK/n)B_1^{Kn}}\sum_{k\in[K], j\in[n]}X_{k,j}w_{k,j}.
	\end{equation}It remains to find a high probability upper bound on the latter term. We can use the following lemma to that end.
	
	\begin{Lemma}\label{lem:weighted_sum_Bernstein_var}
		Let $X_{k,j}=\sum_{t\in\cC_k}\Psi_{tj}$ for $(k,j)\in[K]\times [n]$. For all $0\leq R\leq Kn$, if $\floorsup{R}\geq 2eKn \exp(-(9/32)n\rho/K)$ then with probability at least $1-\left(\floorsup{R}/(2eKn)\right)^{\floorsup{R}}$,
		\begin{equation*}
		\sup_{w\in B_\infty^{Kn}\cap R B_1^{Kn}}\sum_{(k,j)\in[K]\times [n]} X_{k,j} w_{k,j}\leq 4\sqrt{8 c_0}R \sqrt{\frac{n\rho}{K}\log\left(\frac{2eKn}{\floorsup{R}}\right)}.
		\end{equation*} 
	\end{Lemma}
	\textbf{Proof of Lemma~\ref{lem:weighted_sum_Bernstein_var}.} Let $N=Kn$ and assume that $1\leq R \leq N$. We denote by $X_1^*\geq X_2^*\geq \cdots, \geq X_N^*$ (resp. $w_1^*\geq \cdots\geq w_N^*$) the non-decreasing rearrangement of $|X_{k,j}|$ (resp. $|w_{k,j}|$) for $(k,j)\in[K]\times [n]$. We have  
	\begin{align*}
	&\sup_{w\in B_\infty^N\cap R B_1^N}\sum_{(k,j)\in[K]\times [n]} X_{k,j} w_{k,j}\leq \sup_{w\in B_\infty^N\cap R B_1^N}\sum_{i=1}^N X_i^* w_i^*\leq \sup_{w\in B_\infty^N}\sum_{i=1}^{\lceil R \rceil} X_i^*w_i^* + \sup_{w\in RB_1^N}\sum_{i=\lceil R \rceil+1}^N X_i^*w_i^*\\
	&\leq  \sum_{i=1}^{\lceil R \rceil} X_i^* + R X_{\lceil R \rceil+1}^*\leq 2\sum_{i=1}^{\lceil R \rceil} X_i^*.
	\end{align*}
	
	Moreover, for all $\tau>0$, using a union bound, we have
	\begin{align*}
	&\bP\left( \sum_{i=1}^{\floorsup{R}} X_i^* > \tau \right)  = \bP\left( \exists I \subset [K]\times[n]:|I| = \floorsup{R} \mbox{ and } \sum_{(k,j)\in I}|X_{k,j}| > \tau \right)\\
	& = \bP\left( \max_{I\subset[K]\times[n]:|I| = \floorsup{R}}\max_{u_{k,j}=\pm1, (k,j)\in I} \sum_{(k,j)\in I} u_{k,j} X_{k,j} > \tau \right)\\
	& \leq \sum_{I\subset[K]\times[n]:|I| = \floorsup{R}} \sum_{u \in \{ \pm 1 \}^{\floorsup{R}}} \bP\left( \sum_{(k,j)\in I}X_{k,j}u_{k,j} > \tau \right) = \sum_{I\subset[K]\times[n]:|I| = \floorsup{R}} \sum_{u \in \{ \pm 1 \}^{\floorsup{R}}} \bP\left( \sum_{(k,j)\in I}\sum_{t\in\cC_k}\Psi_{t,j}u_{k,j} > \tau \right).
	\end{align*} Let us now control each term of the latter sum thanks to Bernstein inequality. The random variables $(\Psi_{t,j}:t,j\in[n])$ are independent with variances at most $\rho = \delta\max(1-\delta(2p-1)^2, 1-\delta(2q-1)^2)$ since $\mathrm{Var}(\Psi_{ij}) =0$ when $i> j$ and $\mathrm{Var}(\Psi_{ij}) = \mathrm{Var}(\Aobs_{ij} - \bE[\Aobs_{ij}]) = \mathrm{Var}(\Aobs_{ij}) \leq \rho $ for $j\geq i$. Moreover, $|\Psi_{ij}| = 0$ when $j<i$ and $|\Psi_{ij}|=|W_{ij}| = |A_{ij} - \E A_{ij}|\leq 2$ for $j\geq i$ because $A_{ij}\in\{-1, 0, 1\}$. It follows from Bernstein's inequality that for all $I\subset[K]\times[n]$ satisfying $|I| = \floorsup{R}$ and $u \in \{ \pm 1 \}^{\floorsup{R}}$ that
	\begin{align*}
	\bP\left( \sum_{(k,j)\in I}\sum_{t\in\cC_k}\Psi_{t,j}u_{k,j} > \tau \right) \leq \exp\left( \frac{-\tau^2}{2\floorsup{R}l_k\rho + 4\tau/3}\right)\leq \exp\left( \frac{-\tau^2}{2\floorsup{R}c_0 n\rho/K + 4\tau/3}\right)\leq \exp\left( \frac{-\tau^2}{4\floorsup{R}c_0 n \rho/K}\right)
	\end{align*}when $\tau \leq (3/2)\floorsup{R}c_0 n \rho/K$. Therefore, $\sup_{w\in B_\infty^N\cap R B_1^N}\sum X_{k,j} w_{k,j}\leq 2\tau$  with probability at least 
	\begin{align*}
	1- \binom{N}{\floorsup{R}}2^{\floorsup{R}} \exp\left( \frac{-\tau^2}{4\floorsup{R}c_0 n\rho/K }\right)\geq 1-\exp\left( \frac{-\tau^2}{8\floorsup{R}c_0 n \rho/K}\right) 
	\end{align*}when
	\begin{equation*}
	(3/2)\floorsup{R}c_0 n \rho/K\geq \tau \geq\sqrt{8c_0}\floorsup{R}  \sqrt{\frac{n\rho}{K}\log\left(\frac{2eN}{\floorsup{R}}\right)}  
	\end{equation*}which is a non vacuous condition since $\floorsup{R}\geq 2eN \exp(-(9/32)n\rho/K)$. The result follows, in the case $1\leq R\leq N$, by taking $\tau = \sqrt{8c_0}\floorsup{R} \sqrt{(n\rho/K)\log\left(2eN/\floorsup{R}\right)}$ and using that $2R\geq\floorsup{R}$ when $R\geq1$.
	
	For $0\leq R \leq1$, we have
	\begin{equation*}
	\sup_{w\in B_\infty^{Kn}\cap R B_1^{Kn}}\sum_{(k,j)\in[K]\times [n]} X_{k,j} w_{k,j} = R \max_{(k,j)\in[K]\times [n]}|X_{k,j}|
	\end{equation*}and using Bernstein inequality as above we get that with probability at least $1-\exp(-K\tau^2/(8c_0 n\rho))$, $\max_{(k,j)\in[K]\times [n]}|X_{k,j}|\leq \tau$ when $3c_0n\rho/(2K)\geq \tau \geq \sqrt{8c_0 n \rho\log(nK)/K}$ which is a non vacuous condition when $9 c_0 n \rho \geq 4 K \log(nK)$. By taking $\tau = \sqrt{8c_0 n \rho\log(nK)/K}$, we obtain, that for all $0\leq R\leq 1$, if $9 c_0 n \rho \geq 4 K \log(nK)$ then with probability at least $1-1/(nK)$, 
	\begin{equation*}
	\sup_{w\in B_\infty^{Kn}\cap R B_1^{Kn}}\sum_{(k,j)\in[K]\times [n]} X_{k,j} w_{k,j} \leq R \sqrt{\frac{8c_0 n\rho \log(nK)}{K}}.
	\end{equation*}
	\endproof
	
	We apply Lemma~\ref{lem:weighted_sum_Bernstein_var} for $R=c_1 r K/n$ to control \eqref{eq:inter_S1_rearrangement}: 
	\begin{equation*}
	\bP\left[\sup_{Z\in\cC\cap (Z^*+r B_1^{n\times n})}\inr{V, Z-Z^*}\leq c_2  r \sqrt{\frac{K\rho}{n}\log\left(\frac{2e K n}{\floorsup{\frac{c_1 r K}{n}}}\right) }\right]\geq 1-\left(\frac{\floorsup{\frac{c_1 r K}{n}}}{2eKn}\right)^{\floorsup{\frac{c_1 r K}{n}}}
	\end{equation*}when $\floorsup{c_1 r K/n}\geq 2eKn \exp(-(9/32)n\rho/K)$.

	Using the same methodology, we can prove exactly the same result for $\sup_{Z\in\cC\cap (Z^*+r B_1^{n\times n})}\inr{UU^\top \Psi^\top, Z-Z^*}$. We can also use the same method to upper bound $\sup_{Z\in\cC\cap (Z^*+r B_1^{n\times n})}\inr{UU^\top \Psi^\top, (Z-Z^*)UU^\top}$, we simply have to check that the weights vector $w^\prime = (w^\prime_{kj}:k\in[K], j\in[n])$ where $w^\prime_{kj} = (1/l_k)\sum_{i\in\cC_k} [(Z-Z^*)UU^\top]_{ij}$ is also in $B_\infty^{Kn}\cap (c_1 r K/n)B_1^{Kn}$ for any $Z\in\cC$ such that $\norm{Z-Z^*}_1\leq r$. This is indeed the case, since we have for all $i\in[n]$, $k^\prime\in[K]$ and $j\in\cC_{k^\prime}$, $[(Z-Z^*)UU^\top]_{ij}=\sum_{p=1}^n (Z-Z^*)_{ip}(UU^\top)_{pj} = \sum_{p\in\cC_{k^\prime}}(Z-Z^*)_{ip}/l_{k^\prime}$ which is therefore constant for all elements in $j\in\cC_{k^\prime}$. Therefore, we have
	\begin{align*}
	&\sum_{k \in [K]}\sum_{j \in [n]} |w_{kj}^\prime| =  \sum_{k \in [K]}\sum_{k^\prime \in [K]}\sum_{j\in\cC_{k^\prime}} \left|\frac{1}{l_kl_{k^\prime}}\sum_{i \in \cC_k} \sum_{p\in\cC_{k^\prime}} (Z-Z^*)_{ip}\right|\\
	& \leq \sum_{k \in [K]}\sum_{k^\prime \in [K]}\sum_{j\in\cC_{k^\prime}} \frac{1}{l_kl_{k^\prime}}\sum_{i \in \cC_k} \sum_{p\in\cC_{k^\prime}} \left|(Z-Z^*)_{ip}\right| \leq \norm{Z-Z^*}_{1} \frac{c_1K}{n} \leq \frac{c_1rK}{n}
	\end{align*}and for all $k^\prime\in[K]$ and $j\in\cC_{k^\prime}$,
	\begin{equation*}
	|w_{kj}^\prime|=\left| \frac{1}{l_kl_{k^\prime}}\sum_{i\in \cC_k} \sum_{p\in\cC_{k^\prime}}(Z-Z^*)_{ip}\right| \leq \norm{Z-Z^*}_{\infty} \leq 1.
	\end{equation*}Therefore, $w^\prime\in B_\infty^{Kn}\cap (c_1 r K/n)B_1^{Kn}$ and we obtain exactly the same upper bound for the three terms in \eqref{eq:S1_begin}. This concludes the proof of Proposition~\ref{prop:control_S1}.

	% \noindent Plugging those inequalities into \ref{ineqS1}, we obtain that:
	% \begin{align*}
	% \bP\left(S_1(Y) \leq 6 D \floorsup{\beta}\sqrt{\rho l \log(\frac{nk}{\beta})} \right) > 1 - 3(\frac{\beta}{kn})^{\floorsup{\beta}}
	% \end{align*}
	% Under assumption \ref{betaBound}, we have $nke^{-\frac{\rho}{C_e l}} \leq \floorsup{\beta} \leq \floorsup{nkv} \leq nkv + 1 \leq \frac{nk}{3}$, where the last inequality holds for $v$ small enough.\\
	% \noindent Let us consider the function $f:x>0\rightarrow x\log(\frac{nk}{x}) - \log(nk)$. The derivative of $f$ if $f^{'}:x>0\rightarrow \log(nk)-\log(x)-1$. Then, for $x\leq \frac{nk}{3}$, we have $f^{'}(x) \geq \log(nk) - \log(\frac{nk}{3}) - 1 = \log(3) - 1 >0$, which means that $f^{'}$ is increasing on $[nke^{-\frac{\rho}{C_e l}} , \frac{nk}{3}]$. But $f(1) = 0$, then $f$ is non-negative on $[1, \frac{nk}{3}]$.
	% As $\floorsup{\beta}\geq 1$, we conclude that $f(\floorsup{\beta})\geq 0$, that is $(\frac{\floorsup{\beta}}{nk})^{\floorsup{\beta}} \leq \frac{1}{nk}$, and then $(\frac{\beta}{nk})^{\floorsup{\beta}} \leq (\frac{\floorsup{\beta}}{nk})^{\floorsup{\beta}} \leq \frac{1}{nk}$\\
	
	% \noindent Finally, we have:
	% \begin{align}\label{S1_result}
	% \bP\left(S_1(Y) \leq 6 D \floorsup{\beta}\sqrt{\rho l \log(\frac{nk}{\beta})} \right) > 1 - \frac{3}{nk}
	% \end{align}

	\subsection{Proof of Proposition~\ref{prop:S2_dense_case}: control of the $S_2(Z)$ term from \cite{MR3901009}}
	In this section, we prove Proposition~\ref{prop:S2_dense_case}. We follow the proof from \cite{MR3901009} but we only consider the ``dense case'' which is when $n\delta\nu\geq \log n$ -- we recall that $\nu=\max(2p-1,1-2q)$. For a similar uniform control of $S_2(Z)$ in the ``sparse case '', when $c_0\leq n\delta\nu\leq \log n$ for some absolute constant $c_0$, we refer the reader to \cite{MR3901009}. 
	
	For all $Z\in\cC$, we have $S_2(Z) = \inr{\porth{Z - Z^*}, W} = \inr{\porth{Z}, W}$ because, by construction of the projection operator, $\porth{Z^*} = 0$. Therefore, $S_2(Z)\leq \norm{\porth{Z}}_*\normop{W}$ where $\norm{\cdot}_*$ denotes the nuclear norm (i.e. the sum of singular values) and  $\normop{\cdot}$ denotes the operator norm (i.e. maximum of the singular value). In the following Lemma~\ref{lem:S2_trace}, we prove an upper bound for $ \norm{\porth{Z}}_*$ and then,  we will obtain a high probability upper bound onto $\normop{W}$. 
	
	\begin{lemma}\label{lem:S2_trace}
		For all $Z\in\cC\cap (Z^*+r B_1^{n\times n})$, we have
		\begin{align*}
		\norm{\porth{Z}}_* = \tr{\porth{Z}} \leq \frac{c_1 k}{n} \norm{Z - Z^*}_1\leq \frac{c_1 K r}{n}. 
		\end{align*}
	\end{lemma}
	\begin{proof}Since $Z\succeq0$ so it is for $(\bI_n - UU^\top)Z(\bI_n - UU^\top)$ and  so $\porth{Z}=(\bI_n - UU^\top)Z(\bI_n - UU^\top ) \succeq 0$ therefore  $\norm{\porth{Z}}_* =  \tr{\porth{Z}}$. Next, we bound the trace of $\porth{Z}$. 
		
		Since $\bI_n - UU^\top$ is a projection operator, it is symmetric and $(\bI_n - UU^\top)^2=\bI_n - UU^\top$, moreover, $\tr{Z} = n = \tr{Z^*}$ when $Z\in\cC$ so
		\begin{align*}
		&\tr{\porth{Z}}   =  \tr{\porth{Z - Z^*}} = \tr{(\bI_n - UU^\top)(Z - Z^*)(\bI_n - UU^\top)} \\
		& = \tr{(\bI_n - UU^T)^2(Z - Z^*)} = \tr{(\bI_n - UU^T)(Z - Z^*)} =  \tr{Z}-\tr{Z^*} + \tr{UU^T(Z^*-Z)} \\
		& = \sum_{i, j} (UU^T)_{ij} (Z^*-Z)_{ij} = \sum_{k\in [K]} \sum_{i, j \in \cC_k} \frac{1}{l_k} (Z^* - Z)_{ij} \overset{\mbox{(i)}}{=} \sum_{k\in [K]} \frac{1}{l_k} \sum_{i, j \in \cC_k} |(Z^* - Z)_{ij}| \\
		& \leq  \frac{c_1 K}{n} \sum_{k\in [K]} \sum_{i, j \in \cC_k} |(Z^* - Z)_{ij}| \leq \frac{c_1K}{n} \norm{Z - Z^*}_1 
		\end{align*}
		where we used in \textit{(i)} that for $i$ and $j$ in a same community, we have $ Z^*_{ij} = 1$ and $ Z_{ij} \in [0, 1]$, thus $(Z^* - Z)_{ij} = |(Z^* - Z)_{ij}|$. Finally, when  $Z$ is in the localized set $\cC\cap (Z^*+rB_1^{n\times n})$, we have $\norm{Z-Z^*}_1\leq r$ which concludes the proof.
	\end{proof}

	Now, we obtain a high probability upper bound on  $\normop{W}$. In the following, we apply this result in the ``dense case'' (i.e. $n\delta\nu\geq \log n$) to get the uniform bound onto $S_2(Z)$ over $Z\in\cC\cap (Z^*+rB_1^{n\times n})$.

	\begin{lemma}[Lemma~4 in \cite{MR3901009}]\label{lem:normopW_dense}
		With probability at least $1-\exp(-\delta\nu n)$, $\normop{W} \leq 16 \sqrt{\delta\nu n} + 168 \sqrt{\mathrm{log}(n)}$.
	\end{lemma}

	\begin{proof}
		Let $A^\prime$ be an independent copy of $\Aobs$ and $R \in \bR^{n\times n}$ be a random symmetric matrix independent from both $\Aobs$ and $A^\prime$ whose sub-diagonal entries are independent Rademacher random variables. Using a symmetrization argument  (see Chapter~2 in \cite{MR2829871} or Chapter~2.3 in \cite{vanderVaartWellner}), we obtain for $W=A-\bE A$,
		\begin{align*}
		\bE\normop{W} = \bE\normop{\Aobs - \bE A^\prime} \overset{\mbox{(i)}}{\leq} \bE\normop{\Aobs - A^\prime}
		\overset{\mbox{(ii)}}{=} \bE\normop{(\Aobs - A^\prime)\circ R} \overset{\mbox{(iii)}}{\leq} 2 \bE\normop{\Aobs\circ R}
		\end{align*}
		where $\circ$ is the entry-wise matrix multiplication operation,  \textit{(i)} comes from Jensen's inequality, \textit{(ii)} occurs since $\Aobs - A^\prime$ and $(\Aobs - A^\prime)\circ R$ are identically distributed and \textit{(iii)} is the triangle inequality. Next, we obtain an upper bound onto $\bE\normop{\Aobs\circ R}$.
		
		We define the family of independent random variables $(\xi_{ij}:1\leq i \leq j\leq n)$  where for all $1\leq i\leq j\leq n$
		\begin{equation}\label{eq:def_xi_ij}
		\xi_{ij} = \left\{ \begin{array}{l}
		\frac{1}{\sqrt{|\bE \Aobs_{ij}|}} \mbox{ with probability } \frac{\bE\Aobs_{ij}}{2}\\
		- \frac{1}{\sqrt{|\bE\Aobs_{ij}|}} \mbox{ with probability } \frac{\bE\Aobs_{ij}}{2} \\
		0 \mbox{ with probability } 1-\bE\Aobs_{ij}.
		\end{array}
		\right.
		\end{equation}We also put $b_{ij}:=\sqrt{|\bE\Aobs_{ij}|}$ for all $1\leq i \leq j\leq n$. It is straightforward to see that $(\xi_{ij}b_{ij}:1\leq i\leq j\leq n)$ and $(\Aobs_{ij}R_{ij}: 1\leq i\leq j\leq n)$ have the same distribution. As a consequence, $\normop{\Aobs \circ R}$ and $\normop{X}$ have the same distribution where $X \in \bR^{n\times n}$ is a symmetric matrix with $X_{ij} = \xi_{ij}b_{ij}$ for $1\leq i \leq  j\leq n$. An upper bound on $\bE\normop{X}$ follow from the next result due to \cite{bandeira2016sharp}.
		
		\begin{Theorem}[Corollary~3.6 in \cite{bandeira2016sharp}]\label{theo:ramon}
			Let $\xi_{ij}, 1\leq i\leq j\leq n $ be independent symmetric random variables with unit variance and $(b_{ij}, 1\leq i\leq j\leq n)$ be a family of real numbers. Let $X \in \bR^{n\times n}$ be the random symmetric matrix defined by $X_{ij} = \xi_{ij}b_{ij}$ for all $1\leq i\leq j\leq n$. Let $\sigma := \max_{1\leq i\leq n}\left\{ \sqrt{\sum_{j=1}^n b_{ij}^2}\right\}$. Then, for any $\alpha \geq 3$, 
			\begin{align*}
			\bE \normop{X} \leq e^{\frac{2}{\alpha}} \left[2\sigma + 14\alpha\max_{1\leq i\leq j\leq n}\left\{ \norm{\xi_{ij}b_{ij}}_{2\floorsup{\alpha\log(n)}}\right\}\sqrt{\log(n)}\right]
			\end{align*}
			where, for any $q>0$, $\norm{\cdot}_{q}$ denotes the $L_q$-norm.
		\end{Theorem}
		
		Since  $(\xi_{ij}:1\leq i\leq j\leq n)$ are independent symmetric such that $\mathrm{Var}(\xi_{ij}) = \bE[\xi_{ij}^2] = 1$ we can apply Lemma~\ref{theo:ramon} to upper bound $\bE \normop{X}=\bE\normop{A\circ R}$.  We have $\norm{\xi_{ij}b_{ij}}_{2\floorsup{\alpha\log(n)}} \leq 1$ for any $\alpha \geq 3$ and $b_{ij}^2 = |\bE \Aobs_{ij}|\leq \delta\nu$. It therefore follows from Lemma~\ref{theo:ramon} for $\alpha = 3$ that
		\begin{equation}\label{eq:ramon}
		\bE\normop{W} \leq 2e^{\frac{2}{3}} \left[ 2\sqrt{n\delta\nu} + 42 \sqrt{\log(n)}\right] \leq 8\sqrt{n\delta\nu} + 168 \sqrt{\log(n)}.
		\end{equation}
		
		The final step to prove Lemma \ref{lem:normopW_dense} is a concentration argument showing that $\normop{W}$ is close to its expectation with high probability. To that end we rely on a general result for Lipschitz and separately convex functions from \cite{BouLugMass13}. We first recall that a real-valued function $f$ of $N$ variables is said separately convex when for every $i=1, \ldots, N$ it is a convex function of the $i$-th variable if the rest of the variables are fixed.

		\begin{Theorem}[Theorem 6.10 in \cite{BouLugMass13}]\label{lipConcentration} 
			Let $\cX$ be a convex compact set in $\bR$ with diameter $B$. Let $X_1,\cdots, X_N$ be independent random variables taking values in $\cX$. Let $f:\cX^N \rightarrow \bR$ be a separately convex and $1$-Lipschitz function, w.r.t. the $\ell_2^N$-norm (i.e. $|f(x)-f(y)| \leq \norm{x-y}_2$ for all $x, y \in \cX^N$). Then $Z=f(X_1,\ldots, X_N)$ satisfies, for all $t>0$, with probability at least $1-\exp(-t^2/B^2)$, $Z \leq \bE[Z]+t$.
		\end{Theorem}
		
		We apply  Theorem~\ref{lipConcentration} to $Z:=\normop{W} = f(A_{ij}-\E A_{ij}, 1\leq i\leq j\leq n) = \frac{1}{\sqrt{2}}\normop{\Aobs - \bE\Aobs}$ where $f$ is a $1$-Lipschitz w.r.t. $\ell_2^N$-norm for $N=n(n-1)/2$ and separately convex function and $(A_{ij}-\E A_{ij}, 1\leq i\leq j\leq n)$ is a family  of $N$ independent random variables. Moreover, for each $i\geq j$, $(\Aobs - \bE \Aobs)_{ij} \in  [-1-\delta(2p-1),1-\delta(2q-1)]$, which is a convex compact set with diameter $B=2(1+\delta(p-q)) \leq 4$. Therefore, it follows from Theorem~\ref{lipConcentration} that for all $t>0$, with probability at least $1-\exp(-t^2/16)$, $\normop{W}\leq \E \normop{W} + \sqrt{2}t$.  In particular, we finish the proof of Lemma~\ref{lem:normopW_dense} for $t = 4\sqrt{\delta\nu n}$ and using the bound from \eqref{eq:ramon}. 
	\end{proof}
	
	It follows from Lemma~\ref{lem:normopW_dense} that when $n\nu \delta\geq \log n$, $\normop{W}\leq 184\sqrt{n\delta\nu}$ with probability at least $1-\exp(-\delta\nu n)$. Using this later result together with Lemma~\ref{lem:S2_trace} concludes the proof of Proposition~\ref{prop:S2_dense_case}.

	\section{Annex~2: Solving SDP's in practice}
	The practical implementation of our approach to the problems of synchronization, signed clustering and MAX-CUT resort to solving a convex optimization problem. In the present section, we describe the various algoritms we used for solving these SDP's.

\subsection{Pierra's method}	
	For SDP's with constraints on the entries, we propose a simple modification of the method initially proposed by Pierra in \cite{pierra1984decomposition}.  Let $f$: $\mathbb R^{n\times n}\mapsto \mathbb R$ be a convex function. Let $\mathcal C $ denote a convex set which can be written as the intersection of convex sets $\mathcal C = \mathcal S_1 \cap \cdots\cap \mathcal S_J$.  
	Let us define $\mathbb H= \mathbb R^{n \times n} \times \cdots \times \mathbb R^{n \times n}$ ($J$ times) and let $\mathbb D$ denote the (diagonal) subspace of $\mathbb H$ of vectors of the form $(Z,\ldots,Z)$. In this new formalism, the problem can now be formulated as a minimisation problem over the intersection of two sets only, i.e. 

\begin{equation*}
    \min_{{\sf Z} \in \mathbb H} \left( \frac1{J} \ \sum_{j=1}^J \ f({\sf Z}_j):{\sf Z}=({\sf Z}_j)_{j=1}^J \in \left(\mathcal S_1 \times \cdots \times \mathcal S_J\right) \cap \mathbb D\right).
\end{equation*}

Define $F(\sf Z) = \frac1{J} \ \sum_{j=1}^J \ f({\sf Z}_j)$. 
The algorithm proposed by Pierra in \cite{pierra1984decomposition} consists of performing the following iterations 
\begin{equation*}
{\sf Z}^{p+1}  =\operatorname{Prox}_{I_{\mathcal S_1 \times \cdots \times \mathcal S_J}+\frac{1}{2} \varepsilon F}\left(\sf{B}^{p}\right) \mbox{ and }
\sf{B}^{p+1}  =\operatorname{Proj}_{\mathbb{D}} ({\sf Z}^{p+1}).
\end{equation*}

	\subsubsection{Application to community detection}
	Let us now present the computational details of Pierra's method for the community detection problem. We will estimate its membership matrix $\bar{Z}$ via the following SDP estimator
	\begin{equation*}
	\hat{Z} \in \text{argmax}_{Z \in \mathcal{C}}  \langle A, Z\rangle, 
	\end{equation*}
	where $\mathcal{C} = \{ Z \in \mathbb R^{n\times n}, Z \succeq 0, Z \geq 0, {\rm diag}(Z)\preceq I_n, \sum Z_{ij} \leq \lambda \}$  and $
	    \lambda  = \sum_{i, j=1}^n \bar{Z}_{ij} = \sum_{k=1}^K |\mathcal{C}_k|^2$ denotes the number of nonzero elements in the membership matrix $\bar{Z}$. The motivation for this approach stems from the fact that the membership matrix  $\bar{Z}$ is actually the oracle, i.e., $Z^*=\bar Z$ , where $Z^* \in \text{argmax}_{Z \in \mathcal C} \langle  \mathbb E [A], Z\rangle$. The function $f$ to minimize in the Pierra algorithm is defined as $f(Z)  = -\langle A,Z \rangle$.
	
	Let us denote by $\mathbb S_+$ the set of symmetric positive semi-definite matrices in $\mathbb R^{n \times n}$. The set $\mathcal C$ is the intersection of the sets 
	\begin{align*}
	    S_1&  = \mathbb S_+ \,
	    S_2  = \left\{ Z \in \mathbb R^{n\times n} \mid Z \ge 0\right\} \,
	    S_3  = \left\{ Z \in \mathbb R^{n\times n} \mid \text{diag}(Z) \preceq I\right\} \\
	    &\mbox{ and } S_4  = \left\{ Z \in \mathbb R^{n\times n} \mid \sum_{i,j=1}^n \ Z_{ij} \le \lambda\right\} 
	\end{align*}
	
We now compute for all $\sf{B}=(\sf{B}_j)_{j=1}^4\in (\bR^{n\times n})^4$ and $j=1,\ldots,4$ ($J=4$ here)
\begin{align*}
    \operatorname{Prox}_{I_{S_1 \times \cdots \times S_4}+\frac{1}{2} \varepsilon F}\left(\sf{B}\right)_j & = \operatorname{Prox}_{I_{S_j}+\frac{1}{2 J} \varepsilon f}\left(\sf{B}_j\right)
\end{align*}
We have for $J=4$
\begin{equation*}
    \operatorname{Prox}_{I_{S_j}+\frac{1}{2 J} \varepsilon f}\left(\sf{B}_j\right)  = 
    \text{argmin}_{Z \in S_j} \ -\frac{\epsilon}{2J} \ \langle A,Z\rangle+ \frac12 \ \Vert Z-{\sf B}_j\Vert_F^2  = P_{S_j}\left({\sf B}_j+\frac{\epsilon}{2J} \ A\right)
\end{equation*}
On the other hand, the projections operators $P_{S_j},j=1,2,3,4$ are given by 
\begin{align*}
    P_{S_1}(\sf Z_1) & = U \max\left\{\Sigma,0\right\} U^\top  \hspace{.3cm} \text{where $\sf Z_1$ has SVD } \sf Z_1=U\Sigma U^\top\\
    P_{S_2}(\sf Z_2) & = \max \left\{{\sf Z}_2,0\right\} \,
    P_{S_3}(\sf Z_3)  = {\sf Z}_3 -\text{diag}({\sf Z}_3)+\min \left\{1,\text{diag}({\sf Z}_3)\right\} \\
    P_{S_4}({\sf Z}_4) & =  \frac{\lambda}{\sum_{ij} \ ({\sf Z}_4)_{ij}} \ \sf Z_4.
\end{align*}
To sum up, Pierra's method can be formulated as follows. 

\vspace{.3cm}

For all iterations $k$ in $\mathbb N$, compute the SVD of ${\sf B}_1^k+\frac{\epsilon}{2 \cdot 4}A=U^k \Sigma^k (U^{k})^\top$. Then compute for all $j=1,\ldots,4$
\begin{align*}
    \sf B^{k+1}_j & = \frac1{4} \Bigg(U^k \max\left\{\Sigma^k,0\right\} (U^{k})^\top+\max \left\{ \sf {\sf B}_2^k+\frac{\epsilon}{2 \cdot 4}A,0\right\} + {\sf B}_3^k+\frac{\epsilon}{2 \cdot 4}A -\text{diag}({\sf B}_3^k+\frac{\epsilon}{2 \cdot 4}A)\\
    & \hspace{2cm}+\min \left\{1,\text{diag}({\sf B}_3^k+\frac{\epsilon}{2 \cdot 4}A)\right\}+ 
    \frac{\lambda}{\sum_{ij} \ ({\sf B}_4^k+\frac{\epsilon}{2 \cdot 4}A)_{ij}} \ {\sf B}_4^k+\frac{\epsilon}{2 \cdot 4}A
    \Bigg).
\end{align*}

	\subsubsection{Application to signed clustering}
Let us now turn to the signed clustering problem. 
We will estimate its membership matrix $\bar{Z}$ via the following SDP estimator $\hat{Z} \in \text{argmax}_{Z \in \mathcal{C}}  \ \langle A, Z\rangle,$ where $ \mathcal C = \{Z\in\mathbb R^{n\times n}: Z\succeq 0, Z_{ij}\in[0,1], Z_{ii}=1, i=1, \ldots, n\}$. As in the community detection case the function $f$ to minimize in the Pierra algorithm is defined as $f(Z)  = -\langle A,Z \rangle$.

Let us denote by $\mathbb S_+$ the set of symmetric positive semi-definite matrices in $\mathbb R^{n \times n}$. The set $\mathcal C$ is the intersection of the sets $S_1  = \mathbb S_+$, $S_2  = \left\{ Z \in \mathbb R^{n\times n} \mid Z \in  [0,1]^{n \times n}\right\}$ and  $S_3  = \left\{ Z \in \mathbb R^{n\times n} \mid Z_{ii} =1, \ i=1,\ldots,n\right\}.$

As before, for $j=1,\ldots,3$
\begin{align*}
    \operatorname{Prox}_{I_{S_j}+\frac{1}{2 \cdot 3} \varepsilon f}\left(\sf{B}_j\right) & = P_{S_j}\left({\sf B}_j+\frac{\epsilon}{2 \cdot 3} \ A\right)
\end{align*}and the projection operators $P_{S_j}$, $j=1,2$ are given by 
\begin{equation*}
    P_{S_1}(\sf Z_1)  = U \max\left\{\Sigma,0\right\} U^\top, \,
    P_{S_2}(\sf Z_2)  = \min \left\{\max \left\{{\sf Z}_2,0\right\},1\right\} \mbox{ and }
    P_{S_3}(\sf Z_3)  = {\sf Z}_3-\text{diag}({\sf Z}_3)+I
\end{equation*}
To sum up, Pierra's method can be formulated as follows. 

\vspace{.3cm}

At each iteration $k$, compute the SVD of ${\sf B}_1^k+\frac{\epsilon}{2 \cdot 3}A=U^k \Sigma^k (U^{k})^\top$. Then compute for all $j=1,\ldots,3$
\begin{align*}
    {\sf B}^{k+1}_j & = \frac1{3} \Bigg(U^k \max\left\{\Sigma^k,0\right\} U^{k^t}+\min \left\{\max \left\{ {\sf B}_2^k+\frac{\epsilon}{2 \cdot 3} \ A,0\right\},1\right\}+
    {\sf B}_3^k+\frac{\epsilon}{2 \cdot 3} \ A- \text{diag}\left({\sf B}_3^k+\frac{\epsilon}{2 \cdot 3} \ A \right)+I
    \Bigg).
\end{align*}

\subsection{The Burer-Monteiro approach and the \textsc{Manopt} Solver}

To solve the MAX-CUT and Angular Synchronization problems we rely on \textsc{Manopt}, a freely available Matlab toolbox for optimization on manifolds \cite{manopt}. \textsc{Manopt} runs the Riemannian Trust-Region method on corresponding Burer-Monteiro non-convex problem with rank bounded by $p$ as follows. The Burer-Monteiro approach consists of replacing the optimization of a linear function $ \inr{A,Z}$ over the convex set $ \mathcal{Z} = \{ Z \succeq 0: \mathcal{A}(Z)=b \}$ with the optimization of the quadratic function $\inr{AY,Y} $ over the non-convex set $ \mathcal{Y} = \{ Y \in \mathbb{R}^{n \times p}: \mathcal{A}(YY^T)=b \}$.

In the context of the MAX-CUT problem, the Burer-Monteiro approach amounts to the following steps. Denoting by $Z$ the positive semidefinite matrix $Z = z z^T$, note that both the cost function and the constraints lend themselves to be expressed linearly in terms of $Z$. Dropping the NP-hard rank-1 constraint on $Z$, we arrive at the well-known convex relaxation of MAX-CUT from \cite{goemans1995improved}
\begin{equation}\label{eq:SDP_relax_max_cut}
	\hat Z \in \underset{Z\in \mathcal{C}}{\mathrm{argmin}} \inr{A, Z} 
\end{equation}
where $\mathcal{C} := \{Z \in \bR^{n\times n}:Z \succeq 0, Z_{ii} =  1, \forall i=1,\ldots, n\}$.

If a solution $\hat{Z}$ of this SDP has rank 1, then $\hat{Z} = z^* z^{*^T}$ for some $z^*$, which then gives the optimal cut. Recall that in the general case of higher rank $\hat{Z}$, Goemans and Williamson \cite{goemans1995improved} introduced the celebrated rounding scheme that extracts approximately optimal cuts within a ratio of 0.878 from $\hat{Z}$. The corresponding Burer-Monteiro non-convex problem with rank bounded by $p$ is given by 
\begin{equation}\label{eq:SDP_relax_max_cut_Burer_Monteiro}
	\hat{Y} \in \underset{X\in \mathcal{B}}{\mathrm{argmin}} \inr{AY, Y} 
\end{equation}
where $\mathcal{B} := \{Y \in \bR^{n\times p}: \text{diag}(Y Y^T) = 1\}$.
Note that the  constraint $ \text{diag}(Y Y^T) = 1$ requires each row of $Y$ to have unit $\ell_2^p$ norm, rendering $Y$ to be a point on the Cartesian product of $n$ unit spheres $\cS_2^{p-1}$ in $\bR^p$, which is a smooth manifold. Also note that the search space of the SDP is compact, since all $Z$ feasible for the SDP have identical trace equal to $n$.

If the convex set $ \mathcal{Z}$ is compact, and $m$ denotes the number of constraints, it holds true that whenever $p$ satisfies $ \frac{p(p+1)}{2} \geq m$, the two problems share the same global optimum \cite{Barvinok1995, Burer2005}. Building on pioneering work of Burer and Monteiro \cite{Burer2005}, Boumal et. al \cite{boumal2016bmapproach} showed that if the set $\mathcal{Z}$ is compact and the set $ \mathcal{Y} $ is a smooth manifold, then $ \frac{p(p+1)}{2} \geq m$ implies that, for almost all cost matrices $A$, global optimality is achieved by any $Y$ satisfying a second-order necessary optimality conditions. Following  \cite{boumal2016bmapproach}, \textit{for $ p = \lceil \sqrt{2n} \rceil$, for almost all matrices $A$, even though \eqref{eq:SDP_relax_max_cut_Burer_Monteiro} is non-convex, any local optimum $Y$ is a global optimum (and so is $Z = Y Y^T$), and all saddle points have an escape (the Hessian has a negative eigenvalues)}. Note that for $p > n/2$ the same statement holds true for \textit{all} $A$, and was previously established by \cite{boumal2015riemannian}.

%%%%%%%%%%%%%%%%%%%%%%%%%%%%%%%%%%%%%%%%%%%%%%%%%%%%%%%%%%%%%%%%%%%%%%%%%%-----------------------
%%%%%%%%%%%%%%%%%%%%%%%%%%%%%%%%%%%%%%%%%%%%%%%%%%%%%%%%%%%%%%%%%%%%%%%%%%-----------------------

\section{Numerical experiments}

This section contains the outcome of numerical experiments on the three  application problems considered: 
signed clustering, 
MAX-CUT, 
and angular synchronization. 
% \end{itemize} 

\subsection{Signed Clustering}
To assess the effectiveness of the SDP relaxation, we consider the following experimental setup. We generate synthetic networks following the signed stochastic block model (SSBM) previously described in Section \ref{sec:SSBM}, with $K=5$ communities. To quantify the effectiveness of the SDP relaxation, we compare the accuracy of a suite of algorithms from the signed clustering literature, \textit{before} (that is when we perform these algorithms directly on $A$) and  \textit{after} the SDP relaxation (that is when we perform the very same algorithms on $\hat Z$). To measure the recovery quality of the clustering results, for a given indicator set $x_1, \ldots, x_K$, we rely on the error rate consider in \cite{Chiang_2012_Scalable}, defined as
\begin{equation}  \label{errorSSBM}
   \gamma =  \sum_{c=1}^{K} \frac{x_c^T A_{com}^{-} x_c   +  x_c^T L_{com}^{+}  x_c  } { n^2}, 
\end{equation}
where $x_c$ denotes a cluster indicator vector, $A_{com}$ ($=\bE A$) is the complete $K$-weakly balanced ground truth network  -- with $1$'s on the diagonal blocks corresponding to inter-cluster edges, and $-1$ otherwise -- with $A_{com} = A_{com}^{+} - A_{com}^{-}$, and  $L_{com}^{+}$ denotes the combinatorial graph Laplacian corresponding to $A_{com}^{-}$. Note that $ x_c^T A_{com}^{-} x_c $ counts the number of violations within the clusters (since negative edges should not be placed within clusters) and $ x_c^T L_{com}^{+}  x_c $ counts the number of violations across clusters (since positive edges should not belong to the cut). Overall, \eqref{errorSSBM} essentially counts the fraction of intra-cluster and inter-cluster edge violations, with respect to the full ground truth matrix. Note that this definition can also be easily adjusted to work on real data sets, where the ground truth matrix $A_{com}$ is not available, which one can replace with the empirical observation $A$. 

In terms of the signed clustering algorithms compared, we consider the following algorithms from the literature. One straightforward approach is to simply rely on the spectrum of the observed adjacency matrix $A$. Kunegis et al. \cite{kunegis2010spectral} proposed  spectral tools for clustering, link prediction, and visualization of signed graphs, by  solving a 2-way ``signed'' ratio-cut problem based on the combinatorial Signed Laplacian \cite{HouSignedLap}   $\bar{L} = \bar{D} - A$, 
where $\bar{D}$ is a diagonal matrix with $\bar{D}_{ii} = \sum_{i=1}^{n} |A_{ij}|$. 
The same authors proposed signed extensions for the case of the random-walk Laplacian $\bar{L}_{\text{rw}} = I - \bar{D}^{-1} A$, and the symmetric graph Laplacian $\bar{L}_{\text{sym}} = I - \bar{D}^{-1/2} A \bar{D}^{-1/2}$, the latter of which is particularly suitable for skewed degree distributions.
Finally, the last algorithm we considered  is \textsc{BNC} of Chiang et al. \cite{DhillonBalNormCut}, who   introduced a formulation based on the   \textit{Balanced Normalized Cut}  objective
\begin{equation}
\operatorname{min}_{ \{x_1,\ldots,x_K\} \in \mathcal{I} } \left(  \sum_{c=1}^{K} \frac{x_c^T(D^{+} - A)x_c }{x_c^T \bar{D} x_c}  \right).
\label{obj_BNormC}  
\end{equation}
which, in light of the decomposition $ D^+ - A = D^+ - (A^+ - A^-) = D^+ - A^+ + A^- = L^+ +  A^-$, is effectively minimizing the number of violations in the clustering procedure. 

In our experiments, we first compute the error rate $\gamma_{before}$ of all algorithms on the original SSBM graph (shown in Column 1 of Figure \ref{fig:tableSignedClust_k5}), and then we repeat the procedure but with the input to all signed clustering algorithms being given by the output of the SDP relaxation, and denote the resulting recovery error by $\gamma_{after}$. The third column of the same  Figure \ref{fig:tableSignedClust_k5} shows the difference in errors $ \gamma_{\delta} = \gamma_{before} - \gamma_{after}$ between the first and second columns, while the fourth column contains a histogram of the error differences $ \gamma_{\delta}$. This altogether illustrates the fact that the SDP relaxation does improve the performance of all signed clustering algorithms, except $\bar{L}$, and could effectively be used as a denoising pre-processing step.

\vspace{-2mm}
\newcommand{\wid}{1.6in}
\newcolumntype{C}{>{\centering\arraybackslash}m{\wid}}
\begin{table*}[!htp]\sffamily
\hspace{2mm}
 % LOUT_END:     \begin{tabular}{l*5{C}@{}}
   % LOUT_END:   & I & II & III & IV & V \\ 
% \begin{tabular}{l*4{C}@{}}
\begin{center}
\hspace{-6mm}
\begin{tabular}{l*4{C}@{}}
\hspace{-6mm}
   & Before & After & Delta & Histogram  \\  
$A$ \hspace{-5mm}
& \includegraphics[width=\wid]{{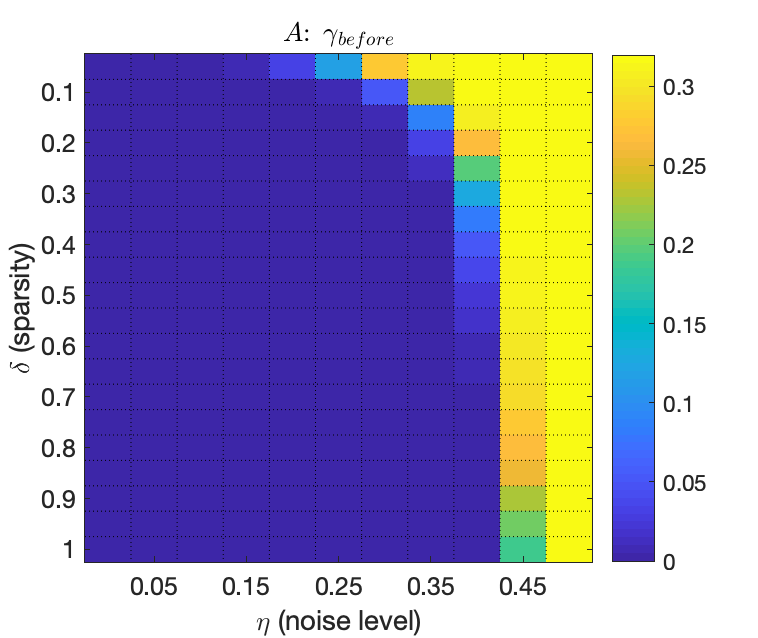}} 
& \includegraphics[width=\wid]{{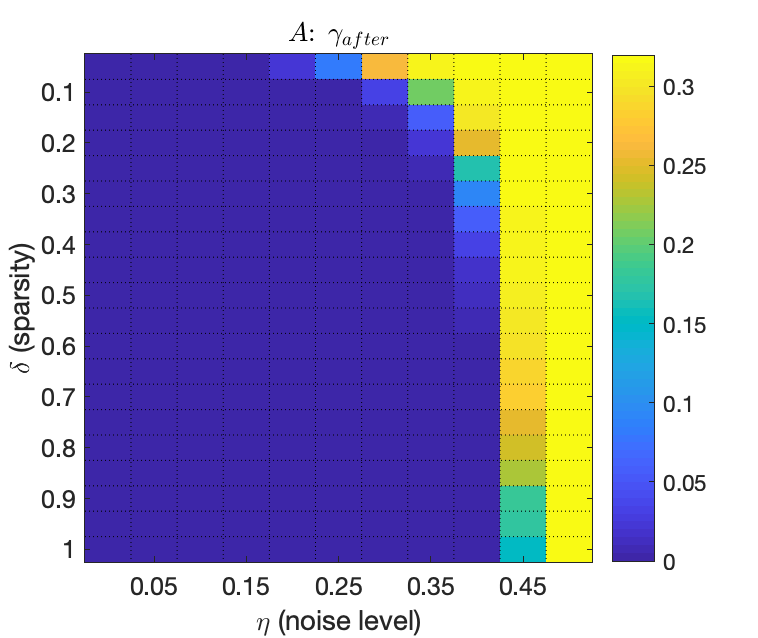}} 
& \includegraphics[width=\wid]{{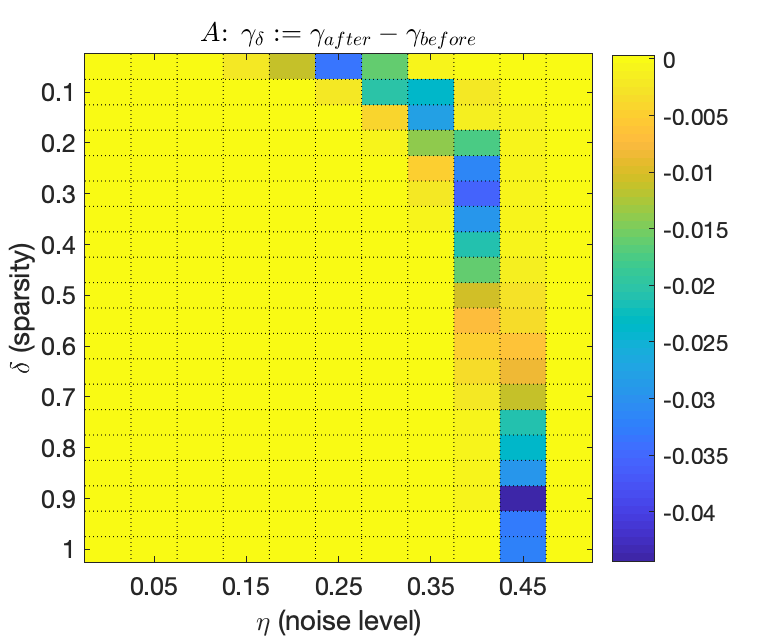}} 
& \includegraphics[width=\wid]{{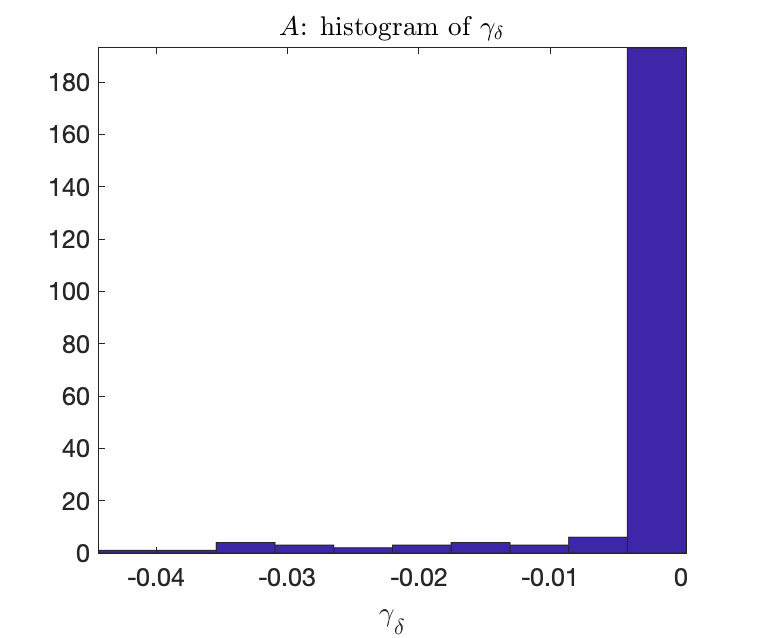}} \\
$\bar{L}$  \hspace{-5mm}
& \includegraphics[width=\wid]{{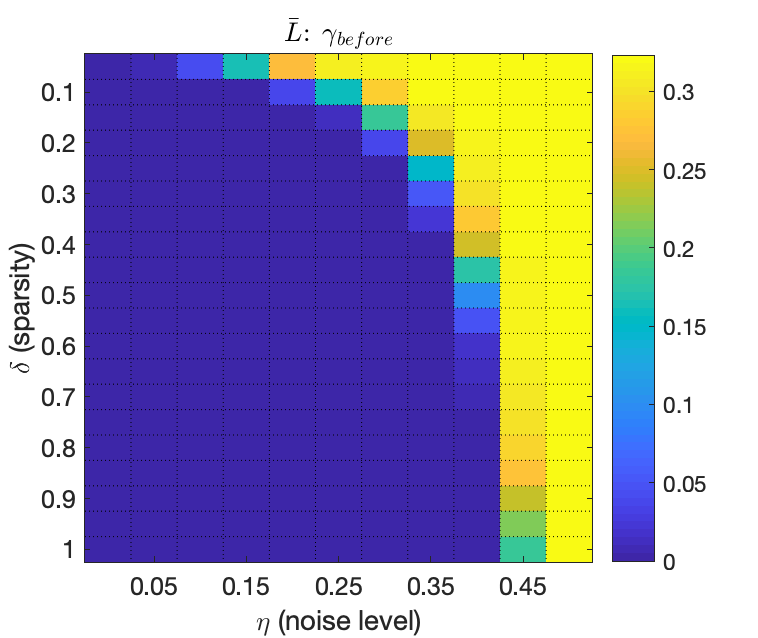}} 
& \includegraphics[width=\wid]{{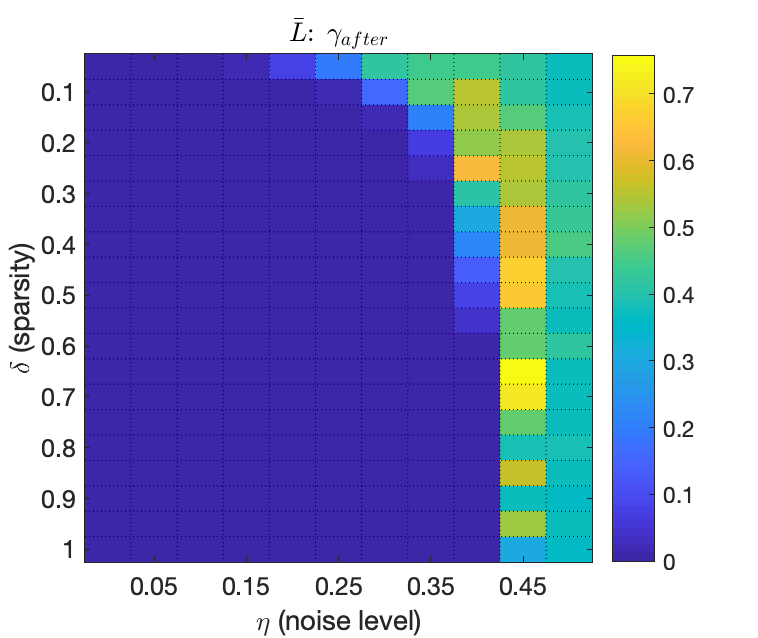}} 
& \includegraphics[width=\wid]{{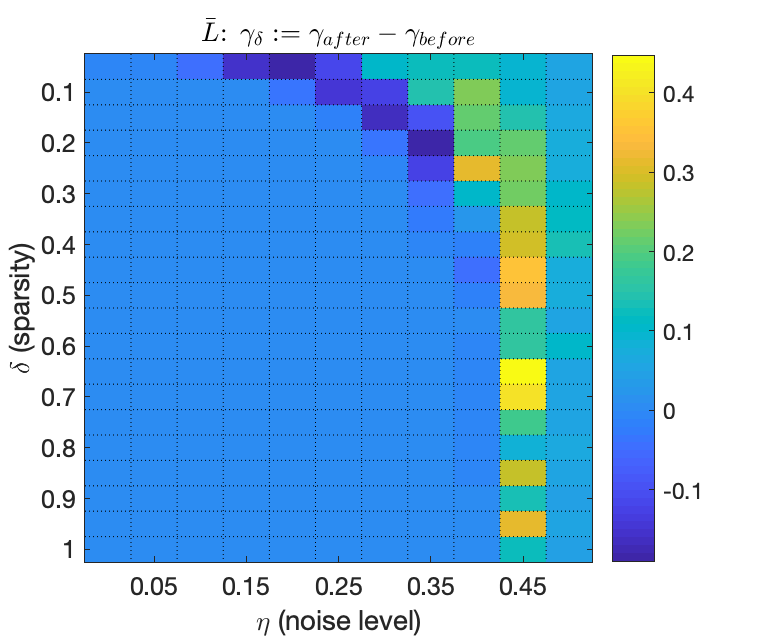}} 
& \includegraphics[width=\wid]{{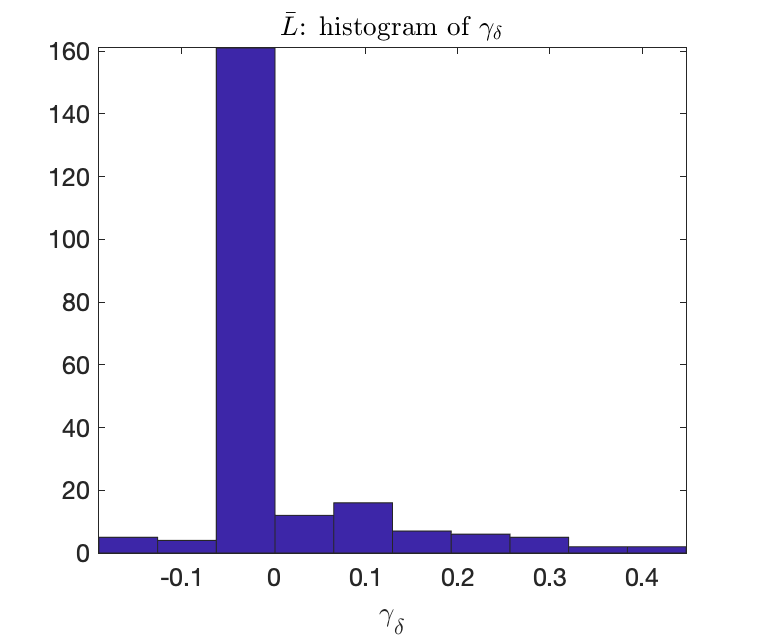}} \\
$\bar{L}_{rw}$ \hspace{-5mm}
& \includegraphics[width=\wid]{{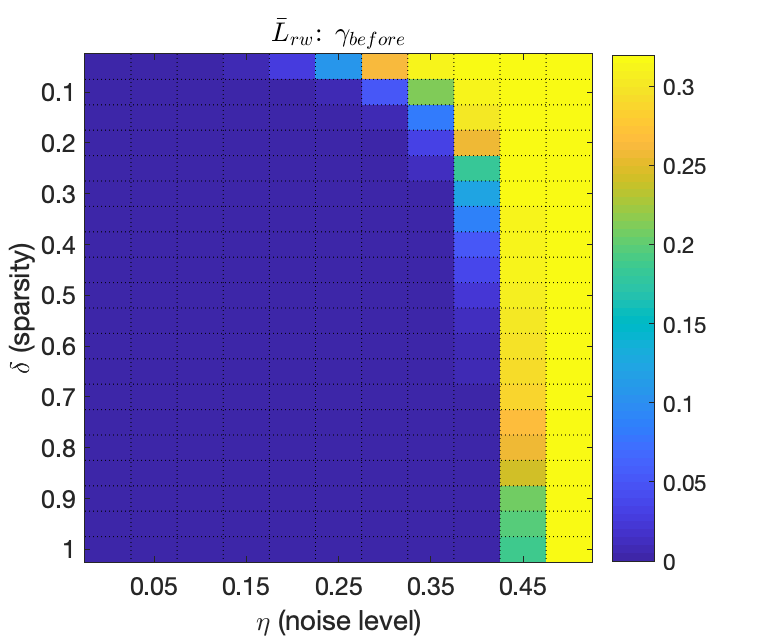}} 
& \includegraphics[width=\wid]{{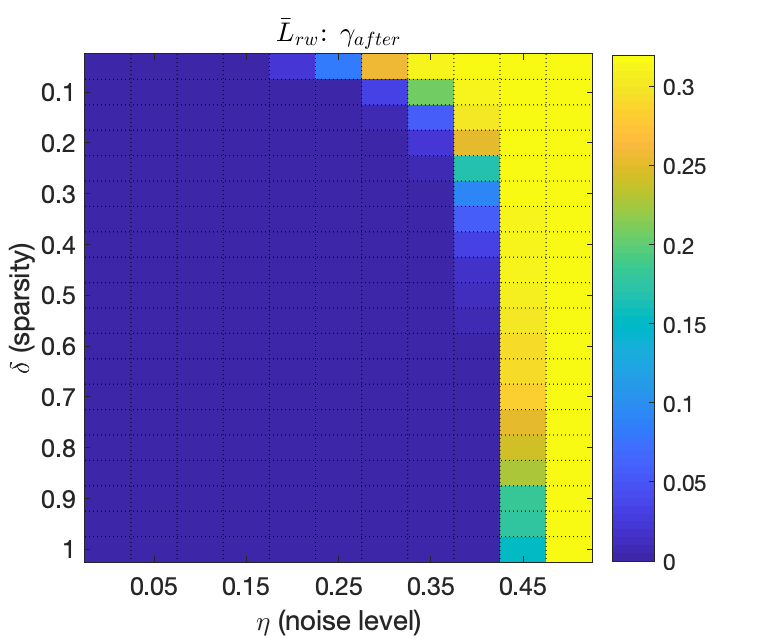}} 
& \includegraphics[width=\wid]{{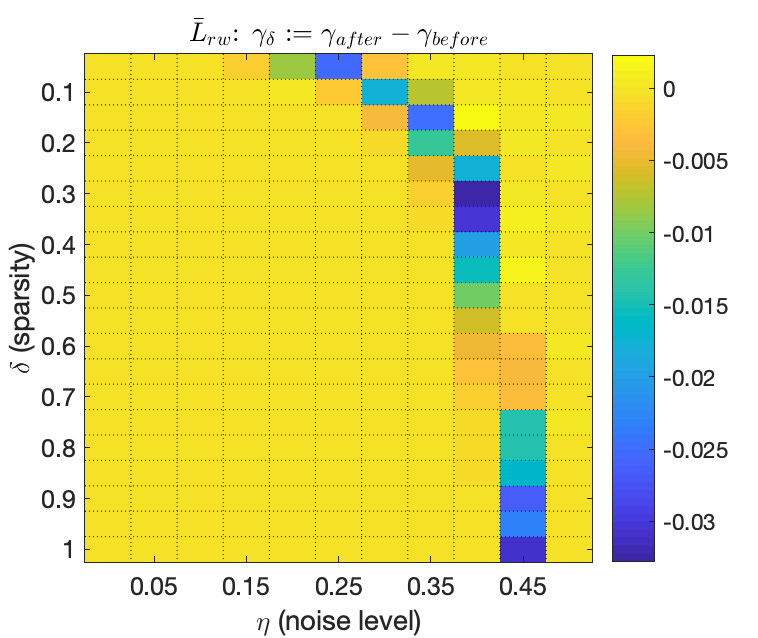}} 
& \includegraphics[width=\wid]{{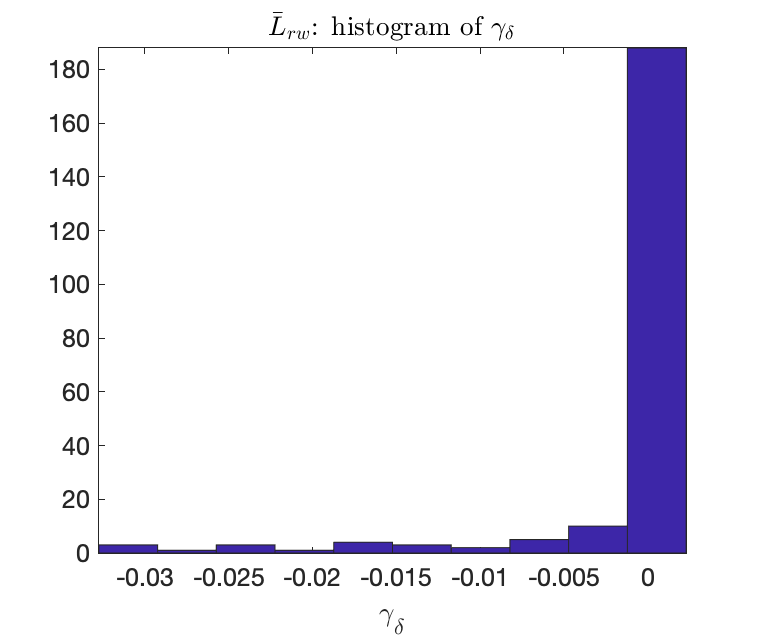}} \\
$\bar{L}_{sym}$ \hspace{-5mm}
& \includegraphics[width=\wid]{{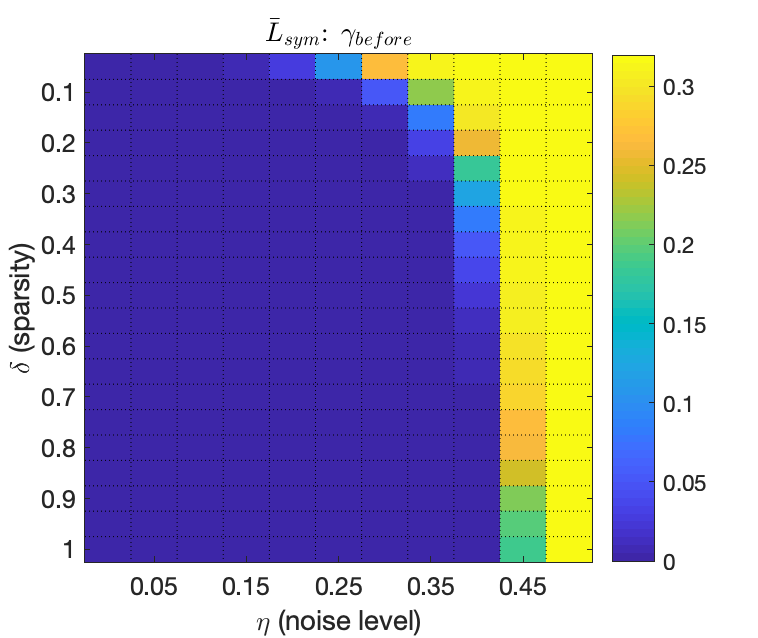}} 
& \includegraphics[width=\wid]{{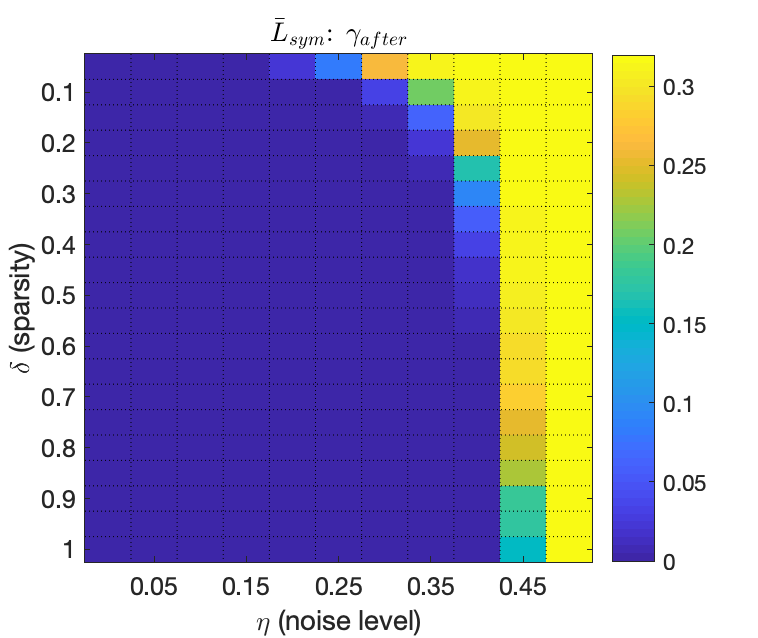}} 
& \includegraphics[width=\wid]{{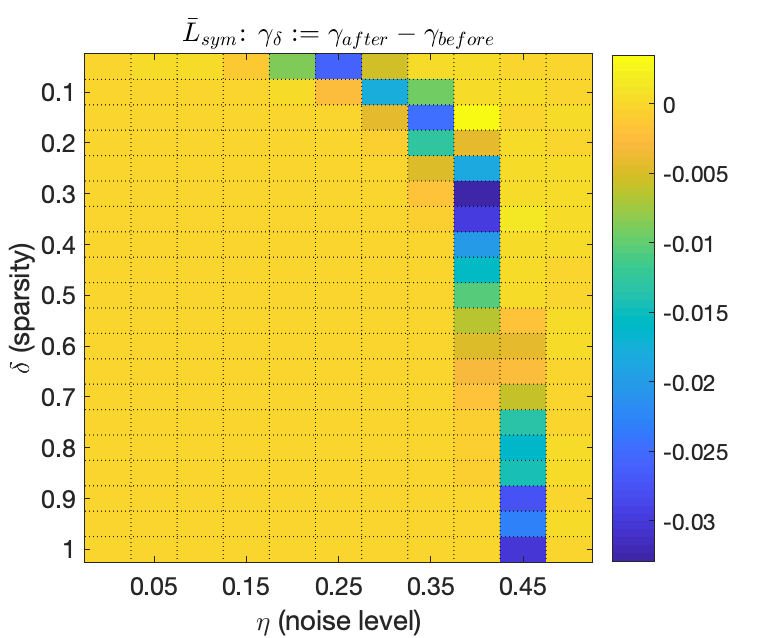}} 
& \includegraphics[width=\wid]{{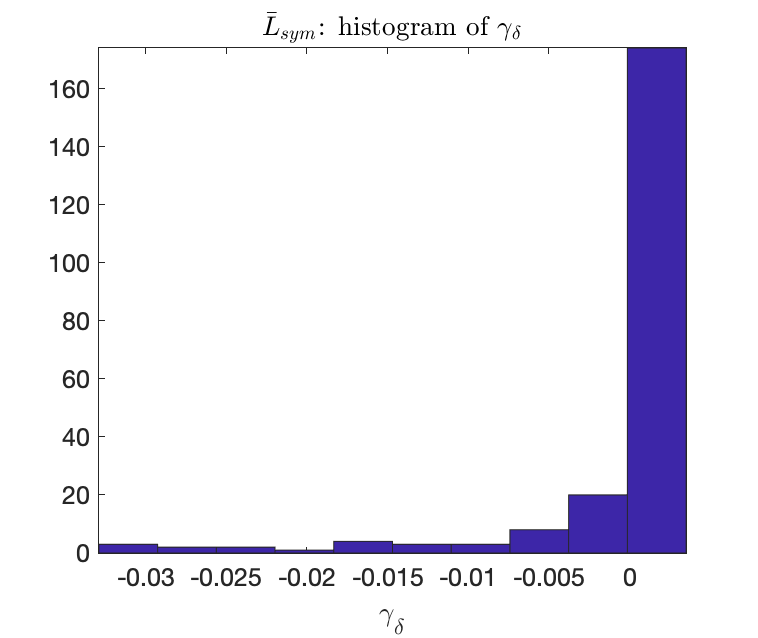}} \\
\textsc{BNC} \hspace{-5mm}
& \includegraphics[width=\wid]{{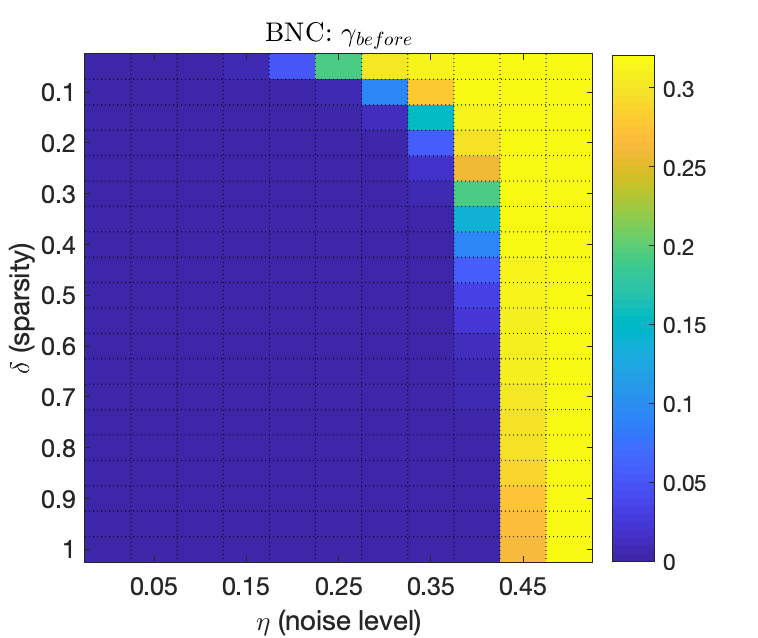}} 
& \includegraphics[width=\wid]{{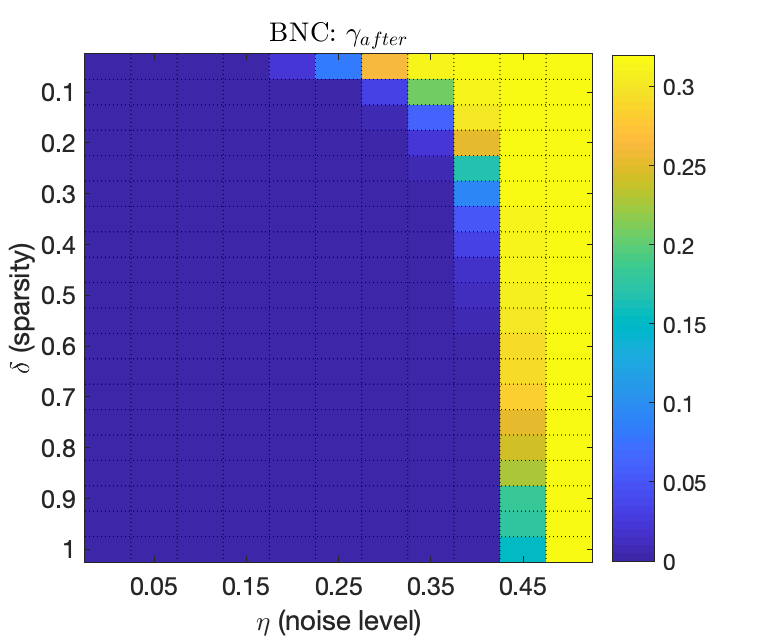}} 
& \includegraphics[width=\wid]{{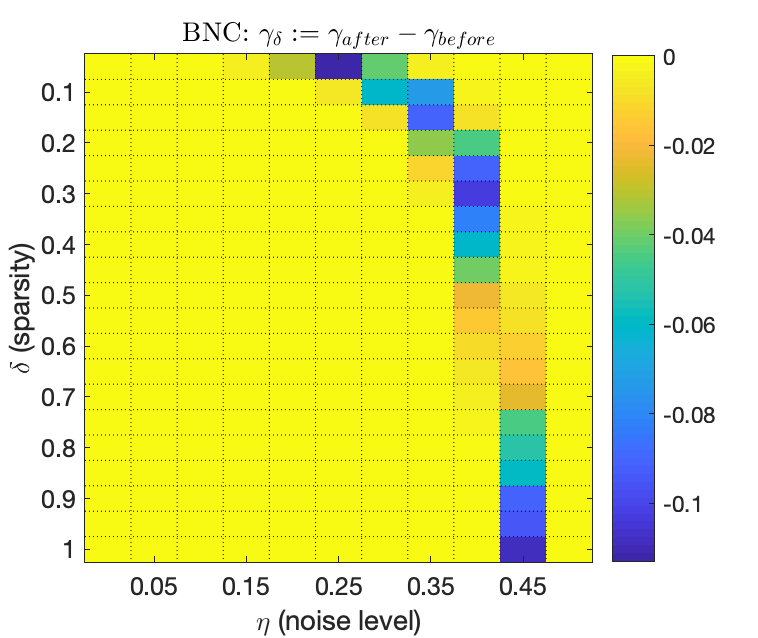}} 
& \includegraphics[width=\wid]{{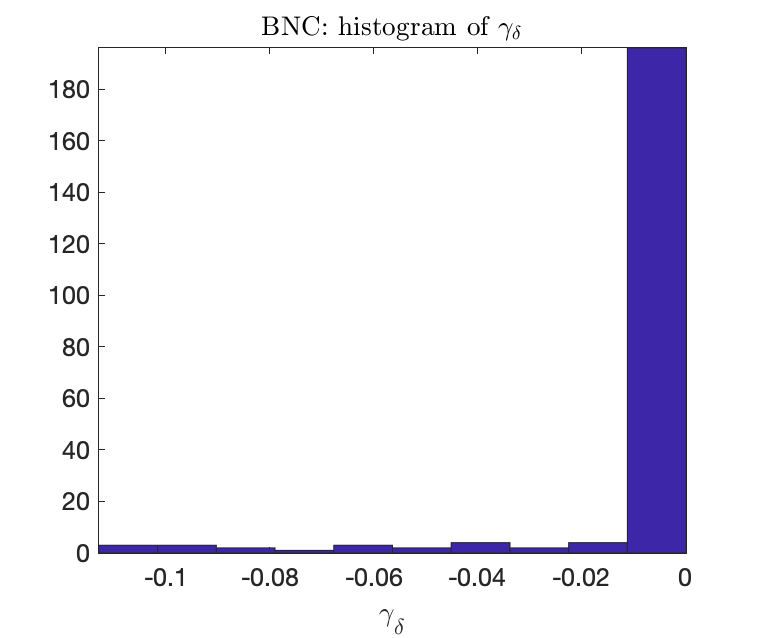}} \\
\end{tabular}
\end{center}
\captionsetup{width=0.99\linewidth}
\vspace{-1mm}
\captionof{figure}{Summary of results for the Signed Clustering problem. The first column denotes the recovery error \textbf{before} the SDP relaxation step, meaning that we consider a number of signed clustering algorithms from the literature which we apply directly the initial adjacency matrix $A$. The second column contains the results when applying the same suite of algorithms \textbf{after} the SDP relaxation. The third column shows the difference in errors between the first and second columns, while the fourth column contains a histogram of the delta errors. This altogether illustrates the fact the SDP relaxation does improve the performance of all signed  clustering algorithms except $\bar{L}$.
}
\label{fig:tableSignedClust_k5}	
\end{table*}

\subsection{Max-Cut}
\begin{figure}[h!]
\begin{center}
\includegraphics[width=0.35\columnwidth]{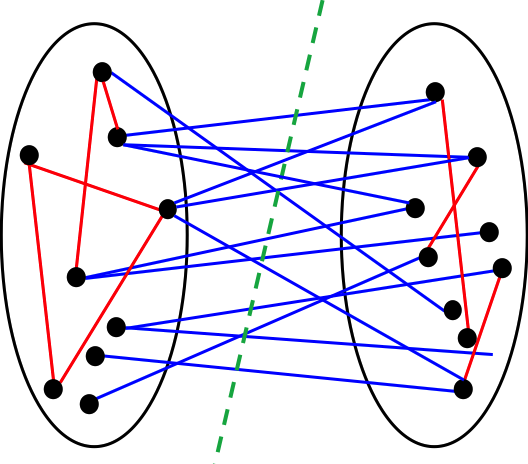}
\end{center} 
\caption{Illustration of Max-Cut in the setting of a perturbation of a complete bipartite graph.}
\label{fig:BipartitePerturbation} 
\end{figure}

\begin{figure}[h!]
\captionsetup[subfigure]{skip=1pt}
\begin{centering}\hspace{2mm}
\subcaptionbox{ Adjusted Rand Index. }[0.45\columnwidth]{\includegraphics[width=0.4\columnwidth]{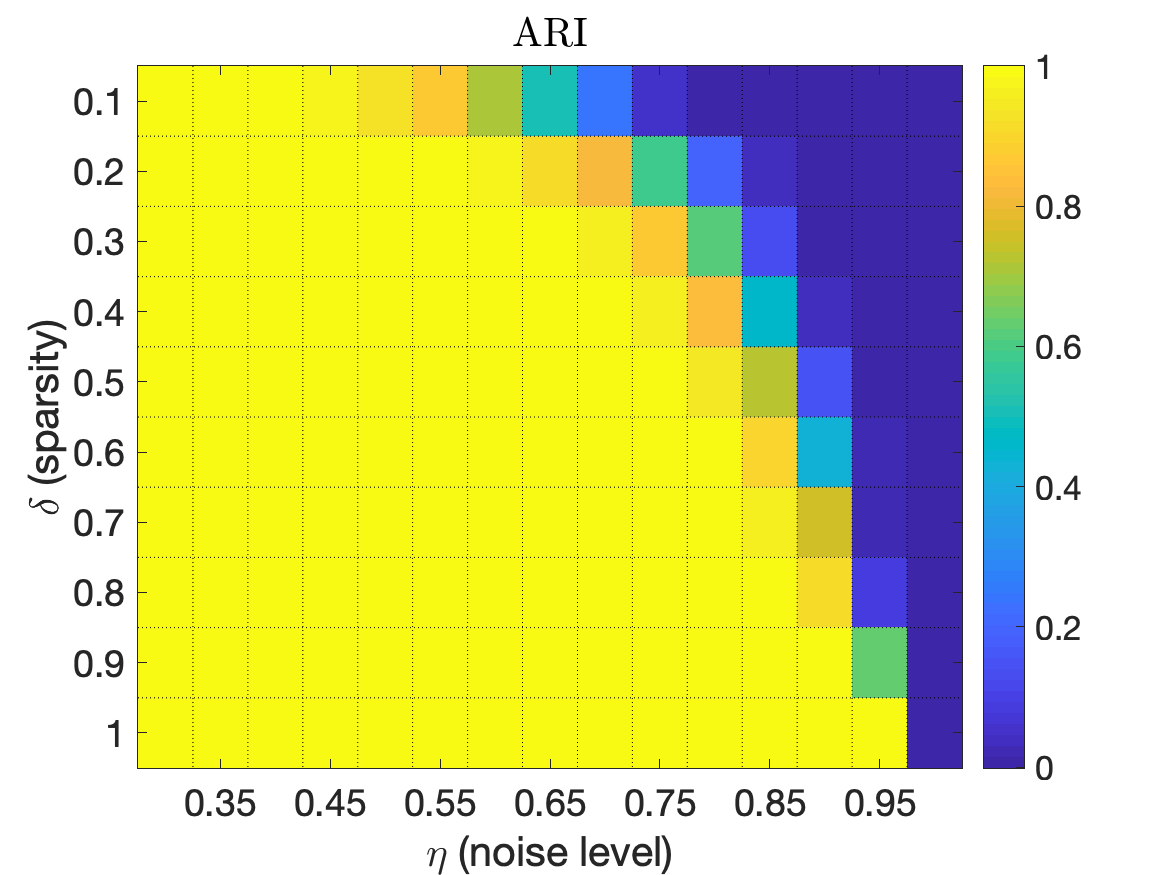}}
\hspace{3mm} 
\subcaptionbox{Running times (\textsc{MANOPT}). % , as a function of noise level and graph sparsity. 
}[0.45\columnwidth]{\includegraphics[width=0.4\columnwidth]{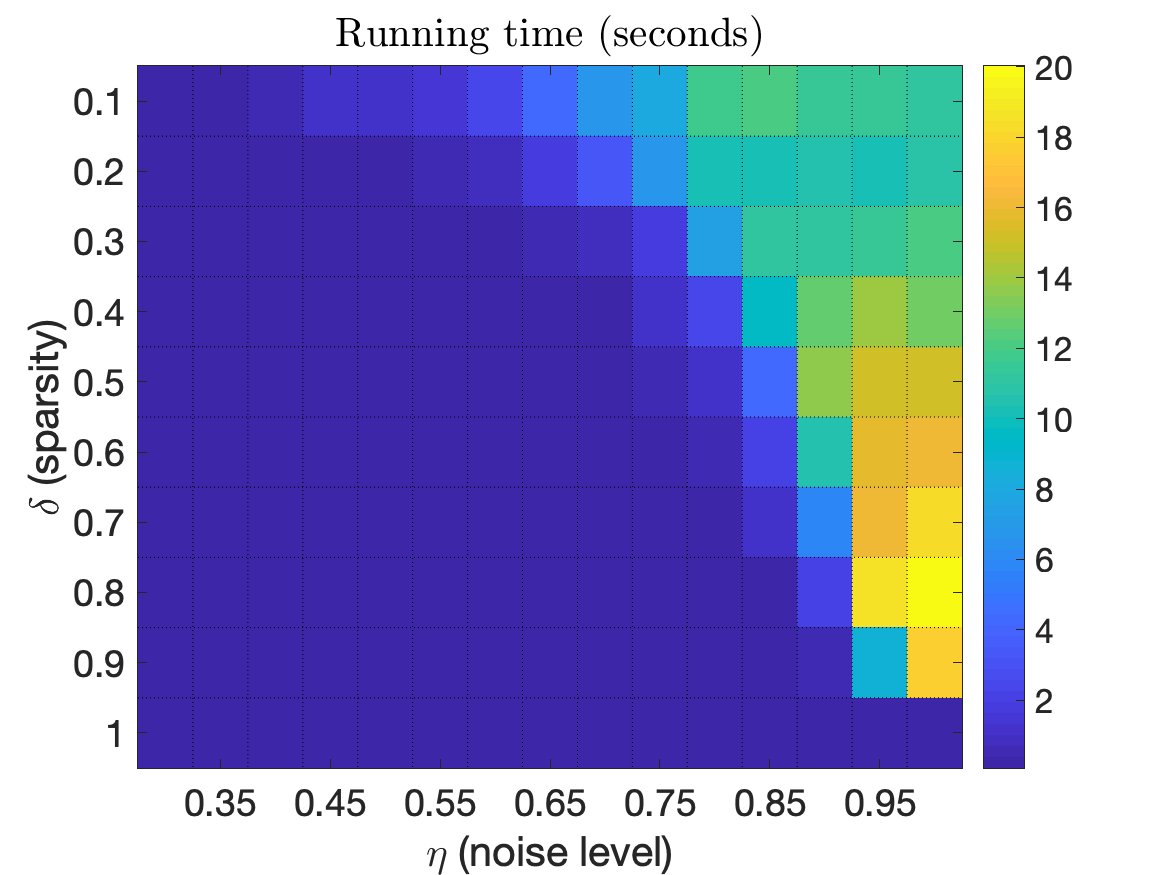}}
\hspace{3mm} 
\end{centering} 
\captionsetup{width=0.99\linewidth} 
\caption{Numerical results for MAX-CUT on a perturbed  complete bipartite graph, as we vary the noise level $\eta$ and the sampling sparsity $\delta$. Results are averaged over 20 runs.}
\label{fig:MaxCut_SBM} 
\end{figure}

For the MAX-CUT problem, we consider two sets of numerical experiments. First, we consider a version of the stochastic block model which essentially perturbs a complete bipartite graph
\begin{equation} \label{biparite} 
	\bs{B}=
\begin{vmatrix}
	\mb{0}_{n_1 \times n_1} & \mb{1}_{n_1 \times n_2} \\
	 \mb{1}_{n_2 \times n_1} & \mb{0}_{n_2 \times n_2} \\
\end{vmatrix}, 
\end{equation}
where $\mb{1}_{n_1 \times n_2}$ (respectively, $\mb{0}_{n_1 \times n_2}$)
denotes an $n_1 \times n_2$ matrix of all ones, respectively, all zeros. In our experiments, we set $n_1 = n_2 = \frac{n}{2}$, and fix $n=500$.
We perturb $\bs{B}$ by deleting edges across the two partitions, and inserting edges within each partition. More specifically, we generated the \textit{full} adjacency matrix $A^0$ from $ \bs{B} $ by adding edges independently with probability $\eta$ within each partition (i.e., along the diagonal blocks in  \eqref{biparite}). Finally, we denote by $A$ the masked version we observe, $A = A^0 \circ S $, 
where $S$ denotes the adjacency matrix of an Erd\H{o}s-R\'enyi($n$, $\delta$) graph. The graph shown in Figure \ref{fig:BipartitePerturbation} is an instance of the above generative model. 
Note that for small values of $\eta$ we expect the maximum cut to occur across the initial partition $\mathcal{P}_{\bs{B}}$ in the clean bipartite graph $\bs{B}$, which we aim to recover as we sparsify the observed graph $A$.  
The heatmap in the left of Figure \ref{fig:MaxCut_SBM} shows the Adjusted Rand Index (ARI) between the initial partition $\mathcal{P}_{\bs{B}}$ and the partition of the Max-Cut SDP relaxation in \eqref{eq:SDP_relax_max_cut}, as we vary the noise parameter $\eta$ and the sparsity $\delta$. As expected, for a fix level of noise $\eta$, we are able to recover the hypothetically optimal Max-Cut, for suitable levels of the sparsity parameter. 
The heatmap in the right of Figure \ref{fig:MaxCut_SBM} shows the computational running time, as we vary the two parameter, showing that the \textsc{Manopt} solver takes the longest to solve dense noisy problems, as one would expect.

\begin{figure}[h!]
\centering
\includegraphics[width=0.4\columnwidth]{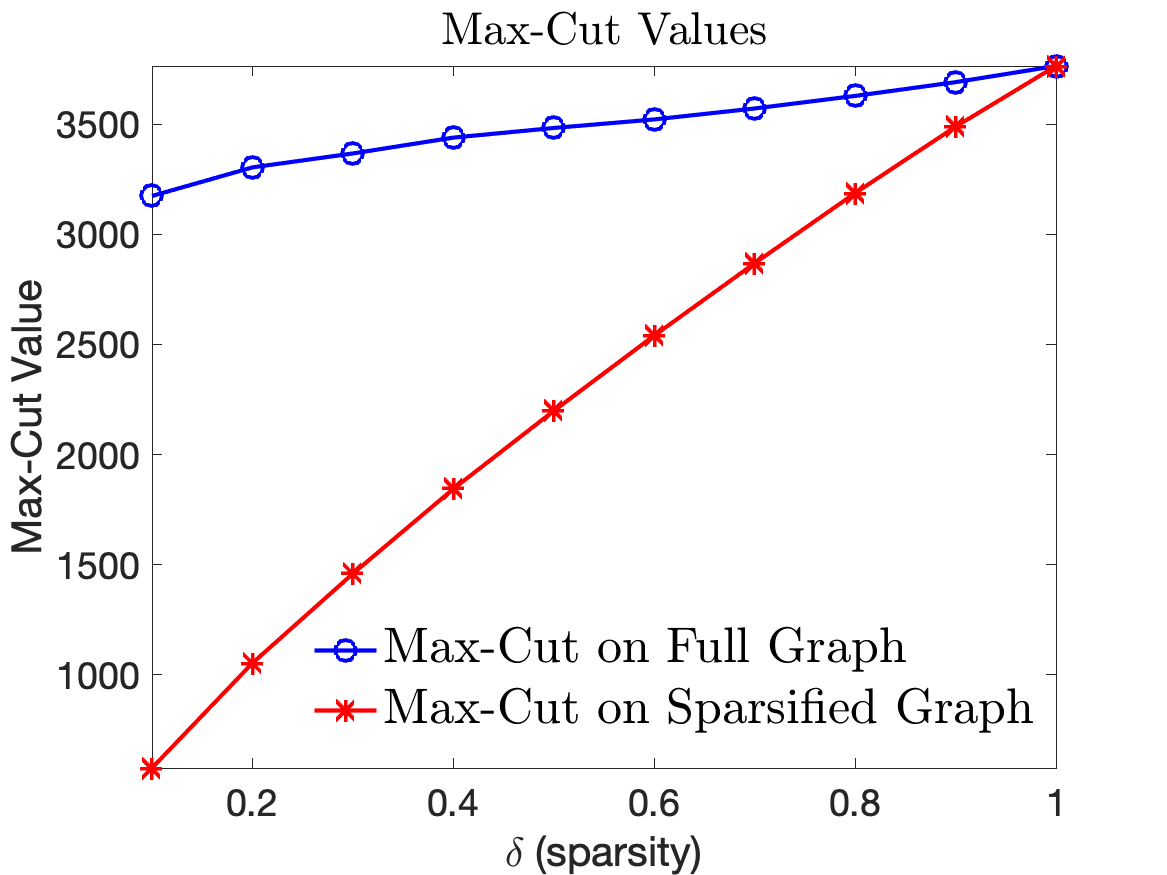}
\caption{Max-Cut results for the \textsc{G53} benchmark graph (from the \textsc{Gset} collection) with $n=1000$ nodes and average degree $\approx 12$. % 11.8
Results are averaged over 20 runs. }
\label{fig:MaxCut_G53_MaxCut}
\end{figure}

In the second set of experiments shown in Figure \ref{fig:MaxCut_G53_MaxCut}, we consider a graph $A^0$ chosen at random from the collection\footnote{\url{http://web.stanford.edu/~yyye/yyye/Gset/}} of graphs known in the literature as the \textsc{Gset}, where we vary the sparsity level $\delta$, and show the Max-Cut value attained on the original full graph $A^0$, but using the Max-Cut partition computed by the SDP relaxation \eqref{eq:SDP_relax_max_cut} on the sparsified graph $A$.

% \begin{figure}[h!]
% \captionsetup[subfigure]{skip=1pt}
% \begin{centering}\hspace{2mm}
% \subcaptionbox{The values of the Max Cuts obtained   }[0.45\columnwidth]{\includegraphics[width=0.4\columnwidth]{PLOTS/MaxCut/MaxCut_G53_n1000_Reps20.png}}
% %
% % \hspace{3mm} 
% \subcaptionbox{Running time (\textsc{MANOPT}) \rednote{perhaps drop?}. }[0.45\columnwidth]{\includegraphics[width=0.4\columnwidth]{PLOTS/MaxCut/MaxCut_G53_n1000_Reps20_time.png}}
% % \hspace{3mm} 
% \end{centering} 
% \captionsetup{width=0.99\linewidth} 
% \caption{Max-Cut results for the \textsc{G53} benchmark graph (from the \textsc{Gset} collection) with $n=1000$ nodes and average degree $\approx 12$. % 11.8
% Results are averaged over 20 runs.} 
% \label{fig:MaxCut_G53_MaxCut} 
% \end{figure} 

\subsection{Angular Synchronization} 
\begin{figure}[h!]
\captionsetup[subfigure]{skip=0pt}
\begin{centering}\hspace{2mm}
\subcaptionbox{ Spectral relaxation. }[0.45\columnwidth]{\includegraphics[width=0.4\columnwidth]{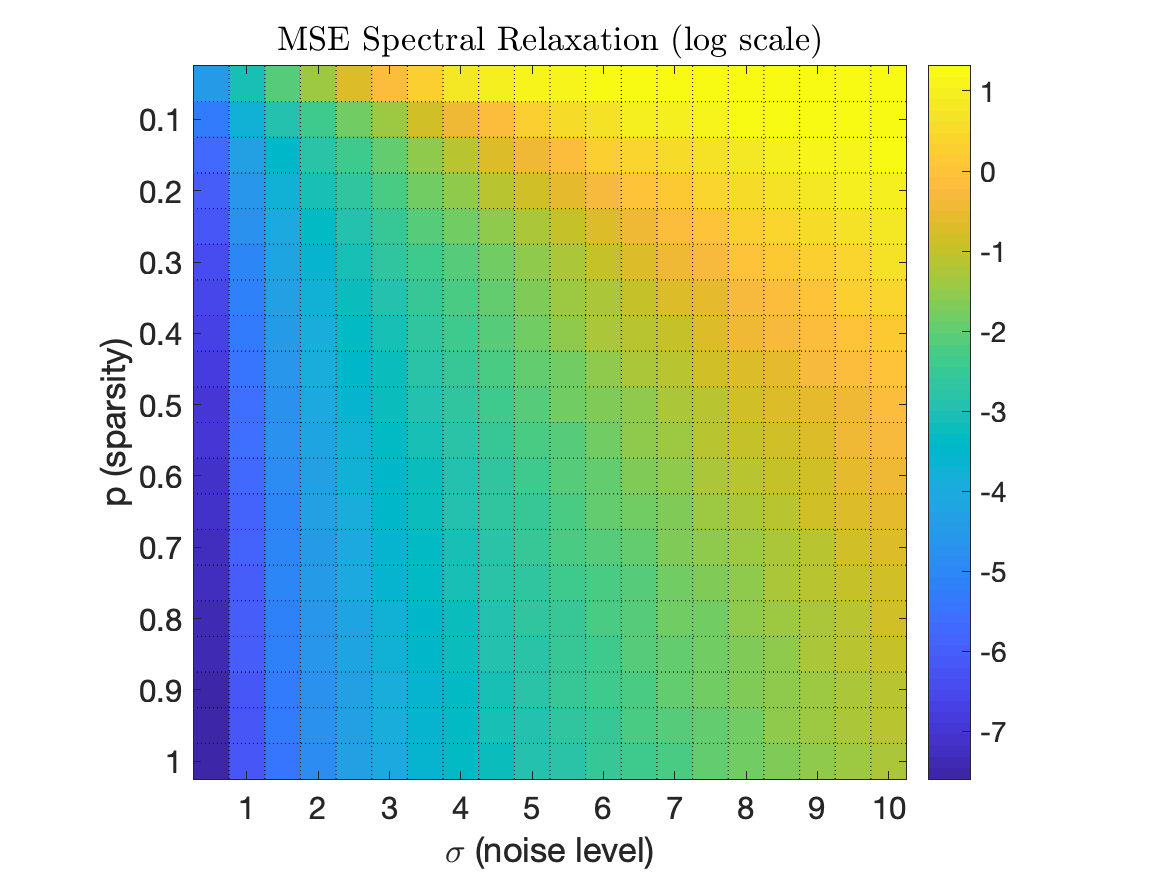}}
\hspace{3mm} 
\subcaptionbox{SDP relaxation (solved via MANOPT). }[0.45\columnwidth]{\includegraphics[width=0.4\columnwidth]{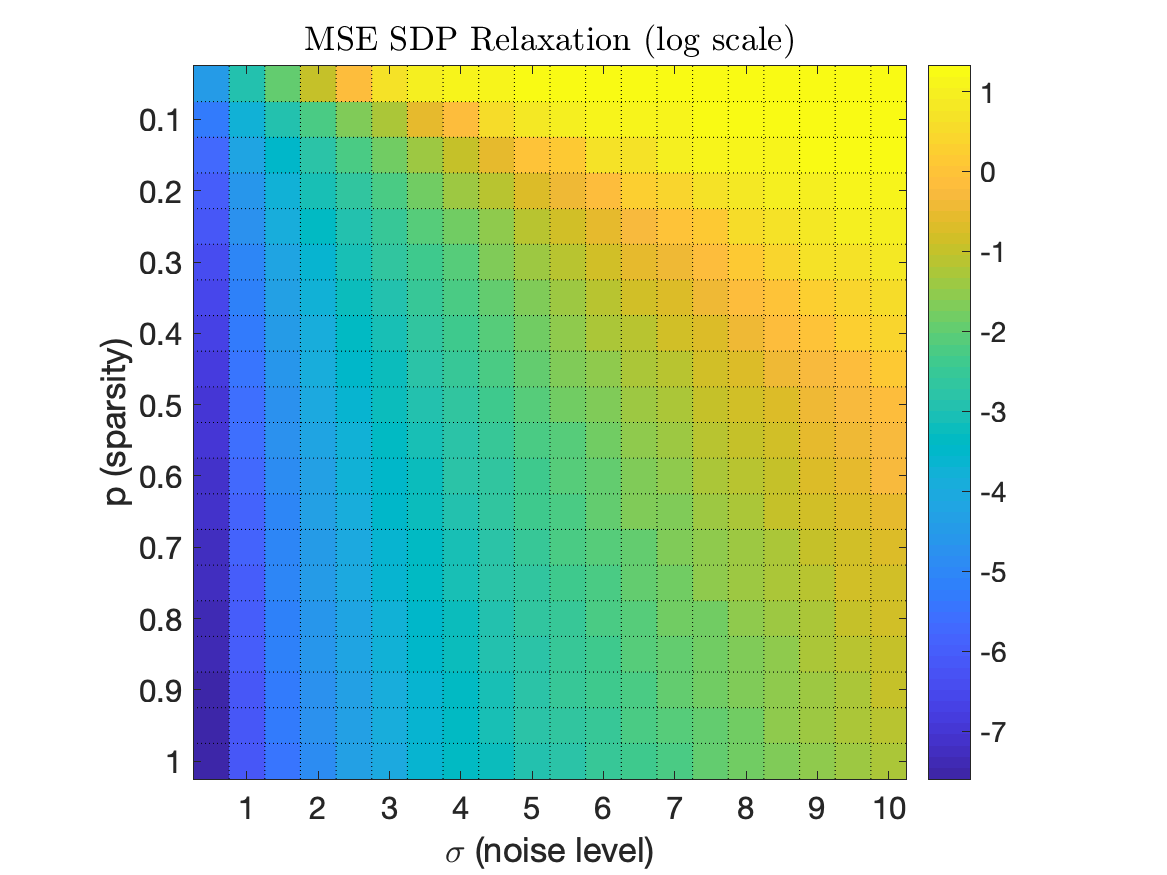}}
 \hspace{3mm} 
\end{centering} 
\captionsetup{width=0.99\linewidth} 
\caption{Recovery rates (MSE  \eqref{def:MSE_sync} - the lower the better) for angular synchronization with $n=500$, under the Gaussian noise model, as we vary the noise level $\sigma$ and the sparsity $p$ of the measurement graph. Averaged over 20 runs.} 
\label{fig:SyncGaussian} 
\end{figure}

\begin{figure}[h!]
\captionsetup[subfigure]{skip=0pt}
\begin{centering}\hspace{2mm}
\subcaptionbox{ Spectral relaxation. }[0.45\columnwidth]{\includegraphics[width=0.4\columnwidth]{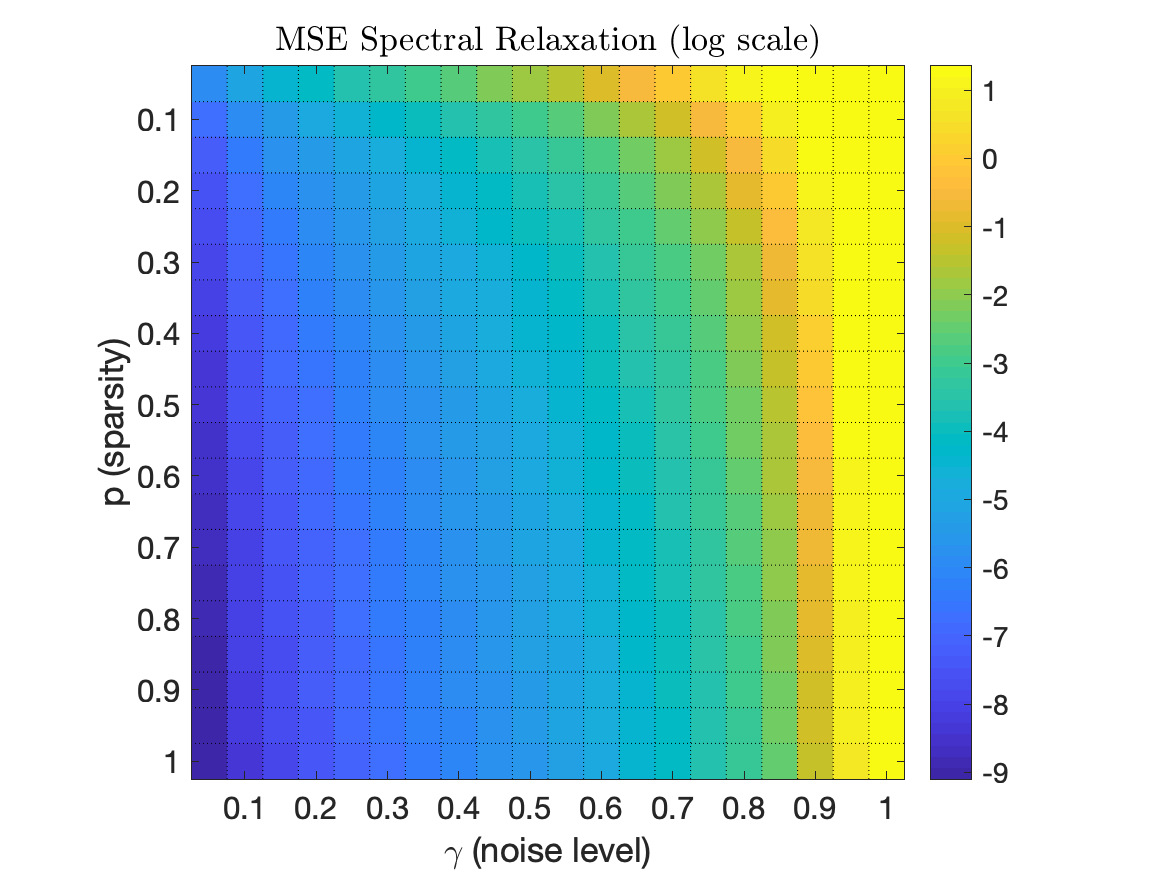}}
\hspace{3mm} 
\subcaptionbox{SDP relaxation (solved via MANOPT).}[0.45\columnwidth]{\includegraphics[width=0.4\columnwidth]{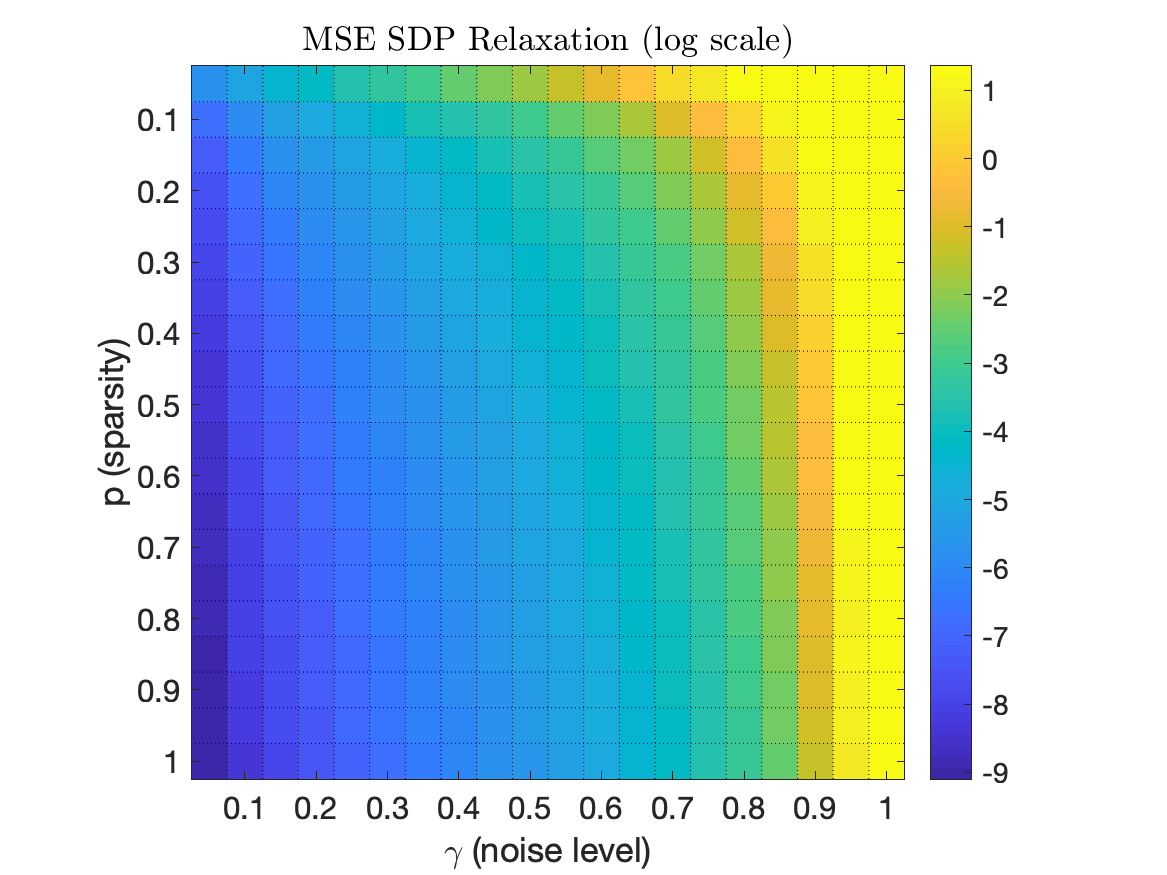}}
 \hspace{3mm} 
\end{centering} 
\captionsetup{width=0.99\linewidth} 
\caption{Recovery rates (MSE \eqref{def:MSE_sync} - the lower the better) for angular synchronization with $n=500$, under the Outlier noise model, as we vary the noise level $\gamma$ and the sparsity $p$ of the measurement graph. Averaged over 20 runs.} 
\label{fig:SyncOutlier} 
\end{figure}

For the angular synchronization problem, we consider the following experimental setup, by assessing the quality of the recovered angular solution from the SDP relaxation, as we vary the two parameters of interest. In the $x$-axis in the plots from Figures \ref{fig:SyncGaussian} and \ref{fig:SyncOutlier} we vary the noise level $\sigma$, under two different noise models, Gaussian and outliers. On the $y$-axis, we vary the sparsity of the sampling graph.

We measure the quality of the recovered angles via the Mean Squared Error (MSE), defined as follows. Since a solution can only be recovered up to a global shift, one needs an MSE error that mods out such a degree of freedom. The following MSE is also more broadly  applicable for the case when the underlying group is the orthogonal group $O(d)$, as opposed to $SO(2)$, as in the present work, where one can replace the unknown angles $\theta_1, \ldots, \theta_n$ with their respective  representation  as $ 2 \times 2$ rotation matrices  $h_1, \ldots, h_n  \in \text{O}(2) $. To that end, we look for an optimal orthogonal transformation $\hat{O} \in \text{O}(2)$ that minimizes the sum of squared distances between the estimated orthogonal transformations and the ground truth measurements 
\begin{equation}
 \hat{O} = \underset{O \in O(2)}{\operatorname{argmin}}  \sum_{i=1}^{n} \|h_i - O \hat{h}_{i} \|_{F}^{2}, 
\end{equation}
where $\hat{h}_{1}, \ldots, \hat{h}_{n} \in \text{O}(2)$ denote the $2 \times 2$ rotation matrix representation of the estimated angles $\hat{\theta}_1, \ldots, \hat{\theta}_n$.
In other words, $\hat{O}$ is the optimal solution to the alignment problem between two sets of orthogonal transformations, in the least-squares sense. Following the analysis of \cite{Singer_Shkolnisky}, and making use of properties of the trace, one arrives at 
\begin{eqnarray}
 \sum_{i=1}^{n} \|h_i - O \hat{h}_{i} \|_{F}^{2}  &=& \sum_{i=1}^{n} \; \mathrm{Trace}  \left[ \left( h_i - O \hat{h}_{i} \right) \left( h_i - O \hat{h}_{i} \right)^T \right]     \nonumber \\
							&=&  \sum_{i=1}^{n}  \;  \mathrm{Trace} \left[ 2 I - 2 O \hat{h}_i h_i^T \right] =   4 n - 2  \;  \mathrm{Trace} \left[ O  \sum_{i=1}^{n}  \hat{h}_i h_i^T \right]. 
\label{MSE_deriv}
\end{eqnarray}
If we let $Q$ denote the $2 \times 2$ matrix
\begin{equation}
  Q = \frac{1}{n} \sum_{i=1}^{n} \hat{h}_i h_i^T,
\end{equation}
it follows from (\ref{MSE_deriv}) that the MSE is given by minimizing
\begin{equation}
   \frac{1}{n}  \sum_{i=1}^{n} \|h_i - O \hat{h}_{i} \|_{F}^{2} = 4 - 2 Tr(O Q).
\end{equation}
In \cite{arun} it is proven that $Tr(O Q) \leq Tr (V U^T Q)$, for all $O \in O(3)$, where $ Q = U \Sigma V^T$ is the singular
value decomposition of $Q$. Therefore, the MSE is minimized by the orthogonal matrix $ \hat{O} = VU^T$ and is given by
\begin{equation}  \label{def:MSE_sync}
 \textsc{MSE} \;  \mydef  \;   \frac{1}{n}  \sum_{i=1}^{n} \|h_i - \hat{O} \hat{h}_{i} \|_{F}^{2} = 4 - 2 \;  \mathrm{Trace}( V U^T  U \Sigma V^T ) =
4 - 2(\sigma_1+\sigma_2), 
\end{equation}
where $\sigma_1, \sigma_2$ are the singular values of $Q$. Therefore, whenever $Q$ is an orthogonal matrix for which $\sigma_1= \sigma_2=1$, the MSE vanishes. Indeed, the numerical experiments (on a log scale) in Figures \ref{fig:SyncGaussian} and \ref{fig:SyncOutlier} confirm that for noiseless data, the MSE is very close to zero. Furthermore, as one would expected, under favorable noise regimes and sparsity levels, we have almost perfect recovery, both by the SDP and the spectral relaxations, under both noise models. 

% We could combine Figure \ref{fig:SyncGaussian} and Figure \ref{fig:SyncOutlier} into the same one row, or drop/move Figure \ref{fig:SyncOutlier} to appendix. 

% \begin{comment}
% \subsection{Community Detection}
% Decide on the experimental setup for community detection.
% \end{comment}

\section{Conclusions and future work}
%\section{Annex~3: synchronization (TODO)}
There are a number of other graph-based problems amenable to SDP relaxations, for which a similar theoretical analysis of their SDP-based estimators could be suitable. For example, the recent work of \cite{mod1JMLR} considered a problem motivated by geosciences and engineering applications of recovering a smooth unknown function $f : G \rightarrow \mathbb{R}$ (where $G = [a,b]$ is known) from noisy observations of its mod 1 values, which is also amenable to a solution  based on an SDP relaxation solved via a Burer-Monteiro approach.
Another potential application concerns the problem of clustering directed graphs, as in the very recent work of \cite{DirectedClustImbCuts} that proposed a spectral algorithm based on Hermitian matrices; this problem is also amenable to an SDP relaxation.

% Could also mention SDPs for planted clique? Standard form \cite{hajek2016achieving} and convexified MLE with nuclear norm constraint, which can be cast as an SDP as well \cite{chen2016statistical}.

Our theoretical and practical findings show that running algorithms (such as spectral methods) directly on $A$ may be improved by using first a SDP estimator such as $\hat Z$ and running the very same algorithms on $\hat Z$ (instead of $A$). Somehow, $\hat Z$ performs a pre-processing de-noising step which improve the recovery of the hidden signal such as community vectors.

\subsection*{Acknowledgements}

Mihai Cucuringu acknowledges support from the EPSRC grant EP/N510129/1 at The Alan Turing Institute. Guillaume Lecué acknowledges support from a grant of the French National Research Agency (ANR), “Investissements d’Avenir” (LabEx Ecodec/ANR-11-LABX-0047).

\begin{footnotesize}
	\bibliographystyle{plain}
	\bibliography{biblio}
\end{footnotesize}

\end{document}

%% file: packages_notes.tex
\usepackage{graphicx}
\usepackage{color}
\usepackage{amsmath}
\usepackage{amssymb}
\usepackage{amscd}
\usepackage{bbm}
\usepackage{tikz}
\usetikzlibrary{patterns}
\usepackage{hyperref}
\hypersetup{
    bookmarks=true,         % show bookmarks bar?
    colorlinks=true,       % false: boxed links; true: colored links
    linkcolor=red,          % color of internal links (change box color with linkbordercolor)
    citecolor=red,        % color of links to bibliography
    filecolor=magenta,      % color of file links
    urlcolor=cyan           % color of external links
}

\usepackage[utf8]{inputenc}
\usepackage{amsthm}
\usepackage{bbm} %indicatrice
\usepackage[linesnumbered,lined,boxed,commentsnumbered]{algorithm2e}

\newtheorem{Theorem}{Theorem}
\newtheorem{Assumption}{Assumption}
\newtheorem{Definition}{Definition}
\newtheorem{Lemma}{Lemma}
\newtheorem{Proposition}{Proposition}
\newtheorem{Corollary}{Corollary}

%% file: commandes_guillaume.tex
\newcommand{\inr}[1]{\bigl< #1 \bigr>}

\newcommand{\norm}[1]{\left\|#1\right\|}%

\newcommand\eps{\epsilon}

\def \endproof
{{\mbox{}\nolinebreak\hfill\rule{2mm}{2mm}\par\medbreak}}

\DeclareMathOperator*{\argmax}{argmax}
\DeclareMathOperator*{\argmin}{argmin}

\def\ds1{\textrm{1\kern-0.25emI}} %{\mathds{1}}

\newcommand{\blue}[1]{{\color{blue} #1}}

\newcommand \E{\mathbb{E}}
\newcommand \R{\mathbb{R}}

\newcommand \cA{{\cal A}}

\newcommand \cC{{\cal C}}

\newcommand \cN{{\cal N}}

\newcommand \cP{{\cal P}}

\newcommand \cS{{\cal S}}

\newcommand \cU{{\cal U}}

\newcommand \cX{{\cal X}}

\newcommand \bC{{\mathbb C}}

\newcommand \bE{{\mathbb E}}

\newcommand \bI{{\mathbb I}}

\newcommand \bN{{\mathbb N}}

\newcommand \bP{{\mathbb P}}

\newcommand \bR{{\mathbb R}}